\documentclass[a4paper,12pt,twoside]{article}
\usepackage[nottoc]{tocbibind} 
\usepackage[latin1]{inputenc}                 
\usepackage{amsfonts,amssymb,amsmath,exscale} 
\usepackage[dvips]{graphicx,color,psfrag}     
\usepackage{fancyhdr}                         
\usepackage[titles]{tocloft}
\usepackage{caption}                          
\usepackage{rotating}
\usepackage{defcml}                           
\usepackage[numbers]{natbib}                  
\Floatboxname{Box}
\usepackage{verbatim}                         

\usepackage{amsmath,amsthm,amsfonts,amssymb,amscd}
\usepackage{subfig}
\usepackage{graphicx,bm,color}
\usepackage{algorithm}
\usepackage{setspace}

\usepackage{algpseudocode}
\usepackage{stmaryrd}
\usepackage[framed,numbered,autolinebreaks,useliterate]{mcode}
\usepackage{xargs}                      
\usepackage[color=gray!20, backgroundcolor=yellow, textwidth=1.85cm]{todonotes}
\newcommandx{\ju}[2][1=]{\todo[linecolor=blue,backgroundcolor=blue!25,bordercolor=blue,#1]{#2}}

\usepackage{mathtools}

\usepackage{booktabs,ctable,multirow,longtable} 
\usepackage[latin1]{inputenc}                 
\usepackage{amsfonts,amssymb,amsmath,exscale} 
\usepackage[dvips]{graphicx,color,psfrag}     
\usepackage{fancyhdr}                         
\usepackage{caption}                          
\usepackage{rotating}
\usepackage{defcml}                           
\usepackage[numbers]{natbib}                  
\usepackage{url}
\Floatboxname{Box}
\newcolumntype{\xx}{\bm{x}}

\newcommand{\xx}{\boldsymbol{x}}
\newcommand{\R}{\mathbb{R}}

\definecolor{gray}{gray}{0.6}

\newtheorem{Remark}{Remark}[section]

\newtheorem{form}{Formulation}[section]

\def\R{\mathbb R}
\usepackage{empheq}

\newcommand{\noii}[1]{{\textcolor{black}{#1}}}

\fancypagestyle{plain}{%
	\fancyhead{}%
	\fancyfoot[c]{\sffamily\thepage}%
}
\makeatletter
\def\cleardoublepage{\clearpage\if@twoside \ifodd\c@page\else
	\hbox{}
	\vspace*{\fill}
	\thispagestyle{empty}
	\newpage
	\if@twocolumn\hbox{}\newpage\fi\fi\fi}

\usepackage{scalerel,stackengine}
\stackMath
\newcommand\reallywidecheck[1]{%
	\savestack{\tmpbox}{\stretchto{%
			\scaleto{%
				\scalerel*[\widthof{\ensuremath{#1}}]{\kern-.6pt\bigwedge\kern-.6pt}%
				{\rule[-\textheight/2]{1ex}{\textheight}}
			}{\textheight}%
		}{0.5ex}}%
	\stackon[1pt]{#1}{\scalebox{-1}{\tmpbox}}%
}

\begin{document}



\thispagestyle{empty}
\vspace*{3cm}
\ce{\bf\Large Phase-Field Modeling of Fracture for Ferromagnetic}
\vspace*{0.45cm}
\ce{\bf\Large  Materials through Maxwell's Equation  \\[3mm]}
\vspace*{0.45cm}
\vskip .35in

\ce{
	Nima Noii\(^{a,b,}\)\footnote{Corresponding author (Nima Noii).\\[3mm] 
		E-mail addresses: nima.noii@dikautschuk.de (N. Noii); 
		mehran.ghasabeh@ifgt.tu-freiberg.de (M. Ghasabeh); 	
    wriggers@ikm.uni-hannover.de (P. Wriggers).
}, Mehran Ghasabeh\(^{c}\), Peter Wriggers\(^{b,d}\)} \vskip .25in
\vspace{-0.2cm}
\ce{\(^a\) Deutsches Institut f\"ur Kautschuktechnologie (DIK e.V.)} \ce{Eupener Stra\ss e 33, 30519 Hannover, Germany}\vskip .25in
\vspace{-0.2cm}
\ce{\(^b\) Institute of Continuum Mechanics} \ce{Leibniz Universit\"at Hannover, An der Universit\"at 1, 30823 Garbsen, Germany}\vskip .25in
\vspace{-0.2cm}
\ce{\(^c\) Chair of Soil Mechanics and Foundation Engineering} \ce{Technische Universit\"at Bergakademie Freiberg, 09599 Freiberg, Germany} \vskip .25in
\vspace{-0.2cm}
\ce{\(^d\) Cluster of Excellence PhoenixD (Photonics, Optics, and
	Engineering - Innovation} \ce{Across Disciplines), Leibniz Universit\"at 
	Hannover, Germany}\vskip .35in

\begin{Abstract}
	Electro-active materials are classified
	as \textit{electrostrictive} and \textit{piezoelectric}
	materials. They deform under the action of an external
	electric field. \textit{Piezoelectric} material, as a special
	class of active materials, can produce an internal electric
	field when subjected to mechanical stress or strain. In
	return, there is the converse piezoelectric response, which
	expresses the induction of the mechanical deformation in the
	material when it is subjected to the application of the
	electric field. This work presents a variational-based
	computational modeling approach for failure prediction of
	ferromagnetic materials. In order to solve this problem, a
	coupling between magnetostriction and mechanics is modeled,
	then the fracture mechanism in ferromagnetic materials is
	investigated. Furthermore, the failure mechanics of
	ferromagnetic materials under the magnetostrictive effects is
	studied based on a variational phase-field model of
	fracture. Phase-field fracture is numerically challenging
	since the energy functional may admit several local minima,
	imposing the global irreversibility of the fracture field and
	dependency of regularization parameters related discretization
	size. Here, the failure behavior of a magnetoelastic solid
	body is formulated based on the Helmholtz free energy
	function, in which the strain tensor, the magnetic induction
	vector, and the crack phase-field are introduced as state
	variables. This coupled formulation leads to a continuity equation for the magnetic vector potential through well-known Maxwell's equations. Hence, the energetic crack driving force is governed by the coupled magneto-mechanical effects under the magneto-static state. Several numerical results substantiate our developments.
	
\end{Abstract} 
\vspace{-0.5cm}
{\tableofcontents} 

\sectpa[Section1]{Introduction} 

Active materials are a particular type of material that undergoes mechanical deformation in response to external effects. The sources of these effects can be pressure,
thermal, electric, and magnetic fields \cite{john+etal07}. 
In the literature, these materials are known as smart materials. 
There are various types of active materials, such as shape memory alloys, electrostrictive elastomers, piezoelectric materials, ferroelectric materials, and electro- and magneto-active polymers.
\cite{bar-cohen+zhang08,kim+tadokoro07}. The basic mechanical
properties of active materials, such as strain generation capability, stiffness, strain,
Hysteresis and electrical impedance vary widely. The energy output
of active materials used in actuators is principally dominated by stiffness and the amount of strain energy generated by the material \cite{john+etal07}. The materials can be classified according to their response when they are subjected to external stimuli. To this end, active materials can be generally categorized as electro-active and magneto-active materials. 

Electro-active materials are classified as \textit{electrostrictive}, and \textit{piezoelectric materials}. These materials deform mechanically under an external electric field. Piezoelectric materials can produce an internal electric field when subjected to mechanical stress or strain (direct piezoelectricity) \cite{berlincourt71}. Piezoelectric materials also represent the reverse piezoelectric effect (converse piezoelectricity), where stress or strain is generated when they are subjected to an electric
field \cite{bhansali+vasudev12}. The direct piezoelectric response of material is determined by the conversion of the mechanical into the electrical energy; however, the converse piezoelectric response
expresses the induction of mechanical deformation as a consequence of the applied electric field \cite{damjanovic+newham92}. 
\begin{figure}[!t]
	\centering
	{\includegraphics[clip,trim=0cm 2cm 3cm 1cm, width=13cm]{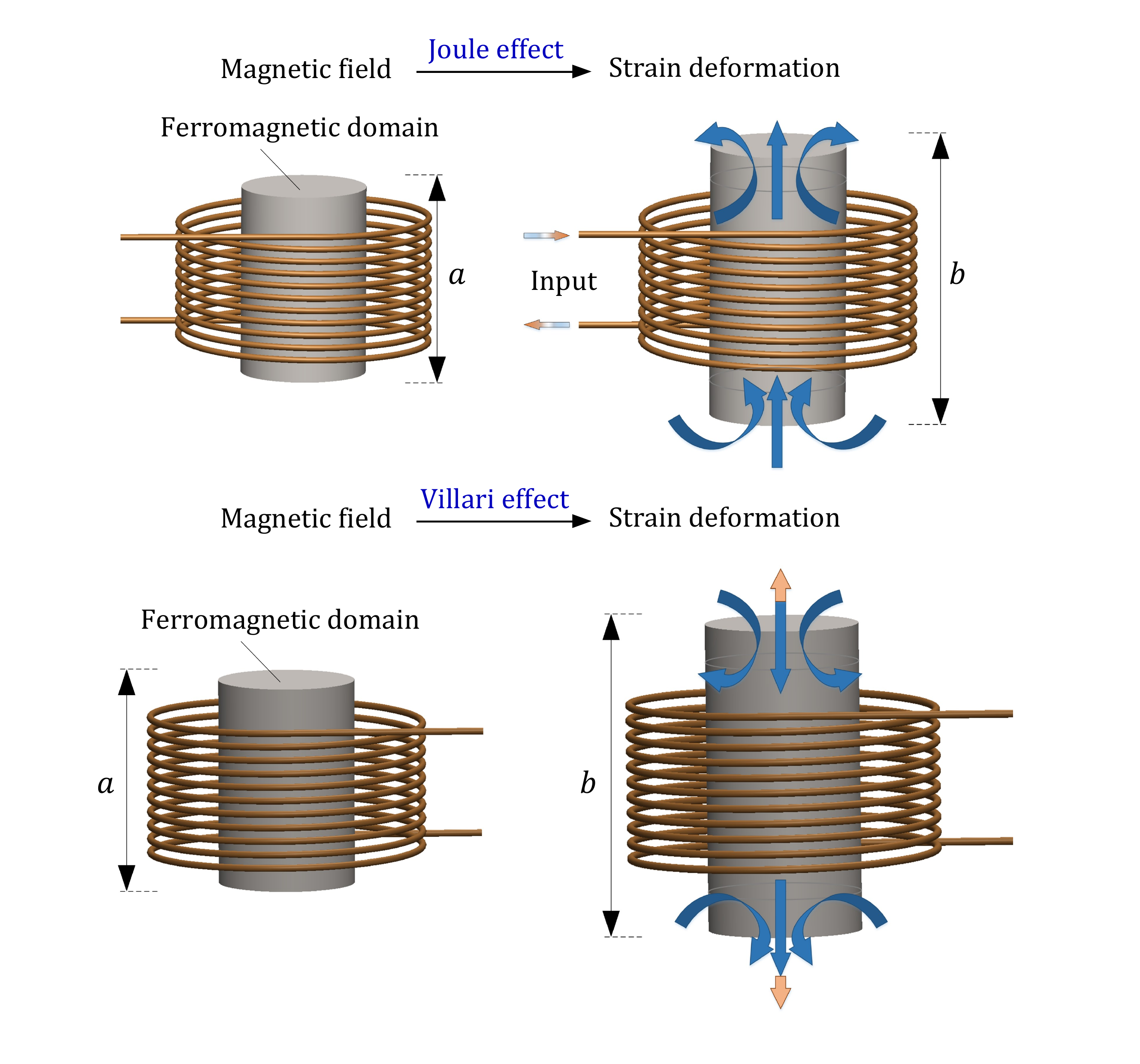}}  
	\caption{Schematic representation of Joule effect (direct
	magnetostrictive effect) in the first raw and Villari effect
	(reversible magnetostrictive effect) in the second raw.}
	\label{Figure:sc:into}
\end{figure}
Piezoelectric materials are applied within a great
range, including piezoelectric motors, actuators in industrial
sector, sensors in the medical sector, actuators in consumer electronics
(printer, speakers), piezoelectric buzzers, piezoelectric igniters, microphones, nanopositioning in
atomic force microscope (AFM) and the scanning tunneling microscope
(STM), and micro-robotics (defense) \cite{sekhar+etal21}.

Electrostrictive materials generate stress or
strain when subjected to an external electric field
\cite{john+etal07}. However, there is a principal difference between 
piezoelectric and electrostrictive materials. The piezoelectric effect is possible only in
non-centrosymmetric materials. However, the electrostrictive
effect is not limited by symmetry and is present in all materials, even those that are amorphous. Therefore, the electrostrictive effect
exhibits a nonlinear (second-order or quadratic)
dependency of the strain on the applied electric
field \cite{damjanovic+newham92}. 

Magnetostriction, a key feature of magnetic materials, is defined
as the alteration in shape and dimension of the material during
the magnetization process. It continues till the magnetic saturation of
the material is attained. The magnetostriction phenomenon was first
introduced in the work of Joule \cite{joule1842}. It is also known as
the Joule effect. In 1846, Villari discovered the
reverse effect of magnetization (see Figure
\ref{Figure:sc:into})\cite{villari1864}. Ferromagnetic material and ferri-magnetic
are known for their magnetostrictive characteristics.
A magnetostrictive material consists of tiny fragments. These
fragments, usually iron, nickel, or cobalt, have small magnetic moments
as a result of their "3d" shells that are not completely filled with
electrons. Ferromagnetics fundamentally acts like tiny permanent
bar magnets. Magnetic materials have a property of magnetostriction
that causes them to change their shape or dimensions during the
magnetization process. Such materials can convert magnetic energy into
kinetic or reverse and thus are mainly used to construct
actuators and sensors. The applicability of the materials is
quantified by the magnetostrictive coefficient. It is defined as the
fractional change in length as the magnetization of the material
increases from zero to the saturation value. This coefficient
can be positive and negative. The magnetostriction
characteristic induces energy loss due to
frictional heating in susceptible ferromagnetic cores. It also
produces the low-pitched humming sound
that can be heard coming from transformers, where an existing
magnetic field is changed by oscillating
AC currents.

In the literature, computational approaches to modeling the electromagnetic response of the electro- and magneto-active materials are developed within a framework of thermodynamics. For some active-material, it is sufficient to take into account only a quasi-static state for modeling the magnetic response \cite{birk+reichel+schroeder22}. The work proposes a hybrid SBFEM-FEM approach for efficiently calculating stray magnetic fields in unbounded
domains. Most of the works aim at constructing a coupling between electromagnetics and mechanics. The key feature of these approaches is to investigate the stress response of the materials under the electro- and magneto-striction effects. The coupled formulation within the framework of magnetoelasticity is found in the work of Brigadnov and Dorfmann \cite{brigadnov+dorfmann03} where a mathematical model is developed to investigate the mechanical characteristics of the magneto-sensitive elastomer under the application of a magnetic field. Further studies related to magneto-sensitive elastomers are conducted in the work of Dorfmann, and his colleagues \cite{dorfmann+ogden03, dorfmann+ogden+saccomandi04}. The investigation of an electromagnetic forming process is conducted in \cite{thomas+triantafyllidis09} by formulating a fully coupled electromagnetic-thermomechanical model. A computational approach incorporating the energy-based magneto-mechanical model is developed in \cite{belahcen+etal06, fonteyn10, fonteyn2010fem} for the electric electrical steel sheets. 
An energy-density functional method is proposed in \cite{rasilo+etal19} based on an isotropic spline-based thermodynamic approach to model 
magneto-mechanical behavior in ferromagnetic material.

Electro-magneto-mechanical coupling effects are formulated, based on
incremental variational principles, within a general continuum
mechanics framework and numerical
implementation of dissipative functional materials in \cite{miehe+rosato+kiefer11}. Micromagnetics is another concept where a geometrically consistent incremental variational formulation is extended in a  micro-magneto-elastic model that accounts for the micro-structural evolution of both magnetically- and mechanically-driven magnetic domains in ferromagnetic materials \cite{miehe2012geometrically}.  
The variational-based computational modeling of ferromagnetics and magnetorheological elastomers, which present a mutually-coupled magneto-mechanical response, is proposed in \cite{ethiraj14}. Further study on the behavior of magnetosensitive polymers based on a multiplicative magneto-elastic model within a micromechanically-based framework is found in
\cite{ethiraj+miehe16}. Moreover, models regarding the
thermomechanical effects are developed to investigate the coupling
between electromagnetics and thermomechanics
in \cite{hanappier+charkaluk+triantafyllidis21} for electric
motors. The multi-physics problem regarding
magnetic-thermal-mechanical modeling of the electromagnetic rail
launching is investigated in \cite{zhang+etal21}. The principles of
designing an in situ fatigue testing system, including structural resonances, transient response, grip design, and thermal insulation
performance of the device, are studied in detail in \cite{ma+etal16} by developing a coupled thermo-mechanical model. The Cell-based smoothed finite element method (CS-FEM) and the asymptotic homogenization method (AHM) is incorporated in \cite{zhou+etal21} to simulate the magneto-electro-elastic response of a structure under dynamic load accurately.

In engineering applications, predicting crack initiation and
propagation in structures under mechanical loading and environmental
conditions is greatly important.
Besides classical approaches for fracture mechanics, non-local damage models and variational
approaches have been
developed in the literature. The non-local damage models have been formulated to
overcome the ill-posedness due to spatial discretization in finite
element simulations. The variational approaches are introduced based
on energy minimization \cite{francfort+marigo98, bourdin+francfort+marigo08, dal+toador02}, and their regularization is obtained by
$\Gamma$-convergence, which is fundamentally inspired by the work of image segmentation conducted by Mumford and Shah
\cite{mumford+shah89}. The model is then improved by formulating a Ginzburg-Landau-type evolution
equation of the fracture phase-field \cite{hakim+karma09}. In recent years, the variationally-based
phase-field approach to fracture within the thermodynamically consistent framework is proposed in the
work of Miehe et al.  \cite{miehe+hofacker+welschinger10}.
In the latest work, the robust algorithmic formulation of the evolution of diffuse crack topology in time is proposed by introducing a local history field that determines the maximum tensile strain obtained in history. It then acts as a crack driving force in the evolution of the crack phase field. In the literature, few works investigate the material's cracking response; for example, in \cite{miehe+welschinger+hofacker101}, the failure mechanism of piezoelectric is studied by extending electromechanics to the phase-field approach. Ferroelastic ceramics is another class of active material in which the crack initiation and propagation are examined in \cite{linder+miehe12}, which proposes a computational framework that regards the electric displacement saturation and its effect on the hysteretic behavior of ferroelectric ceramics and the resulting cracking response.

This work applies the phase-field approach to the magneto-mechanically driven fracture in ferromagnetic materials. In the formulations for examining the response of the electromagnetic materials, finite element method (FEM) is applied to discretize
the time-dependent Maxwell's equations on a bounded domain in
three-dimensional space. At first, the formulation is developed by
deriving a weak formulation for the electric and magnetic fields with
approximate Neumann and Dirichlet boundary conditions, and the problem
is discretized both in time and space. In general, the electric and
the magnetic fields are discretized by adapting \textit{N\'ed\'elec}
curl-conforming and \textit{Raviart-Thomas} div-conforming finite element
approaches, respectively.   

The main objective of this contribution is to introduce:

\begin{itemize}
	\item Extension of Maxwell's equations towards a coupled
	magneto-mechanical model to investigate the stress response of the
	ferromagnetic materials;
	\item Transition rule for the electromagnetic material properties from undamaged to fully damaged states through the fracture phase-field, which acts as a geometric interpolation variable;
	\item Investigation of the magnetostrictive-induced cracking
	in the ferromagnetic materials by developing a magneto-mechanical
	model coupled with the phase-field model;
	\item Representation of numerical examples to substantiate our developments in predicting the fracture response
	of the ferromagnetic material.
\end{itemize}

The paper is structured as follows: For a better insight into
phase-field modeling of fracture for ferromagnetic materials, primary
fields for the multi-field problem are provided in Section 2. Through
a regularized fracture phase-field formulation, we further
outline the theoretical framework for magnetostrictive-induced fracture in ferromagnetic materials. The transition rule for the electromagnetic material properties from undamaged to fully damaged states is also discussed. Next, Section 3 outlines the constitutive energy density for the magneto-mechanical formulation of fracturing solids, which suffices to define the variational formulation setting. In Section 4, three numerical simulations are performed to demonstrate the correctness of our algorithmic developments. Finally, the paper is concluded with some remarks.

\sectpa[Section2]{Phase-field modeling of fracture for ferromagnetic materials}
This section outlines a mathematical model for magnetostrictive-induced cracking in ferromagnetic materials by developing a magneto-mechanical model, considering small deformations. The fracture process is modeled by employing the well-developed phase-field formulation to resolve the sharp crack surface topology in the regularized concept. We specifically elaborate \textit{three governing equations} to characterize the constitutive formulations for the mechanical deformation, electromagnetic as well as the fracture phase-field. This coupled formulation is derived by obtaining a continuity equation for the magnetic vector potential through well-known \textit{Maxwell's equations}. Thus, the comprehensive objective of the following section is to advance a continuum theory of cracking response in ferromagnetic materials within a framework of thermodynamics.

\sectpb[Section21]{Primary fields for the multi-field problem}
\label{sec:basic_def}

\begin{figure}[!t]
	\centering
	{\includegraphics[clip,trim=2cm 10.3cm 1cm 18cm, width=16cm]{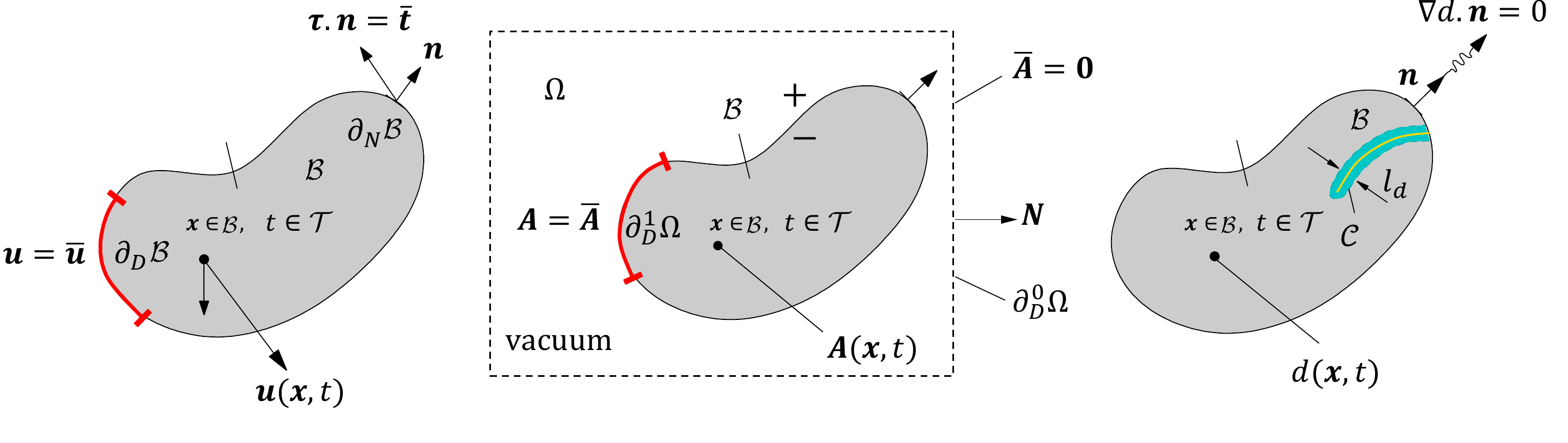}}  
	\vspace{-1cm}
	\caption*{\hspace*{0.1cm}(a)\hspace*{4cm}{(b)}\hspace*{5cm}{(c)}}
	\caption{Primary variable fields in magneto-elasticity, for a solid
		body $\calB \subset \calR^{\delta}$ with dimension $\delta \in [2,3]$. $(a)$
		The displacement field $\Bu$ defined on $\calB$ and Neumann-type
		boundary condition for traction $\bar{\Bt}=\Btau \cdot
		\Bn$ on $\partial_N\calB$. $(b)$ The magnetic vector potential $\BA$ in $\Omega$ that is
		continuous across the interface $\partial{\calB}$,
		i.e. $\llbracket \BA \rrbracket:=\BA^{+}-\BA^{-}=\textbf{0}$ on
		$\partial{\calB}$, and zero on the boundary $\partial{\Omega}$ free
		space box (vacuum). $(c)$ The crack phase-field is determined by
		Dirichlet-type boundary condition $d=1$ on $\mathcal{C}$ and Neumann-type
		boundary condition $\nabla d \cdot \Bn = 0$ on the full surface $\partial{\calB}$.}
	\label{Figure1}
\end{figure}

Let $\Omega\subset{\mathbb{R}}^{\delta}$ be an arbitrary vacuum free space box with dimension $\delta=\{2,3\}$ and $\calB\subset{\mathbb{R}}^{\delta}$ denote the domain occupied by the material solid,
as demonstrated in Figure \ref{Figure1}. Vacuum space $\Omega$ is considered to be large enough such that
the magnetic field induced by the magnetization of the body $ \calB$ is
decayed at its surface $\partial \calB\subset{\mathbb{R}}^{\delta}$. Also, in the following, let a Dirichlet boundaries conditions of vacuum defined as $\partial_D\Omega=\partial^0_D\Omega\cup\partial^1_D\Omega$ which contains inner surface (it has an intersect with the solid body $ \calB$) and outer boundary surface, respectively, as shown in Figure \ref{Figure1}b . We further assume Dirichlet boundary conditions on $\partial_D\calB$ and complete Neumann boundary conditions on $\Gamma_N := \partial_N \calB \cup \mathcal{C}$, where $\partial_N \calB $  denotes the outer domain boundary (where the traction imposed) and {$\calC\in \mathbb{R}^{\delta-1}$} is the crack boundary, as illustrated in Figure \ref{Figure1}c. Therefore, we have three different domains required for the electromagnetically induced fracture to be used as follows:
\begin{equation*}
\texttt{Vacuum region:} \;\Omega\subset{\mathbb{R}}^{\delta},   \;\;
\texttt{Solid region:}  \;\;  \calB\subset{\mathbb{R}}^{\delta} \AND
\texttt{Damage region:}  \;\; \calC\in \mathbb{R}^{\delta-1}.
\end{equation*}

The response of the fracturing solid at material points $\Bx\in\calB$ at time $t\in \calT= [0,T]$ is described by the displacement field $\Bu(\Bx,t)$ and the crack phase-field $d(\Bx,t)$ as
\begin{equation}
\Bu: 
\left\{
\begin{array}{ll}
\calB \times \calT \rightarrow \mathbb{R}^\delta \\[2mm]
(\Bx, t)  \mapsto \Bu(\Bx,t)
\end{array}
\right.
\AND
d: 
\left\{
\begin{array}{ll}
\calB \times \calT \rightarrow [0,1] \\[2mm]
(\Bx, t)  \mapsto d(\Bx,t)
\end{array}
\right.\;.
\end{equation}
Intact and fully fractured states of the material are characterized by $d(\Bx,t)=0$ and $d(\Bx,t)=1$, respectively. In order to derive the variational formulation, the following space is first defined. For an arbitrary $A\subset\mathbb{R}^\delta$, we set
\begin{align}
\mathrm{H}^1(\calB,A):=\Big\{v:\calB\times\calT\rightarrow A \; \colon\;\; v\in \mathrm{H}^1(A)\Big\}.
\end{align}
We also denote the vector valued space $\mathbf{H}^1(\calB,A):=\left[\mathrm{H}^1(\calB,A)\right]^\delta$ and define
\begin{equation}
\calW^{\Bu}(\text{div};\calB):=\Big\{\Bu\in\mathbf{H}^1(\calB,\mathbb{R}^\delta) \;\ \colon\; \ \Bu=\overline{\Bu} \ \text{on} \ \partial_D\calB\Big\},
\label{eq:spaces_u}
\end{equation}
and correspondingly for displacement test function, we have
\begin{equation}
	\calW_0^{\Bu}(\text{div};\calB):=\Big\{\Bu\in\mathbf{H}^1(\calB,\mathbb{R}^\delta)\}.
	\label{eq:spaces_u0}
\end{equation}
Concerning the crack phase-field, we set
\begin{equation}
\calW^{d}:=\mathrm{H}^1(\calB) \AND  \calW^{d}_{d^n}(\text{div};\calB):=\Big\{d\in \mathrm{H}^1(\calB,{\color{black}[0,1]}) \; \colon\; \ d \geq d^n \Big\},
\label{eq:spaces_d}
\end{equation}
where $d^n$ is the damage value in a previous time instant. Note that $\calW^{d}_{d^n}(\text{div};\calB)$ is a non-empty, closed, and convex subset of $\calW^{d}$ and introduces the evolutionary character of the phase-field, incorporating an irreversibility condition in incremental form. Correspondingly the phase-field test function reads
\begin{equation}
	\calW^{d}_{0}(\text{div};\calB):=\Big\{d\in \mathrm{H}^1(\calB,{\color{black}[0,1]}) \; \colon\; \ d \geq 0 \Big\}.
	\label{eq:spaces_d0}
\end{equation}
%
In order to formulate a wide variety of the failure mechanism of the electro- and magneto-active materials into the variational equations, two different sets of fields need to be introduced. To this end, the electric and magnetic primary fields are introduced:
\begin{equation}\label{eq:EB}
\BE:
\left\{
\begin{array}{ll}
\Omega \times \calT \rightarrow \mathbb{R}^\delta \\
(\Bx, t)  \mapsto \BE(\Bx,t)
\end{array}
\right.
\AND
\BH:
\left\{
\begin{array}{ll}
\Omega \times \calT \rightarrow \mathbb{R}^\delta \\
(\Bx, t)  \mapsto \BH(\Bx,t)
\end{array}
\right.\;,
\end{equation}
such that, the intensity fields $(\BE,\BH)$ belong to the space:
\begin{equation}
\calS^{\Bv}(\text{curl};\Omega):=\Big\{\Bv\in\mathbf{L}^2(\Omega,\mathbb{R}^\delta) \; \colon\; \nabla \times \Bv \in\mathbf{L}^2(\Omega,\mathbb{R}^\delta) \Big\},
\label{eq:spaces_E}
\end{equation}
along with their test function space as
\begin{equation}
\calS_0^{\Bv}(\text{curl};\Omega):=\Big\{\Bv\in	\calS^{\Bv}(\text{curl};\Omega) \; \colon\; \Bv \times \Bn =0 \in\mathbf{L}^2(\Omega,\mathbb{R}^\delta) \Big\}.
\label{eq:spaces_E0}
\end{equation}
The generalized fluxes conjugate to \req{eq:EB} for the electric and magnetic primary fields read:
\begin{equation}
\BD:
\left\{
\begin{array}{ll}
\Omega \times \calT \rightarrow \mathbb{R}^\delta \\
(\Bx, t)  \mapsto \BH(\Bx,t)
\end{array}
\right.
\AND
\BB:
\left\{
\begin{array}{ll}
\Omega \times \calT \rightarrow \mathbb{R}^\delta \\
(\Bx, t)  \mapsto \BB(\Bx,t)
\end{array}
\right.\;.
\end{equation}
Correspondingly, the flux fields $(\BD,\BB)$ are defined in the space:
%
\begin{equation}
	\calS^{\Bv}(\text{div};\Omega):=\Big\{\Bv\in\mathbf{L}^2(\Omega,\mathbb{R}^\delta) \; \colon\; \nabla.\Bv \in\mathbf{L}^2(\Omega,\mathbb{R}^\delta) \AND \Bv=\overline{\Bv} \ \text{on} \ \partial^0_D\Omega\Big\},
	\label{eq:spaces_B}
\end{equation}
and augmented with their test function space, as
\begin{equation}
	\calS_0^{\Bv}(\text{div};\Omega):=\Big\{\Bv\in\calS^{\Bv}(\text{div};\Omega) \; \colon\; \Bv.\BN=0\; \;\text{on}\; \ \partial_D\Omega:=\partial^0_D \Omega\cup\partial^1_D\Omega \Big\}. 
	\label{eq:spaces_B0}
\end{equation}
Here, $\BN$ is the outward unit normal vector on the surface $\partial \Omega$, as illustrated in Figure \ref{Figure1}. 

These electric and magnetic primary fields and their fluxes describe the fundamental macroscopic fields governing all electromagnetic phenomena denoted as \textit{Maxwell's equations}. These are described in detail in Section 2.2.2.
To complete this section, we need to introduce the vector potential formulation of Maxwell's equations denoted as $\BA(\Bx,t)$. Since,  the $div$ and $curl$ operator are surjective, following \cite{monk2003finite, gross2004electromagnetic} the de Rham complex reads
\begin{equation}
	\mathbb{R} \xrightarrow{id} {H^{1}(\Omega)}\xrightarrow{^\text{grad}}  H(\text{curl},\Omega) \xrightarrow{^\text{curl}}  {H(\text{div},\Omega)} \xrightarrow{^\text{div}}  {L_2(\Omega)} \xrightarrow{^\text{0}} 0.
\end{equation}
For a magnetic primary field $\BB\in\calS^{\Bv}(\text{div};\Omega)$ it holds
\begin{equation}\label{eq:B_curl}
	\nabla.\BB=0\quad\xrightarrow{\text{de Rham}} \quad \exists \;\BA \in\calS^{\Bv}(\text{curl};\Omega) \;\;\colon
	\nabla\times\BA=:\BB,
\end{equation}
in which, $\BA$ denotes the magnetic vector potential. To summarize, we have following sets:
\begin{equation}
	\begin{aligned}
		(\Bu,d)\in\Big(\calW^{\Bu}(\text{div};\calB),\;&\calW^{d}(\text{div};\calB)\Big),\qquad (\BE,\BH,\BA) \in\calS^{\Bv}(\text{curl};\Omega)\\[2ex]
		&(\BD,\BB)\in\calS^{\Bv}(\text{div};\Omega),
	\end{aligned}
\end{equation}
and correspondingly it follows for their test functions:
\begin{equation}
	\begin{aligned}
		(\delta\Bu, \delta d)\in\Big(\calW_0^{\Bu}(\text{div};\calB), \;&\calW_0^{d}(\text{div};\calB)\Big), \qquad (\delta\BE,\delta\BH,\delta\BA) \in\calS_0^{\Bv}(\text{curl};\Omega),\\[2ex] &(\delta\BD,\delta\BB)\in\calS_0^{\Bv}(\text{div};\Omega).
	\end{aligned}
\end{equation}

\sectpb[Section22]{Governing equations of the failure mechanism for magnetostrictive effects}
\label{sec:energy_functions}
In this section, we present a theoretical development associated with the failure mechanism of the ferromagnetic material based on the variational phase-field model by taking into account their fully coupled magneto-mechanical characteristics. This includes mechanical, electromagnetic, and fracture contributions to the model. Additionally, the transition rule from undamaged to fully damaged states is presented. This section is used as a departure point for Section 3 in which to formulate
the variational framework for magnetostrictive-induced cracking is formulated.
 
\sectpc{Elastic contribution}
The standard elastic energy density, so-called the effective strain energy density \cite{noii2021quasi,noii2020adaptive} is expressed in terms of the scalar-valued function $\psi_e$. For isotropic materials, $\psi_e$ is defined in terms of the strain tensor $\Bve$. In our formulation, in order to preclude fracture in compression, a decomposition of the effective strain energy density into \textit{damageable} and \textit{undamageable} parts are employed. Thus, we perform
additive decomposition of the strain tensor into \textit{volume-changing}
(volumetric) and \textit{volume-preserving} (deviatoric) counterparts:
\begin{equation}\label{eq:additive_strain}
\bm\varepsilon(\bm u)=\bm\varepsilon^{vol}(\bm u)+\bm\varepsilon^{dev}(\bm
u) \WITH \bm\varepsilon(\bm u):=\frac{1}{2}(\nabla{\Bu}+\nabla^{T}{\Bu}),
\end{equation}
where the volumetric strain is denoted as $\bm\varepsilon^{vol}(\bm
u):=\frac{1}{3}(\bm \varepsilon(\bm u):\BI)\BI$ and the deviatoric strain is
denoted as $\bm\varepsilon^{dev}(\bm u):=\mathbb{P}:\bm \varepsilon$. The
fourth-order projection tensor
$\mathbb{P}:=\mathbb{I}-\frac{1}{3}\text\BI\otimes\BI$ 
is introduced to map the full strain tensor onto its deviatoric component.
Therein,
$\mathbb{I}_{ijkl}:=\frac{1}{2}\big(\delta_{ik}\delta_{jl}+\delta_{il}\delta_{jk}\big)$
is the fourth-order symmetric identity tensor. 

Next we construct the mechanical BVP. The solid geometry $\calB$ is
loaded by prescribed deformations $\overline{\Bu}$ and an external traction vector $\overline{\bm{t}}$ on the boundary, defined by time-dependent Dirichlet conditions and Neumann conditions as
	\begin{equation}\label{BC_M}
		\Bu  = \overline{\Bu}(\Bx,t) \ \textrm{on}\ \partial^{\Bu}_D\calB
		\AND
		\Btau \cdot \Bn 
		= \overline{\Bt}(\Bx,t)  \ \textrm{on}\ \partial^{\Bu}_N\calB ,
	\end{equation}
	where $\Bn$ is the outward unit normal vector on the surface $\partial_N \calB$. The stress tensor $\Btau$ is the thermodynamic dual to the strain $\Bve$.

The global mechanical form of the equilibrium equation for the solid body can be represented by a second-order PDE for the multi-field system as
\begin{equation}
	\fterm{
		\Div\,\Btau(\Bu, \BB, d) + \overline\Bb = \textbf{0} \, ,
		\label{equil:defo}
	}
\end{equation}
which is valid for a quasi-static response and where denote $\overline\Bb$ is the prescribed body force.

\sectpc{Electromagnetic contribution} 
%
Maxwell's equations can principally express the macroscopic electromagnetic phenomena in a ferromagnetic material. These equations are related to the calculations of electromagnetic fields, including the electric and the magnetic fields with eddy currents (their corresponding fluxes)
\begin{equation} \label{eq1}
	\begin{split}
		\text{Maxwell-Faraday law}:& &\nabla \times \BE &= -\frac{\partial{\BB}}{\partial{t}} \hspace*{1.5cm} &\text{in} \;\;\Omega \times(0,\text{T}), \\[0.25cm]
		\text{Maxwell-Amp\`ere law}:& &\nabla \times \BH &= \frac{\partial{\BD}}{\partial{t}} + \BJ  \hspace*{1.5cm} &\text{in} \;\;\Omega \times(0,\text{T}), \\[0.25cm]
		\text{Maxwell-Gauss law}:& &\nabla \cdot \BD &= \rho_0(d)  \hspace*{1.5cm} &\text{in} \;\;\Omega\times(0,\text{T}), \\[0.25cm]
		\text{No magnetic monopole law}:& &\nabla \cdot \BB &= 0  \hspace*{1.5cm} &\text{in} \;\;\Omega\times(0,\text{T}),
	\end{split}
\end{equation}

where $\BE(\Bx,t)$ and $\BH(\Bx,t)$ are the electric and magnetic fields, $\BD(\Bx,t)$ and
$\BB(\Bx,t)$ are the corresponding electric and magnetic flux densities, respectively, along with $\BJ(\Bx,t)$ that is the
electric current density, and $\rho$ is the electric charge density. In
addition to these equations, there are constitutive equations that
describe the macroscopic characteristics of a electromagnetic
materials. These equations construct relations between the magnetic and electric flux
densities and the magnetic and the electric fields
\begin{equation}\label{eq2}
	\BD =\epsilon_0(d)  \BE, \quad \quad \BB =  {\mu_0}(d)\BH \quad \AND \quad \BJ
	= \sigma_0(d) \BE + \BJ_s\;.
\end{equation}

\begin{Remark}
	\label{PPT}	
	We note that polarization is to be considered in the case of electromagnetic wave propagation that contains the electric and magnetic fields whose oscillations are perpendicular to the wave direction. The polarization is a feature of an electromagnetic wave that determines the geometrical orientation of the	oscillations \cite{shipman12}. When an electromagnetic wave propagates in a homogeneous isotropic medium, the oscillations of electric and magnetic fields are perpendicular to each other and perpendicular to wave propagation's direction. The polarization expresses the direction of the electric field. Depending on the oscillation direction of the fields, there are two types of polarization, linear and circular the latter can be specified as the left circular and right circular polarization \cite{emery+camps17}. The polarization effect is regarded in Maxwell's equation by modifying \req{eq2} in terms of the electric magnetic field $\BD$ as \cite{dorfmann+ogden03}
	\begin{equation}
	\BD = \epsilon_0(d)\BE + \BP,
	\end{equation}
	where $\BP$ represents the polarization of the electromagnetic wave propagation. In our work we further assume $\BP$ has negligible effect, thus $\BP=\textbf{0}$. Therefore, this effect is not anymore take into account in our study.	
\end{Remark}
%
Here, $\BJ_s$ denotes the specific current density. Also, material constants ${\epsilon_0}$, ${\mu_0}$, and $\sigma_0$ represent the permittivity,
the permeability, and the conductivity of the medium,
respectively. We further note that ${\epsilon_0}$, $ {\mu_0}$, and $\sigma_0$ are piecewise smooth, real, bounded, and positive and may vary in different Cartesian coordinates. We need to mention that dependency of these variables to the crack phase-field is further elaborated in Section 2.2.4. Thus, we have for every $\Bx\in\Omega$:
\begin{equation}
	\begin{aligned}
	&0<{\epsilon_{0,min}}<{\epsilon_0}<{\epsilon_{0,max}}<\infty,\\[2ex]
	&0<{\mu_{0,min}}<{\mu_0}<{\mu_{0,max}}<\infty,\\[2ex]
	&0<{\sigma_{0,max}}<\infty.\\[0.25cm]
	\end{aligned}
\end{equation}
%
Additional to \req{eq1} and \req{eq2}, the Dirichlet
boundary condition are for the electric field and magnetic vector potential:
\begin{equation}
	\BE = \bar{\BE}(\Bx,t) \ \textrm{on}\ \partial^1_ D\Omega~~\text{and}~~
	\BA = \bar{\BA}(\Bx,t) \ \textrm{on}\ \partial^1_ D\Omega.
	\label{BC.F}
\end{equation}
On the boundary of the free space, outer Dirichlet and Neumann boundary conditions have to be satisfied as:
\begin{equation}\label{eq:BC:A}
	\bar{\BA}=\textbf{0} \quad \text{on}\;\;\partial^0_ D{\Omega} \AND \nabla \bar{\BA} \cdot \BN = \textbf{0}\quad \text{on}\;\;\partial^0_ D{\Omega},
\end{equation}
respectively, such that: 
\begin{equation*}
	\displaystyle\lim\limits_{\Bx\to \infty}\bar{\BA}(\Bx,t) \rightarrow \textbf{0}.
\end{equation*}
In the case of the magnetized material, the constitutive relation
between magnetic flux $\BH$ and magnetic field $\BB$ needs to be
further modified. This adjustment will help to accommodate our proposed formulation of a magneto-mechanical problem of fracturing solids. To do so, we describe three common types of macroscopic magnetization responses as follows:
%
\begin{figure}[!t]
	\centering
	{\includegraphics[clip,trim=0cm 23.3cm 0cm 0cm, width=15cm]{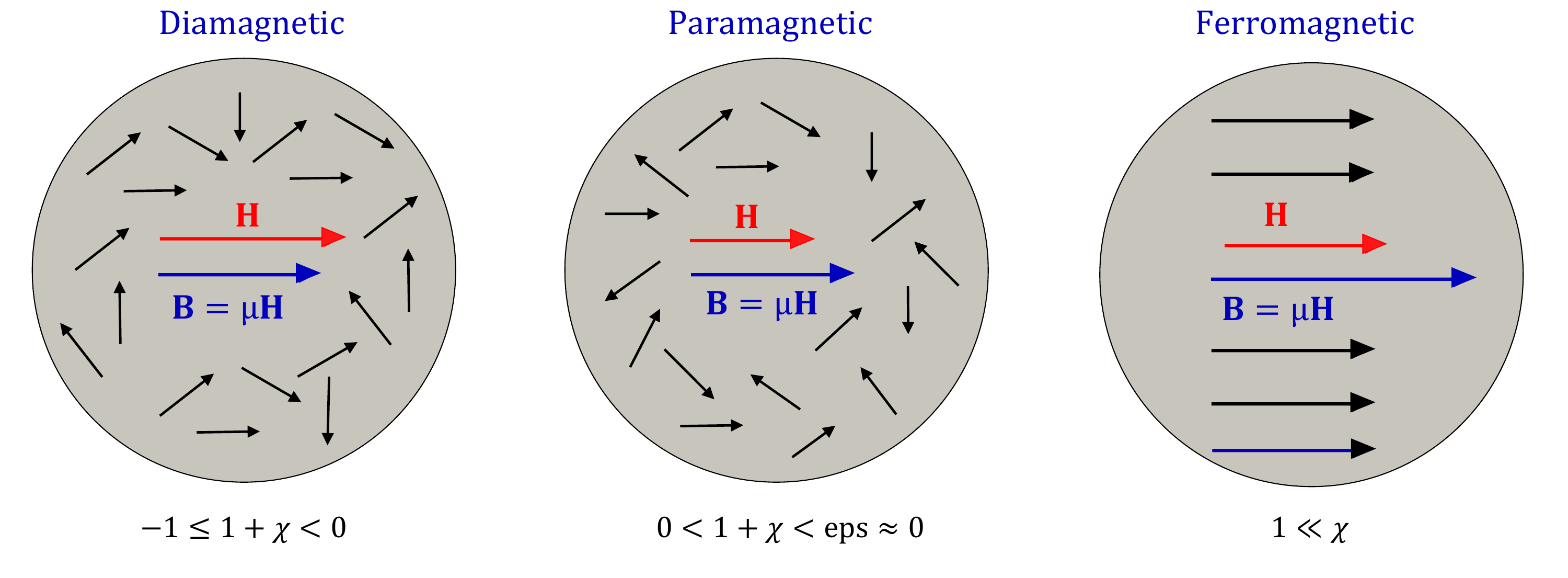}}  
		\vspace{-0.5cm}
	\caption*{\hspace*{0.1cm}(a)\hspace*{4.5cm}{(b)}\hspace*{4.5cm}{(c)}}
	\caption{Partial alignment of the magnetic moments in the presence of the external magnetic field for classified based on relative magnetic permeability of materials, (a) diamagnetic with $-1\le1+\chi<0$, (b) for paramagnetic with $0\le1+\chi<\text{eps}$, and (c) for ferromagnetic $1\ll\chi$  in parallel orientation of magnetic moment in the presence of the external magnetic field.
	}
	\label{intro_magnetic}
\end{figure}
\begin{itemize}
	\item \textbf{Diamagnetic.} These types of materials are those which become weakly magnetized when they are affected by an external magnetic field. The tendency of magnetization for paramagnetic materials is to move in the direction of a strong magnetic field towards weak parts of the external magnetic field \cite{spaldin2010magnetic}. So, diamagnetic materials are magnetized in the opposite direction of the external magnetic field. Moreover, these material present no magnetic hysteresis. 
	\item \textbf{Paramagnetic.} Similar to the diamagnetic response, this type of material becomes weakly magnetized. Unlike diamagnetic materials, the tendency of magnetization for paramagnetic materials is to move in the direction of weak magnetic fields to strong ones \cite{spaldin2010magnetic}. So, paramagnetic materials are magnetized in the direction of the external magnetic field. Here, small magnetic hysteresis could be observed . 
	\item \textbf{Ferromagnetic.} This class of magnetic materials is described through a permanent magnetization effect and mainly has a profound response to magnetic fields. These materials are characterized by a highly nonlinear response in which magnetic hysteresis typically occurs \cite{fliegans2020modeling}.
\end{itemize}
To formulate a different class of magnetic material, the constitutive relation between $\BB$ and $\BH$ is described through: 
\begin{equation}\label{magnetic_B}
	\BB = \boldsymbol{\mu}\cdot\BH \WITH \boldsymbol{\mu}=\mu_0(\BI+\boldsymbol{\chi}),
\end{equation}
in which magnetic susceptibility $\boldsymbol{\chi}$ (as a
dimensionless quantity), measures the  degree of magnetization
response of material if there is an external magnetic field leading to:
\begin{equation}\label{magnetic_susceptibility}
	\BM=\boldsymbol{\chi}\cdot\BH \WITH \boldsymbol{\chi}=\frac{\partial \BM}{\partial \BH}.
\end{equation}
The term $\BI+\boldsymbol{\chi}$ in \req{magnetic_B} typically refers to the \textit{relative magnetic permeability} of materials. Thus, the coupling of \req{magnetic_B} and \req{magnetic_susceptibility} results in:
\begin{equation}\label{eq_BH}
	\BB =\mu_0(\BI+\boldsymbol{\chi}):\BH={\mu_0}(\BH+\BM). 
\end{equation}
On the basis of the magnitude of magnetic susceptibility
$||\boldsymbol{\chi}||_2$  in \req{eq_BH}, magnetized material is
classified into the above mentioned:
\begin{equation}
	\begin{aligned}
		&\texttt{Diamagnetic: }  -1\le||\boldsymbol{\chi}||_2<0,\quad\\[2ex] 
		&\texttt{Paramagnetic: }  0<||\boldsymbol{\chi}||_2<\text{eps}\approx0,\quad\\[2ex] 
		&\texttt{Ferromagnetic: }  ||\boldsymbol{\chi}||_2\gg1,
	\end{aligned}
\end{equation}
see Figure \ref{intro_magnetic}. In our case, we use the magnetic susceptibility as an isotropic tensor for the ferromagnetic response: $\boldsymbol{\chi}=\chi_0\BI$ with $\chi_0\gg1$, which is valid for a variety of amorphous solids or materials with a uniform crystal structure \cite{bruzzese2022theory}. 
Since $\BB$ is solenoidal (divergence-free)
according to Maxwell's equations, we know through \req{eq:B_curl}
that $\BB$ must be the curl of the magnetic vector potential field $\BA$, so $\BB=\nabla\times\BA$. In a common way, the electromagnetic fields, $\BE(\Bx,t)$, $\BH(\Bx,t)$, $\BD(\Bx,t)$ and $\BB(\Bx,t)$ are presumed to be time-harmonic. It means that the fields can harmonically oscillate with a single frequency $\omega$. In such cases, they can be written as  
\begin{equation}
	\begin{aligned} 
	\label{harmonic}	
	\BE(\Bx,t) & =\exp(-i\omega t)\BE(\Bx,\omega),\\[2ex] 
	\BH(\Bx,t) & =\exp(-i\omega t)\BH(\Bx,\omega),\\[2ex]
	\BJ_s(\Bx,t) & =\exp(-i\omega t)\BJ_s(\Bx,\omega).
	\end{aligned}	
\end{equation}
The first and the second-order derivatives are
\begin{equation}
\begin{aligned}
\frac{\partial \BE(\Bx,t)}{\partial t}= -i\omega\BE(\Bx,t), \quad
 && \frac{\partial^2 \BE(\Bx,t)}{\partial t^2}= -\omega^2 \BE(\Bx,t), \\[2ex]
\frac{\partial \BH(\Bx,t)}{\partial t}= -i\omega\BH(\Bx,t), \quad
&& \frac{\partial^2 \BH(\Bx,t)}{\partial t^2}= -\omega^2\BH(\Bx,t), \\[2ex]
\frac{\partial \BJ_s(\Bx,t)}{\partial t}= -i\omega\BJ_s(\Bx,t), \quad
&& \frac{\partial^2 \BJ_s(\Bx,t)}{\partial t^2}= -\omega^2\BJ_s(\Bx,t).
\end{aligned}
\label{derivatives}
\end{equation}

\begin{figure}[!t]
	\centering
	{\includegraphics[clip,trim=0cm 27cm 0cm 0cm, width=16cm]{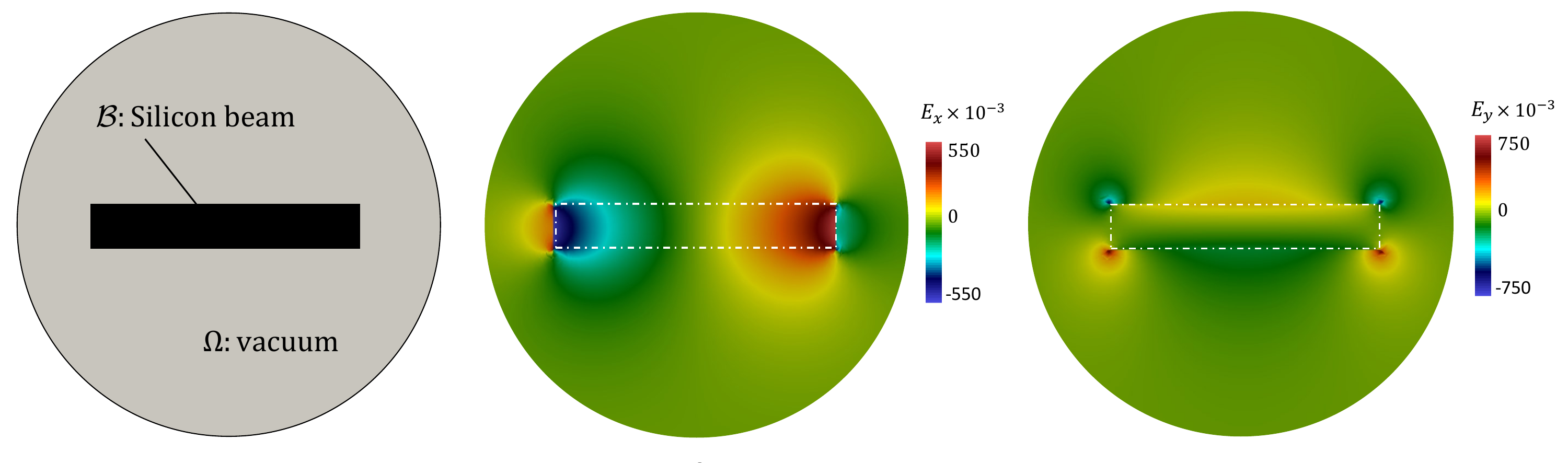}}  
	\vspace{-0.9cm}
	\caption*{\hspace*{0.3cm}(a)\hspace*{4.2cm}{(b)}\hspace*{5.2cm}{(c)}\hspace*{1cm}}
	\caption{The geometry of BVP, and the the representation of the $x$- and $y$- component of the electric field $\BE$.
	}
	\label{example0_BVP}
\end{figure}

Therefore, the time-harmonic Maxwell's equations are obtained by substituting \req{harmonic} and \req{eq2} into \req{eq1}
\begin{equation}
\begin{cases} \label{eq25}
	\nabla \times \BE(\Bx,t)  = i \omega \mu_0 \BH(\Bx,t) & \quad \text{in}\; $\Omega \times(0,\text{T})$\\\\
	\nabla \times \BH(\Bx,t)  =-i \omega \left(\epsilon_0 + i \frac{\sigma_0}{\omega}\right)\BE(\Bx,t) + \BJ_s(\Bx,t) & \quad \text{in}\; $\Omega \times(0,\text{T})$
\end{cases}.
\end{equation}
We now take a curl function, i.e., $\nabla \times$, from the first equation of \req{eq25}, which results in:
\begin{equation}\label{eq25_curl1}
	\nabla \times \BH(\Bx,t)=\nabla \times \left(\frac{1}{i \omega \mu_0} \nabla \times \BE(\Bx,t)\right).
\end{equation}
Next, we use the curl function from the second expression
of \req{eq25} and with \req{eq25_curl1}, we obtain the equation for electric field $\BE(\Bx,t)$ as 
\begin{equation}
	\nabla \times \left(\frac{1}{\mu_0}\nabla \times \BE(\Bx,t)\right) - \omega^2 \left(\epsilon_0 +i \frac{\sigma_0}{\omega} \right)\BE(\Bx,t)=i\omega \BJ_s(\Bx,t).
	\label{eqEH}
\end{equation}
In a similar manner, we can derive the equation for the magnetic field $\BH(\Bx,t)$ by taking the curl from the second expression of \req{eq25} as
\begin{equation} \label{eq28}
	\nabla \times (\nabla \times \BH(\Bx,t))=-i\omega\left(\epsilon_0+i\frac{\sigma_0}{\omega}\right)\nabla \times \BE(\Bx,t) + \nabla \times \BJ_s(\Bx,t),
\end{equation}
where $\nabla \times \BE(\Bx,t)=i\omega\mu_0\BH(\Bx,t)$. Therefore, we obtain
\begin{equation}\label{eq29}
	\nabla \times (\nabla \times \BH(\Bx,t))-\omega^2 \mu_0\left(\epsilon_0 + i\frac{\sigma_0}{\omega}\right)\ \BH(\Bx,t) = \nabla \times \BJ_s(\Bx,t).
\end{equation}
The time-domain vector wave equations in terms of the electric field
$\BE(\Bx,t)$ and the magnetic field $\BH(\Bx,t)$ are further
simplified by substituting \req{derivatives} into \req{eqEH}
and \req{eq29}. This yields:
\begin{equation}\label{eq5}
\begin{cases}
	\epsilon_0 \frac{\partial^2 \BE}{\partial t^2} + \sigma_0
	\frac{\partial \BE}{\partial t} + \nabla \times
	\left(\frac{1}{\mu_0} \nabla \times \BE \right) = -\frac{\partial
		\BJ_s}{\partial t} & \quad \text{in} \quad $\Omega\times(0,\text{T})$  \\\\
	\epsilon_0 \frac{\partial^2 \BH}{\partial t^2} + \sigma_0
	\frac{\partial \BH}{\partial t} + \nabla \times
	\left(\frac{1}{\mu_0} \nabla \times \BH \right) = \frac{1}{\mu_0}\nabla \times \BJ_s & \quad \text{in} \quad $\Omega\times(0,\text{T})$\\
\end{cases}.
\end{equation}
In order to demonstrate the modeling capacity of \req{eq5}, a
representative simulation is presented in
Figure \ref{example0_BVP}. In this boundary-value problem, the
evolution of the electric field $\BE$ is represented in the $x$ and
$y$ directions under the variation of the specific current $\BJ_s$. In
this example, the specific current $\BJ_s$ is defined as a bi-linear
function of the nodal coordinates and time. As a boundary condition,
the electric field is set to $\BE=\textbf{0}$ at the boundaries of the
vacuum. The results show that the variation of the electric field is
symmetric along the $x$ and $y$ directions as the electric specific current is uniformly applied on the silicon beam. In fact, this is expected due to symmetric response of electric field imposed to domain.

Next, we reduce \req{eq5} toward Maxwell's equations in terms of vector potential. The Maxwell-Faraday law in \req{eq1}, based on magnetic vector potential $\BA(\Bx,t)$ given in \req{eq:B_curl}, is then related to electric field through:
\begin{equation}\label{eq:E_A}
	\nabla \times \BE = -\frac{\partial({\nabla \times \BA})}{\partial{t}} \quad \rightarrow \quad \exists\;\Phi(\Bx,t)\;:\;\; -\nabla \Phi=\BE+\frac{\partial \BA}{\partial t}.
\end{equation}
Here, $\Phi(\Bx,t)$ is so-called the electric potential, see for a
detailed
discussion \cite{jackson1999classical,melia2001electrodynamics}. From
the mathematical point of view,  the homogeneous Maxwell equations for
Maxwell-Faraday law \req{eq1}$_1$, and Maxwell-Amp\`ere law
in \req{eq1}$_2$ results in the existence set of $(\BA,\Phi)$, and thus implies \req{eq:E_A}.  We refer interested readers to \cite{van2021comment} for more details.
Considering Maxwell-Amp\`ere law in \req{eq1}, and constitutive equations \req{eq2}, yields:
\begin{equation}
	\frac{1}{\mu_0}\nabla \times \BB=\epsilon_0(d)\frac{\partial \BE}{\partial t}+\BJ
\end{equation}
that is
\begin{equation}\label{eq:curl_divA}
	\nabla \times \left(\frac{1}{\mu_0}\nabla \times \BA\right)
	=-\epsilon_0(d)\left[\frac{\partial \nabla \Phi}{\partial t}
	+\frac{\partial^2 \BA}{\partial t^2}\right]
	-\sigma_0(d) \left[\nabla \Phi+\frac{\partial \BA}{\partial t}\right]
	+ \BJ_s
\end{equation}
where \req{eq:E_A} and \req{eq2} are used. 
Next, we assume that for a time dependent magnetic vector potential $\BA(\Bx,t)$, and the electric potential $\Phi(\Bx,t)$, we have:
\begin{equation}\label{Lorenz_gauge}
	\Div \BA + \mu_0\epsilon_0\frac{\partial \Phi(\Bx,t)}{\partial t}=0,
\end{equation}
which is the so-called Lorenz gauge condition \cite{maudlin2018ontological,melia2001electrodynamics}. Thus, \req{eq:curl_divA} together with  \req{Lorenz_gauge} reads:
\begin{equation}\label{eq:final_A}
		\fterm{
	\frac{1}{\mu_0}\nabla^2{\BA}
	=\epsilon_0(d)\left[\frac{\partial^2 \BA}{\partial t^2}\right]
	+\sigma_0(d) \left[\nabla \Phi+\frac{\partial \BA}{\partial t}\right]
	- \BJ_s
}
\end{equation}
Here, the identity $\nabla \times 	\nabla \times \BA= \nabla(\Div \BA)-\nabla^2{\BA}$ is used. 

\noii{We note that, the magnetic potential $A(x,t)$ as a solution of \req{eq:final_A} is related to crack parameter $d$ by PDE coefficients. So,  the permittivity, the permeability, and the conductivity of the medium $\big(\epsilon_0(d),\mu_0(d),\sigma_0(d)\big)$ are function of $d$, see Section 2.2.4.}

The setting presented in \cite{powell2008two} constitutes the starting point to formulate the mechanical response within the electromagnetic analysis where the electric scalar potential $\nabla{\Phi}$ in \req{eq:final_A} is a \textit{given quantity}. 
Following  \cite{powell2008two, zhang2021multi}, the magneto-static
state $\nabla{\Phi}$ has the direction
orthogonal to the $x$-$y$ plane, i.e., $\nabla{\Phi}(t)=[0, 0, \nabla\Phi_0(t)]$, where $\Phi_0(t)$ is the pulse
power supply parameter which is approximated as:
\begin{equation}\label{eq:phi_poission}
	\nabla\Phi_0(t) = a + b\,e^{-c\,t}.
\end{equation}
It is determined as a function of the constants $a$, $b$ and $c$
\cite{powell2008two, zhang2021multi} \footnote{
	Clearly, the representation on electric scalar potential given in \req{eq:phi_poission} can be derived through Maxwell's equations, as well, and \textit{not} as a given quantity. By means of Maxwell-Gauss law in \req{eq1}$_3$, and constitutive equation \req{eq2}$_1$, we have:
	\begin{equation}\label{eq_phi1}
		-\nabla^2{\Phi}-\nabla.\left(\frac{\partial \BA}{\partial t}\right)=-\frac{\rho_0(d)}{\epsilon_0(d)}
	\end{equation}
	where \req{eq:E_A} is used. Then, by imposing Lorenz gauge condition in \req{eq_phi1}, results in the wave-like equation by:
	\begin{equation}\label{eq_phi2}
		c^2_{\Phi}\frac{\partial^{2} \Phi}{\partial^{2} t} = \nabla^{2}\Phi + f_{\Phi}
		\quad \text{in} \quad \Omega\times(0,\text{T}),
	\end{equation}
	with
	\begin{equation}
		f_{\Phi}= -\frac{\rho_0(d)}{\epsilon_0(d)}
		\AND c^2_{\Phi}=\mu_0(d)\epsilon_0(d)
	\end{equation}
	The the electric potential PDE given \req{eq_phi2} is a linear second-order differential equation with the variable $\Phi\in H^1(\Omega,\mathbb{R})$ represents the wave-like propagates in the free
	space with the constant speed of $c_\Phi$, and $f_{\Phi}$ states the
	source of the wave propagation. Thus, set of equations \req{eq:final_A} and \req{eq_phi2} are replaced by \req{eq5} to represent the Maxwell's equations in terms of vector potential.}. In all our representative
numerical examples in Section 4, these values are set as $a=-1\times10^{-3}$ V/m,
$a=-1\times10^{-3}$ V/m, and $c=1\times10^{-2}$ s$^{-1}$.

\sectpc{Fracture contribution} 
In the smeared fracture framework, a sharp crack interface denoted by $\mathcal{C}$ for satisfying the continuity of the crack topology is further regularized, which is denoted as $\mathcal{C}_{l}$ as outlined \cite{bourdin+francfort+marigo08}. To govern a regularized fracture surface $\mathcal{C}_{l}$, it is required to incorporate a continuous field variable -- the so-called order parameter -- denoted by $d$, which differentiates between multiple physical phases within a given system through a smooth transition. In the context of fracture,
such an order parameter (termed the crack phase-field) describes the smooth transition
between the fully broken ($d=1$) and intact material phases ($d=0$),
thereby approximating the sharp crack discontinuity. This geometrical
perspective is in agreement with the framework
of~\cite{bourdin2000numerical}, which was conceived as a
$\Gamma$-convergence regularization of the variational approach to
Griffith fracture~\cite{francfort1998revisiting}. 
A variety of research studies have recently extended the phase-field approach to fracture toward the
cohesive-frictional materials~\citep{kienle2019}
including thermal effects~\citep{dittmann2020,ruan2022thermo,peng2023meso}, ductile failure \cite{ambati+kruse+lorenzis16,yin2020ductile,noii2021bayesian1}, hydraulic fracture~\citep{heider2020phase,ulloa2022variational,noii2019phase}, stochastcic analysis \cite{noii2022bayesian1,noii2022probabilistic, noii2021bayesian,noii2022bayesian}, degradation of the fracture toughness \cite{yin2020ductile}, topology optimization \cite{noii2023level}, and multi-scale approach \cite{liu2022phase,noii2020adaptive,aldakheel2020global,aldakheel2021multilevel}, and electro-mechanical approach \cite{hageman2022electro,wu2022crack, zhao2022phase} among others. 
In this manuscript, for the case of
isotropic materials, the regularized functional is given by
\begin{equation}
	\calC_{l}(d) = \int_{\calB} \gamma(d, \nabla d) \, \text{d}\B, 
	\label{s2-gamma_l}
\end{equation}
with possitivness for crack disspation as:
\begin{equation}
	\frac{d}{dt} \calC_l(d)  \ge 0\;.
	\label{gamma-evol}
\end{equation}
In line with standard phase-field models, a general surface density function for the isotropic part $\gamma(d, \nabla d)$ is defined as 
\begin{equation}
	\gamma(d, \nabla d):=\frac{1}{c_d}\, \bigg(\frac{f(d)}{l_d} +\frac{l_d}{2} \nabla d \cdot \nabla d \bigg) \WITH c_d:=4\int_0^1\sqrt{\omega(b)}\,db ,
\end{equation}
where $\omega(d)$ is a monotonic and continuous \emph{local fracture energy function} such that $\omega(0)=0$ and  $\omega(1)=1$. A variety of suitable choices for $\omega(d)$ are available in the literature~\citep{kuhn2015,wu2017,wu2018}. Here, the widely adopted linear and quadratic formulations are considered, which yield, models with and without an elastic stage, respectively. Specifically, we define 
\begin{equation}
	f(d) : = \begin{cases}  d\phantom{^2} \implies c_d=8/3 \quad &\text{model with an elastic stage}, \\ d^2 \implies c_d=2 \quad &\text{model without an elastic stage}.  \end{cases}
	\label{wd}
\end{equation}
Following \cite{miehe2010phase}, the local evolution of the crack phase-field equation in the domain 
$\calB$ is:
\begin{equation}\label{euler-eq-d}
	\fterm{
		[ \, d - l_d^2 \Delta d \, ] + \eta_d \dot{d} + (1-\kappa)(d-1) {\calH} = 0 
		\, ,
	}
\end{equation}
augmented with its homogeneous Neumann boundary conditions, i.e., $\nabla d \cdot 
\Bn = 0$ on $\partial\calB$.
Following \cite{noii2021quasi,noii2020adaptive}, the small residual 
scalar
$0<\kappa\ll1$ is introduced to prevent numerical instabilities,
which states the {\it third} equation in the coupled system. Additionally, the damage viscosity material parameter denoted by $\eta_d  \ge 0$ is used to characterize 
\textcolor{black}{the viscosity term of the crack propagation}. 
The maximum absolute value for the crack driving state function denoted by $\calD$ is defined by the crack driving force $\calH$, which reads
\begin{equation}
	\calH = \max_{s\in [0,t]} \calD(\Bx,s) \ge 0,
	\label{driving-force}
\end{equation}
that accounts for the irreversibility of the crack phase-field evolution by filtering out a maximum value of $\calD$. We define the crack driving state function in Section 3. 

In summary, following formulation has to be solved for three-filed $(\Bu,d,\BA)$ multi-physics problem.

\begin{form}[Strong form of the Euler-Lagrange equations]
	\label{form_4}
	Let $\epsilon_0$, $\sigma_0$, $\mu_0$, $G_c$, $K$, and $\mu$ be given with the initial conditions $\bm u^0=\bm u(\bm{x},0)$,  $d^0=d(\bm{x},0)$, and $\BA^0=\BA(\bm{x},0)$. For the loading increments $n=1,2,\ldots, N$, we solve a displacement equation where we seek $\bm u:= \bm u^n: \calB \rightarrow \mathbb{R}^{\delta}$ such that
	\begin{align*}
		-\div (\Btau )= \Bzero \quad \quad \quad \quad &in\; \calB,\\
		{\bm u} =\bar \Bu \quad \quad \quad \quad &on \;  \partial_D\calB,\\
		\Btau\cdot \Bn = \bar{\Bt}\quad \quad \quad \quad &on \;  \partial_N\calB,
	\end{align*}
	in terms of the stress tensor $\Btau $ defined in \req{eq17} and the given displacement field $\bar \Bu$. The phase-field system consists of four parts: the PDE, the
	inequality constraint and a compatibility condition along with the Neumann-type boundary conditions.
	Find $d:=d^n : \calB \rightarrow [0,1]$ such that
	\begin{align*}
		(1-\kappa)(1-d) {\calH(I_1,I_2,t)}- [ d - l_d^2 \Delta d \, ]=\eta_d \dot{d} \qquad &in \; \calB,\\
		\dot{d}\geqslant 0 \qquad &in \; \calB,
	\end{align*}
	in terms of the crack driving force given \req{driving-forceD}, along with $\nabla d \cdot \Bn=0 \; \text{on}\; \partial \Omega$. along with a second-order hyperbolic problem for the following magnetic vector potential. Find $\BA:= \BA^n: \Omega \rightarrow \mathbb{R}^{\delta}$ such that
	\begin{align*}
		\quad\quad\quad\;		
		\frac{1}{\mu_0}\nabla^2{\BA}
		=\epsilon_0(d)\left[\frac{\partial^2 \BA}{\partial t^2}\right]
		+\sigma_0(d) \left[\nabla \Phi+\frac{\partial \BA}{\partial t}\right]
		- \BJ_s \qquad \text{in} \;\Omega.
	\end{align*}
	where the Lorenz gauge condition is considered.  
\end{form}

\sectpc{Transition structure from undamaged to fully damaged states} 
In this section, we provide a transition rule to connect the
transition between the intact (solid) region and the fracture
domain. In fact, the coupling of electromagnetic response to the crack phase field is achieved by introducing constitutive functions, which are characterized by degraded related material constants. Following \cite{miehe2015phase}, two state functions need to be defined. The first set of the formulation is denoted as $\texttt{P}^s(\Bx) \in \calB \backslash \calC_l$, which exists in the solid part of materials to describe the degradation of elastic stored contribution in energy. In order to extend $\texttt{P}^s(\Bx)$ to the entire domain $\calB$, we define $\widehat{\texttt{P}}^s(\Bx,d)$ as: 
\begin{equation}
	\texttt{P}^s(\Bx) \in \calB \backslash \calC_l \quad \xrightarrow \quad\quad
	\widehat{\texttt{P}}^s(\Bx,d) \in \calB  \WITH \widehat{\texttt{P}}^s(\Bx,d)=\calG^s(d)\texttt{P}^s(\Bx)\;,
\end{equation}
with decreasing crack phase-field. Here, $\calG^s(d)$ has a descending
order in terms of the crack phase-field. This set of functions can be
used to describe, e.g., heat permeability,  fluid permeability, or in
our study, the permittivity and the permeability of the ferromagnetic
medium. More precisely, all sets of solid response are a subset of
$\widehat{\texttt{P}}^s(\Bx,d)$ such that in the limit case this state
function
$\widehat{\texttt{P}}^s(\Bx,d\rightarrow0)\rightarrow\texttt{P}^s(\Bx)$
thus resembling classical constitutive state variables for the solid response.

The second set of the formulation is denoted as
$\texttt{P}^d(\Bx) \in \calC_l$, which exists in the damaged domain to
describe the fractured constitutive response. To further extend
$\texttt{P}^f(\Bx)$ to the entire domain, we define:
\begin{equation}
	\texttt{P}^f(\Bx) \in \calC_l \quad \xrightarrow \quad\quad
	\widehat{\texttt{P}}^f(\Bx,d) \in \calB  \WITH \widehat{\texttt{P}}^f(\Bx,d)=\calG^f(d)\texttt{P}^f(\Bx)\;,
\end{equation}
with increasing crack phase-field. Here, $\calG^f(d)$ has an ascending
order in terms of the crack phase-field, which has a inverse effect of $\calG^s(d)$ into the constitutive response. As a result, we have the following properties:
\begin{equation}\label{prop_deg}
	\begin{aligned}
		&\texttt{solid part:}  \quad\quad\;\;\calG^s(d\rightarrow0)\rightarrow1, \quad
		\calG^s(d\rightarrow1)\rightarrow0\WITH \partial_{d}\calG^s\le0\\
		&\texttt{fracture part:} \quad \calG^f(d\rightarrow0)\rightarrow0,\quad
		\calG^f(d\rightarrow1)\rightarrow1  \WITH \partial_{d}\calG^f\ge0
	\end{aligned}
\;\;\;.
\end{equation}
Thus, $\calG^f(d=0)$ is zero in the solid domain, while $\calG^s(d=1)$ is zero in the fracture domain. 
In this contribution, we propose the following formulation for
$\calG^s(d)$ and $\calG^f(d)$ of order $m$ for a given set of
parameters $(c_1,c_2)$, for the solid contribution as:
\begin{equation}\label{Gs}
	\begin{aligned}
    \calG^s(d)=\Big(\frac{c_2-d}{c_2-c_1}\Big)^m \WITH \calG^s(d)=1 \;\; \text{if}\;\; d<c_1 \AND  \calG^s(d)=0  \;\; \text{if}\;\;  d>c_2\;,
	\end{aligned}
\end{equation}
and the fracture contribution as:
\begin{equation}\label{Gf}
	\begin{aligned}
		\calG^f(d)=\Big(\frac{d-c_1}{c_2-c_1}\Big)^m \WITH \calG^f(d)=0 \;\; \text{if}\;\; d<c_1 \AND  \calG^f(d)=1  \;\; \text{if}\;\;  d>c_2,
	\end{aligned}
\end{equation}
\begin{figure}[!t]
	\centering
	{\includegraphics[clip,trim=0.6cm 25.3cm 0cm 0cm, width=16.3cm]{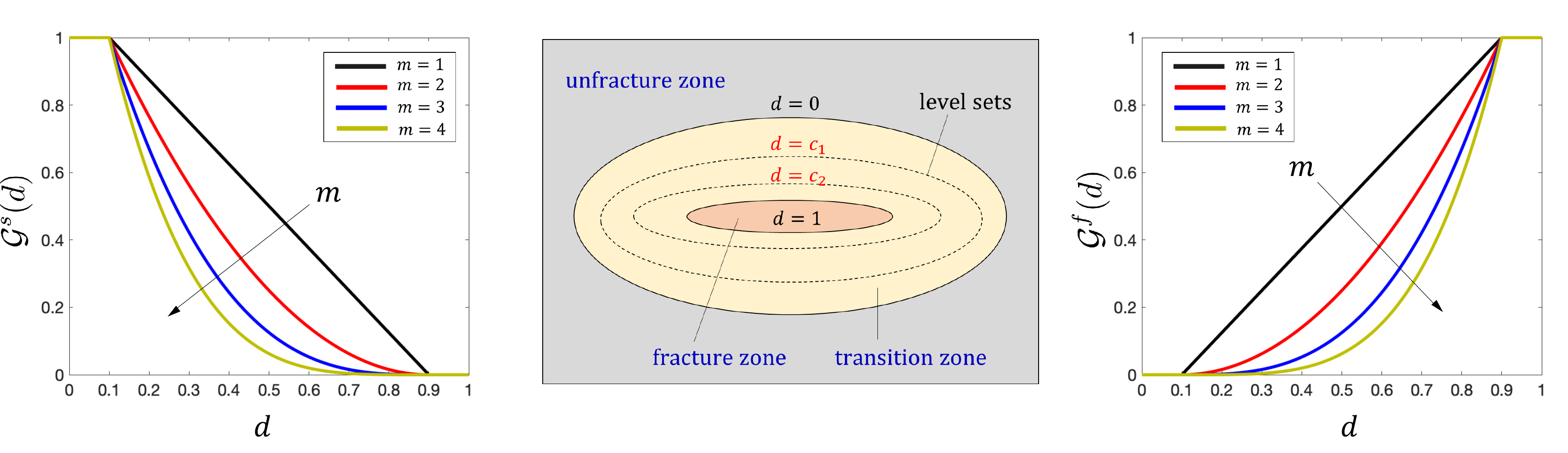}}  
	\vspace{-0.8cm}
    \caption*{\hspace*{1cm}(a)\hspace*{4.9cm}{(b)}\hspace*{5.3cm}{(c)}\hspace*{0.2cm}}
	\caption{ (a) Degradation function for solid material denoted as $\calG^s(d)$ (b) with representation the fracture domain if $d>c_2$ and unfracture  domain  if $d<c_1$, and (c) Degradation function for regularized crack phase-field $\calG^f(d)$. Here, we used $m=\{1,2,3,4\}$ with adjustable constants $(c_1,c_2)=(0.1,0.9)$. 
	}
	\label{tran_rule}
\end{figure}
which is in line properties \req{prop_deg}. In the following numerical study we use the quadratic order for both $\calG^s(d)$ and $\calG^f(d)$, so $m=2$, and we further set $(c_1,c_2)=(0.1,0.9)$, see Figure \ref{tran_rule}.

\begin{Remark}
	\label{remk1}	
For the linear order $m=1$ our proposed degradation function is
similar with linear indicator functions given
in \cite{lee2016pressure}, and for the quadratic order $m=2$ with
constants $(c_1,c_2)=(0,1)$ it mimics quadratic indicator functions given in \cite{miehe2015phase}. Hence, \req{Gs} and \req{Gf} are generalized formulation versus \cite{miehe2015phase} and \cite{lee2016pressure} .
\end{Remark}

To denote the effect of the magnetization-induced
cracking in the material, we propose an anisotropic electromagnetic
material constants in Maxwell's equations. This provides a simple
constitutive assumption for an electromagnetic material to represent the
transition rule. These constants can be re-formulated through crack phase-field $d$. Thus, by means of state functions for solid and fracture regions described in \req{Gs} and \req{Gf}, the definition of the
electromagnetic material constants is then proposed as follows:
\begin{equation}\label{aniso_permability}
	\begin{array}{ll}
		 &\texttt{permeability:} \quad{\mu_0}(d) =\calG^s(d)\mu^s_0+\calG^f(d)\mu^f_0, \\\\
		 &\texttt{permittivity:} \quad{\epsilon_0}(d) =\calG^s(d)\epsilon^s_0+\calG^f(d)\epsilon^f_0, \\\\
		  &\texttt{conductivity:} \quad{ \sigma_0}(d) =\calG^s(d) \sigma^s_0+\calG^f(d) \sigma^f_0, \\\\
		  &\texttt{electric charge density:} \quad{\rho_0}(d) =\calG^s(d)\rho^s_0+\calG^f(d)\rho^f_0.
	\end{array}
\end{equation}
%
%
The first term in the anisotropic materials expression
in \req{aniso_permability} for a set of
$(\mu_0,\epsilon_0, \sigma_0,\rho_0)$ can be indicated as a classical
isotropic intrinsic material property for permeability,  permittivity,
conductivity, and electric charge density, respectively. Also, in
general, the number of constants
$(\mu^f_0,\epsilon^f_0, \sigma^f_0,\rho^f_0)$ are much less than
material constants corresponding to solid parts and greater than the
number of constants describing vacuum counterparts. We note that in
our study, the second term in \req{aniso_permability} is also a material property, which one can relate to classical crack jump, e.g.,
in failure mechanics of a thermo-elastic solid
in \cite{miehe2015phase}, or in hydraulic fracturing of fluid-saturated porous media \cite{miehe2016phase,noii2021bayesian,noii2022bayesian}.
\sectpa[Section3]{Energy quantities and variational principles}
\label{Bayesian}
To outline the variational formulation setting, it suffices to define
the constitutive energy density functions $W_{elas}$, $W_{mag}$,
$W_{mos}$, and $W_{frac}$ corresponds to elastic, magnetostrictive,
magnetization, and fracture contributions, respectively. These energy
quantities lead to establishing the multi-field evolution problem in
terms of the primary fields described in Section 2. To this end, the
related constitutive relations are provided. On this basis, a coupling
between the magnetostriction and the mechanics is formulated to
investigate the mechanical deformation of materials under the magnetic
forces. Furthermore, the cracking response under the magnetic forces is
examined. The energetic crack driving force is governed by the
coupling magneto-mechanical effects under the magneto-static state. We
extend the variational formulations of the coupled multi-physics
system within the magneto-mechanical formulation of fracturing
solids. We complete this section by providing a compact algorithmic
framework for magnetostrictive-induced cracking.
\sectpb[Section22]{Constitutive functions}
\label{sec:Constitutive}
The coupled BVP is formulated by introducing three specific fields
for magneto-mechanically-induced cracking of the magneto-active materials
%
\begin{equation}
	\mbox{Global Primary Fields}: \  
	\BfrakU := \{ \Bu, d ,\BA \}.
	\label{primary_m3}
\end{equation}
Here, $\Bu$ is the displacement (mechanical deformation), $\BA$ denotes the 
magnetic vector potential field, and $d$ is the crack phase-field ($0\le d\le1$). From the numerical implementation standpoint, to guarantee $0\le d\le1$ holds, we project $d>1$ to 1 and $d<0$ to 0 to avoid unphysical crack phase-field solution \cite{noii2019phase}. The constitutive equations for the magnetostrictive phase-field fracture are written in terms of the set
%
\begin{equation}
	\mbox{Constitutive State Variables}: \  
	\BfrakC := \{ \Bve, d , \nabla d,\BB \}
	,
	\label{state_m3}
\end{equation}
representing a combination of a magneto-mechanical model through the Maxwell's equation and a first-order gradient damage model. 
\sectpc[Section23]{Energy quantities}
A pseudo-energy density per unit volume is then defined as $W:=W(\BfrakC)$, which is additively decomposed into an elastic contribution $W_{elas}$, a magnetization contribution $W_{mag}$, a magnetostriction contribution $W_{mos}$, and a (regularized) fracture contribution~$W_{frac}\;$ through:
\begin{equation}
	\boxed{W(\BfrakC):= 
		{W}_{elas}(\Bve ,d ) + {W}_{mag}(\Bve,\BB, d)+ {W}_{mos}(\Bve,\BB,d)+ {W}_{frac}(d, \nabla d ). } 
	\label{psuedo-energy}
\end{equation}
%
We note that $W(\BfrakC)$ is a state function that contains both energetic and dissipative contributions. With this function at hand, a pseudo potential energy functional can be constructed
\begin{equation}
		\begin{aligned}
	\calE(\Bu,\BA,d) &:= \int_{\Omega} {W}(\BfrakC)  \, \text{d}\Bx\\
	&= \int_{\calB} \Big[{{W}_{elas}(\Bve ,d )}+{W}_{frac}(d, \nabla d )\Big]  \, \text{d}\Bx
	+ \int_{\Omega} \Big[{W}_{mag}(\Bve,\BB, d)+ {W}_{mos}(\Bve,\BB,d)\Big]  \, \text{d}\Bx.
		\end{aligned}
	\label{potential-functional}
\end{equation}
All these counterparts will be elaborated in following sections.
%
\sectpc{Magnetization contribution} 
%
Magnetization is defined as a vector field describing the density of
permanent or induced magnetic \textit{dipole moments} in a magnetic
material \cite{pao+hutter75}. The magnetic dipole moments are caused by either microscopic
electric currents resulting from the motion of electrons in atoms or
the spin of the electrons or the nuclei. In a magnetic field, a
strong magnetization response is observed in ferromagnetic
materials. Any ferromagnetic material is capable of being magnetized
in the absence of an external magnetic field and becoming a permanent
magnet. It is not necessary that the magnetization distributes
uniformly within the material, so the material may present an anisotropic
response when subjected to the magnetization effect. Ferromagnetic
materials have inherent strong coupled magnetic and mechanical
behavior, so the material magnetization affects the stress response of
material via magnetostriction phenomenon, and the stress state of
the material also affects magnetization via inverse magnetostriction
\cite{daniel+bernard+hubert20}. Following \cite{fonteyn10}, the magnetization density function is formulated as:
\begin{equation}
	{W}_{mag}(I_1,I_4,d):= 
	\frac{1}{2} \sum^{4}_{i=0}\frac{g_i(I_1)}{i+1}\left(
	\frac{I_4}{B^2_{\text{ref}}} \right)^iI_4,
	\label{mag-part}
\end{equation}
with the fourth invariant $I_4=\BB \cdot \BB$.

Here, $B_{\text{ref}}$ is a reference value of the magnetic flux
density. Additionally, functions $g_i(I_1)$ for $i=\{0,1,2,3,4\}$ can
be further approximated, as it is shown in Appendix A leading to:
\begin{equation}
	\begin{aligned}
		g_0 & = \frac{3}{4}\alpha_0\exp\left(\frac{3}{4} I_1\right) -
		\frac{1}{3}(\frac{1}{\mu_0} - \alpha_5)\\[1ex]
		g_i & = \frac{3(i+1)}{4}\alpha_i\exp\left(\frac{4(i+1)}{3}
		I_1\right)  & \text{where} \quad i=\{1,2,3,4\}.
	\end{aligned}
\end{equation}

The parameters $\alpha_i$ for $i=\{0,1,2,3,4,5\}$
refer to the ferromagnetic structure
due to the magnetization effect, see \cite{fonteyn10}.
%
\sectpc{Magnetostrictive contribution} 
Magnetostriction is a constitutional characteristic of magnetic material
related to electron spin or orbit orientations, their
interactions, and molecular lattice geometry. The
magnetostriction effect does not degrade in time or during usage. A magnetostrictive material does not present the magnetostriction effect
when it is heated above its Curie temperature
\cite{flatau+dapino+calkins00}. Curie temperature is defined as a
transition point between ferromagnetic and paramagnetic behaviors. In
magnetic materials, three types of magnetostriction mechanisms are
determined to express the magnetostrictive feature of the magnetic
material. These mechanisms are paramagnetic magnetostriction,
spontaneous magnetostriction, and magnetostriction induced by magnetic
field \cite{gao+etal22}. Paramagnetism is a weak form of magnetism that is induced by an external field in the direction of the applied
magnetic field. It will disappear when the magnetic field is
removed. The spontaneous magnetization is followed by 
paramagnetism when a ferromagnetic material is cooled below its
Curie temperature, so a transition from paramagnetism to ferromagnetism
occurs, and magnetic moments cause a spontaneous magnetization effect
\cite{dapino04}. In a magnetic material, the action of the external
magnetic field \cite{gao+etal22} induces a magnetostriction effect by
rotating and moving the internal magnetic domain walls.

We first employ two additional deformation-magnetization-dependent
scalar-valued function (invariants) to represent magnetostrictive effect through:
\begin{equation}
	I_5=\BB \cdot \Bve \cdot \BB \AND I_6=\BB \cdot
	\Bve^2 \cdot \BB.
\end{equation}
The fifth and sixth invariants account for the magnetostrictive response as a function of the magnetic flux density $\BB$ and the mechanical strain $\Bve$. The magnetostrictive density function is formulated as a function of $(I_5,I_6)$ by:
\begin{equation}
	{W}_{mos}(I_5,I_6,d):= \frac{1}{2}\alpha_5I_5 + \frac{1}{2}\alpha_6I_6.
\end{equation}\label{mos_density}
Now we are able to define the magnetization vector $\BM$ as follows:
\begin{equation}\label{eq16}
	\begin{aligned}
		{\BM} (\Bve,\BB;d)
		:= - \frac{{\partial W(\BB, \Bve)}}{\partial \BB}
		&=-\displaystyle\sum_{i} \frac{{\partial W(\Bve,\BB)}}{\partial I_i}\frac{{\partial I_i}}{\partial \BB}\\
		&=-\frac{{\partial W_{mag}(\Bve,\BB)}}{\partial I_4}\frac{{\partial I_4}}{\partial \BB}
		-\displaystyle\sum^{6}_{i=5} \frac{{\partial W_{mos}(\Bve,\BB)}}{\partial I_i}\frac{{\partial I_i}}{\partial \BB}\\
		&=-\sum_{i=0}^{4}\frac{g_i}{i+1}\partial_{I_4}\frac{I_4^{i+1}}{B^{2i}_{\text{ref}}}
		-\frac{1}{2}\left(\alpha_5\BB
		\cdot \Bve -\alpha_6\BB \cdot \Bve^{2}\right).
	\end{aligned}
\end{equation}
For a detailed derivation, see Appendix A.

\sectpc{Elastic contribution} 
Since, the fracturing material behaves quite differently in bulk and shear
parts of the domain, we employ a consistent split for the strain energy
density function \req{eq:psi_iso} into \textit{tension} and \textit{compression} counterparts, respectively. 
Following additive split for strain tensor in \req{eq:additive_strain}, the effective strain energy density $\psi_e$ admits the additive decomposition 
\begin{equation}
	\psi_e^{}\big(I_1,I_2\big):=
	\psi_e^{vol}(I_1)+\psi_e^{dev}(I_1,I_2),
	\label{eq:psi_iso}
\end{equation}
with
\begin{align}
	\psi_e^{vol}\big(I_1\big)
	=\frac{K_n}{2}I^2_1
	\AND
	\psi_e^{dev}\big(I_1,I_2\big)
	=\mu\Big(\frac{I_1^2}{3}-I_2\Big), \quad \text{note} \quad K_n,\mu>0.
	\label{eq:psi_iso1}
\end{align}
Here, $K_n= \lambda+\frac{2}{3}\mu$ is the bulk modulus which includes elastic Lam\'e's constant $\lambda$ and the shear modulus $\mu$.  Additionally, $I_1$ and $I_2$ denote the invariants
\begin{equation}
	I_1=\text{tr}[\Bve] \AND I_2=\tr[\Bve^2].
	\label{eq:prin_invar}
\end{equation}
The isotropic strain energy density function given in \req{eq:psi_iso} is additively decomposed into  tension and compression contributions:
\begin{equation}
	\psi_e\big(I_1,I_2\big)=\psi_e^{+}(I_1,I_2)+\psi_e^{-}(I_1),
\end{equation}
where
\begin{align}
	{\psi_e^{+}}(I_1,I_2)={H{^+}[I_1]}\psi_e^{vol}\big(I_1\big)
	+\psi_e^{dev}\big(I_1,I_2\big)~\AND
	{\psi_e^{-}}(I_1)=\big(1-{H{^+}[I_1]}\big)\psi_e^{vol}\big(I_1\big)~.	
	\label{eq232_5}
\end{align}
{Therein, $H{^+}[I_1(\Bve^e)]$ is a \textit{positive Heaviside function} which returns one if $I_1(\Bve^e)>0$ and zero if $I_1(\Bve^e)\leq0$. The total elastic contribution to the pseudo-energy density  in \req{psuedo-energy} finally reads
	\begin{equation}
		{W}_{elas}(\Bve,d):= g_e(d) \; \big[ \psi_e^{+}(I_1,I_2)\big ] + \psi_e^{-}(I_1) ,
		\label{elas-part}
	\end{equation}
	where $g_e(d)$ is a \textit{elastic degradation function}. Here, the standard monotonically decreasing quadrature degradation function, reads as $g_e(d):=(1-\kappa)(1-d)^2 + 
	\kappa$. Here,  $\kappa$ is a very small parameter (to avoid numerical instabilities), and mathematically it is also dependents on the discretization space, see
	\cite{khodadadian2020bayesian}.
	
	\noii{We note that, the current resulting degraded stress tensor by the bulk energy is rely on the $vol/dev$ decomposition. But, it is also possible to enhance \req{elas-part}, for the material which are sensitive to the shear fracture as it is investigated in rock-like materials, e.g. \cite{wang2020phase,wang2023modeling}. So, to explore different characteristic behavior of fracture which is subjected to compression and shear modes for ferromagnetic materials, we could introduce different crack phase-field driving force. This subject is open for further investigation.}

The Cauchy stress tensor is defined as a derivative of the free
energy function with respect to the strain tensor as follows 
\begin{equation}
\begin{aligned}
	{\Bsigma}(\Bve, \BB, d) := & \frac{{\partial W}}{\partial \Bve}
	=\displaystyle\sum_{i} \frac{{\partial W}}{\partial I_i}\frac{{\partial I_i}}{\partial \Bve} \\
	= & \underbrace{g_e(d){\widetilde{\Bsigma}}^{+}({\bm{\varepsilon}})+{\widetilde{\Bsigma}}^{-}({\bm{\varepsilon}})}_{\text{elasticity term}}
	+\underbrace{\frac{1}{2}\sum_{i=0}^{4}\frac{1}{i+1}\partial_{I_1}g_i
		\frac{I_4^{i+1}}{B^{2i}_{\text{ref}}}\textbf{I}}_{\text{magnetization term}}\\
	&+\underbrace{ \frac{1}{2}\alpha_5
		\BB \otimes \BB + \frac{1}{2}\alpha_6\left( \BB \otimes \BB \cdot \Bve
		+ \Bve \cdot \BB \otimes \BB \right)}_{\text{magnetostrictive term}}.
\end{aligned}
\label{eq12}
\end{equation}
with,
\begin{equation}\label{eq232_8}
	{\widetilde{\Bsigma}}^{+}({\bm{\varepsilon}})=K_nH{^+}(I_1)(\bm \varepsilon:\textbf{I})\textbf{I}+2\mu{\bm \varepsilon}^{dev},\quad \text{and} \quad {\widetilde{\Bsigma}^{-}}({\bm{\varepsilon}})=K_n\big(1-H{^+}(I_1)\big)(\bm \varepsilon:\textbf{I})\textbf{I}.
\end{equation}
%
In addition to the Cauchy stress tensor, we have a stress tensor induced
by the electromagnetic effects, denoted as $\Btau_m$ \cite{kovetz00}. This
electromagnetic-induced stress tensor is expressed by the Lorentz
force, which furnishes an association between electromagnetism and
the mechanics of the material.
In accordance to \cite{fonteyn10, kovetz00,fonteyn2010fem}, the second-order electromagnetic stress tensor defined as:
\begin{equation}\label{eq13}
	\begin{aligned}
		\Btau_m(\BB;d) = &\underbrace{\frac{E_M}{\mu_0}\Big( \BB \otimes \BB - \frac{1}{2} (\BB \cdot \BB)\textbf{I}\Big) +
			(\BM \cdot \BB)\textbf{I} - \BB \otimes \BM }_{\text{magnetic response}}+\\
		&\underbrace{\epsilon_0\Big( \BE \otimes \BE - \frac{1}{2} (\BE \cdot \BE)\textbf{I}+\BE\times\BB \otimes\dot{\Bx}\Big)}_{\text{electric response}}.
	\end{aligned}
\end{equation}
We note that in the case of polarization effect, the above equation,
i.e., \req{eq13}, should also include this effect, see \cite{kovetz00}
Chapter 15. Here, $E_M>0$ is a constant material property which is
known as the magnetostrictive viscosity constant \cite{miehe2012geometrically}. Since, in the vacuum/air ($\Bx \in\Omega \backslash \calB$) the magnetization effect disappeares, thus results in $\BM =
\textbf{0}$, in the electromagnetic stress  \req{eq13}. 
Additionally, since we only consider here for ferromagnetic material thus, we only take into account magnetic response in the quasi-static state, so the electric response in \req{eq13} is not considered. This assumption is acceptable for many materials which only represent the magnetic forces, e.g., iron; see for more detailed discussion \cite{fonteyn2010fem, fonteyn10}. For further information related to the derivation of the stress expression, which contains the electromagnetically-dependent fields and parameters, you can refer to \cite{dorfmann+ogden03, brigadnov+dorfmann03, dorfmann+ogden+saccomandi04, hanappier+charkaluk+triantafyllidis21}.

Finally, the reduced second-order stress tensor (by removing the
effect of electric response for the quasi-static state) is defined as
\begin{equation}
	\begin{cases}
		\Btau_m(\BB;d) = \frac{E_M}{\mu_0}( \BB \otimes \BB - \frac{1}{2} (\BB \cdot \BB)\textbf{I}) +
		(\BM \cdot \BB)\textbf{I} - \BB \otimes \BM \quad \text{if} \;\Bx \in\calB\\\\  
		\Btau_m(\BB;d) = \frac{E_M}{\mu_0}( \BB \otimes \BB - \frac{1}{2} (\BB \cdot \BB)\textbf{I})  \hspace*{4.2cm} \quad \text{if} \;\Bx \in\Omega \backslash \calB
	\end{cases}.
\end{equation}
%
Therefore, the total stress tensor is given in the solid domain
$\calB$ for a ferromagnetic material as:
\begin{equation}
	{\Btau}(\Bve,\BB, d) = {\Bsigma}(\Bve,\BB,d) + {\Btau}_m(\BB;d).
	\label{eq15}
\end{equation}
By substituting the Cauchy stress tensor and the electromagnetic
stress considering the magnetization effect, the total stress tensor follows  
\begin{equation}
	\begin{aligned}
		{\Btau}(\Bve, \BB, d) = &g_e(d){\widetilde{\Bsigma}}^{+}({\bm{\varepsilon}})+{\widetilde{\Bsigma}}^{-}({\bm{\varepsilon}})+
		\frac{1}{2}\sum_{i=0}^{4}\frac{1}{i+1}\partial_{I_1}g_i\frac{I_4^{i+1}}{B^{2i}_{\text{ref}}}\textbf{I} \\[1ex]
		& + \left( E_M\mu_{0}^{-1} + \sum_{i=0}^{4}\frac{g_i}{i+1}\partial_{I_4}\frac{I_4^{i+1}}{B^{2i}_{\text{ref}}} +
		\frac{1}{2}\alpha_5\right)\BB \otimes \BB \\[1ex]
		&- \left( \frac{E_M}{2}\mu_{0}^{-1} +
		\sum_{i=0}^{4}\frac{g_i}{i+1}\partial_{I_4}\frac{I_4^{i+1}}{B^{2i}_{\text{ref}}}
		\right)\BB \cdot
		\BB \textbf{I} \\[1ex]
		& + \frac{1}{2}\alpha_5\Big(\BB \otimes \BB \cdot \Bve - (\BB \cdot \Bve \cdot \BB)\BI \Big) \\[1ex]
		& + \frac{1}{2}\alpha_6\Big(\BB \otimes \BB + \BB \otimes \BB \cdot \Bve^2 -(\BB \cdot \Bve^2 \cdot \BB)\BI + (\BB \otimes \BB \cdot \Bve + \Bve \cdot \BB \otimes \BB) \Big).
	\end{aligned}
	\label{eq17}
\end{equation}
Additionally the magnetic flux $\BH$ can be specified
\begin{equation}
	\BH( \Bve,\BB;d) =\frac{1}{\mu_0}\BB-2\sum_{i=0}^{4}\frac{g_i}{i+1}\frac{\partial I_4^{i+1}}{\partial I_4} - \frac{1}{2}(\alpha_5\BB \cdot \Bve + \alpha_6\BB \cdot \Bve^2 ).
\end{equation}

\begin{Remark}
	\label{traction_remark}	
	We note the traction force vector $\bar{\Bt}=\Btau \cdot\Bn$ is generally decomposed through the mechanic, electric, and magnetic fields' so-called \textit{generalized traction vector}. If we define $\bar{\Bt}$ as a mechanical traction force vector \textit{only} only mechanical contribution ${\widetilde{\Bsigma}}$ exist, while the traction force vector for the electric and magnetic fields will be neglected, and set as zero at Neumann boundary conditions \cite{wu2021wrinkling}. Additionally, the traction force due to the electric field is discussed in \cite{trimarco2002stresses}, which one can evaluate the electrostatic traction through the so-called \textit{Maxwell stress tensor}. Also, regarding the magnetic force vector, the resultant magnetic force is divided into magnetic traction and body forces \cite{zheng2011magnetic}. Let's define traction vector as the interaction between two adjacent surface elements. Traction magnetic force is computed based on potential magnetic force inside and outside of the magnetizable body. Thus, magnetic traction can be seen as the jumping value of a physical quantity between two sides of the surface \cite{zheng2011magnetic}. So, the Neumman boundary conditions can be set for the magnetization jump between two adjacent magnetized surfaces. The traction force in magnetic separators is also investigated in \cite{henjes1994traction}. Additionally, the electromagnetic traction due to electromagnetic stress tensor is described in \cite{rinaldi2002body}. They discussed whether the electromagnetic effects enter as a body force or as a surface stress, but not both. This is same as our case where we have only body force. 
\end{Remark}

\sectpc{Fracture contribution} 
The fracture contribution of the pseudo-energy density given in \req{psuedo-energy} reads:
\begin{equation}
	\begin{aligned}
		{W}_{frac}(d, \nabla d ):= G_c \gamma_l(d, \nabla d),
	\end{aligned}
	\label{frac-part}
\end{equation}
where ${G}_c > 0$ is so-called a Griffith's energy release rate. In
this study, we use a small time increment $\Delta t =
t_{n+1}- t_n$ where $t_n$ denotes the previous time step. This reduces the complexity of the problem. We further assume that: 
\begin{equation}\label{aniso_permability_assume}
	\begin{aligned}
	&\mu_0(d,t)\approx\mu_0(d,t_n),\;\;
	\epsilon_0(d,t)\approx\epsilon_0(d,t_n),\\\\
	&\sigma_0(d,t)\approx\sigma_0(d,t_n),\;\;
	\rho_0(d,t)\approx\rho_0(d,t_n).
		\end{aligned}
\end{equation}
This type of assumption has already been employed in hydraulic
phase-field fracture. See, for
instance, \cite{mauthe2017variational,miehe2015minimization}, where
the second-order anisotropic permeability tensor is computed of the
previous time $t_n$. Finally, by taking the first variational derivative $\delta_d W$ of \req{psuedo-energy}, the positive crack driving state function $\calD(\Bx,t)$ follows
\begin{equation}
	\calD:=\frac{2l_d}{G_c} \psi_e^{+}(I_1,I_2)
	\ .
	\label{driving-forceD}
\end{equation}

\sectpb[Section21]{Variational formulations for the coupled multi-field problem}
\label{sec:weak_form}
The variational formulations with respect to the four PDEs described in the previous sections for failure mechanics of ferromagnetic materials under the magnetostrictive effects
are further discussed.
The rate-dependent PDEs models, discussed in the previous section, are
defined in a temporal domain $[t_n, t_{n+1}]$ with $\Delta t=t_{n+1}-t_{n} >0$ holds. 
For the time-dependent problem in electromagnetic-induced fracture, we require initial  conditions at $t=0$ as:
\begin{equation}
	\begin{aligned}
		&\Bu^0=\Bu(\Bx,0),  \quad \;d^0=d(\Bx,0),\\
		&\BE^0=\BE(\Bx,0),  \quad \BH^0=\BH(\Bx,0),\\
		&\BD^0=\BD(\Bx,0),  \quad \BB^0=\BB(\Bx,0),
	\end{aligned}
\end{equation}
where $\BH^0$ has to satisfy the divergence -free assumption
\begin{equation}
	\nabla.\left(\mu_0\BH^0(\Bx)\right)=0\quad \text{in}\;\Omega\;\;\;\; \AND \;\;\;\;\BH^0\cdot\BN=0  \quad\text{on}\;\partial \Omega\;.
\end{equation}

\sectpc{Electromagnetic contribution} 
The magnetostatic response is a special case of the electromagnetic
problem which is obtained when the electric field
$\BE(\Bx,t)$ does not change by time, and thus ${\partial \BD} /{\partial t} =0$.
This results in 
\begin{equation}
	\frac{\partial{\BD(\Bx,t)}}{\partial{t}}= \epsilon_0\frac{\partial{\BE(\Bx,t)}}{\partial{t}}=\bm{0} .
\end{equation}
The variation of the electric displacement in time is 0, and thus we can assume that permittivity is zero, i.e., $\epsilon_0=0$. The time-discretized variational formulation of 
the evolution equation of the magnetic vector potential $\BA(\Bx,t)$ in  \req{eq:final_A} is derived by multiplying it
with the test function $\delta{\BA}$ and then applying integrating by part
\begin{equation}
	\begin{aligned}\label{eq_weak_A}
		\mathcal{E}_{\BA}(\BfrakU, \delta \BA) & = 
		\Delta t\int_{\Omega} \frac{1}{\mu_0(d)} \nabla \BA : \nabla (\delta \BA)
		\,\text{d}\Bx 
		 + \int_{\Omega} \sigma_0(d)\left[(\BA-\BA^n)+\Delta t\nabla \Phi \right] \cdot \delta \BA
		\,\text{d}\Bx  \\
		& - \Delta t\int_{\Omega} \BJ_s(t) \cdot \delta \BA \, \text{d}\Bx
		= 0 \qquad \qquad \forall \,\delta\BA\in\calS_0^{\Bv}(\text{curl};\Omega)\;, 
	\end{aligned}
\end{equation}
where \req{eq:final_A} is imposed.


\sectpc{Elastic contribution} 

The variational formulations of the strong formulation of the Euler-Lagrange equations with respect to the $\Bu(\Bx,t)$ in Formulation 2.a reads: 
\begin{equation}
	\mathcal{E}_\Bve(\BfrakU, \delta \Bu) = \displaystyle\int_\calB \Big[ 
	\Btau:\delta \Bve - \bar{\Bb} \cdot \delta \Bu \Big] \text{d}\Bx  - \int_{\partial^{\Bu}_N\calB} \bar{\Bt} \cdot \delta \Bu  \; \text{d}s
	= 0 \quad \forall \,\delta\Bu\in\calW_0^{\Bu}(\text{div};\calB).
	\label{weakFormE}
\end{equation}
We note that $\mathcal{E}_\Bve(\BfrakU, \delta \Bu)$
unlike the magnetic vector potential given in \req{eq_weak_A}, lives in
$\calB$ (and not $\Omega$). Since $\Btau$ is a function of $\BB$ (and
not $\BA$), the local constitutive equation $\BB=\nabla \times \BA$
has to be imposed  after solving magnetic vector potential from \req{eq_weak_A}.
\sectpc{Fracture contribution} 

The variational formulation of the strong formulation of fracture with respect to the crack phase-field $d(\Bx,t)$ in Formulation 2.a reads: 
\begin{equation}
	\begin{aligned}
		\mathcal{E}_d(\BfrakU, \delta d) &= (1-\kappa)\Delta t\displaystyle\int_{\calB} (d-1) \calH\cdot \delta d \;\text{d}\Bx +
		\displaystyle\int_{\calB}\Big[   \Delta t d\cdot\delta d+\eta_d(d-d^n)\cdot\delta d\Big]\;\text{d}\Bx \\ 
		&+\displaystyle\int_{\calB}\Big[  l_d^2 \Delta t\nabla d\cdot\nabla(\delta d)\Big]\,
		\text{d}\Bx  = 0  \quad \quad \forall \,\delta d\in\calW_0^{d}(\text{div};\calB).
		\label{weakFormF}
	\end{aligned}
\end{equation}
As in the weak form $\mathcal{E}_\Bve(\BfrakU, \delta \Bu)$, the crack phase-field equation lives in $\calB$.

Thus, by means of \req{eq_weak_A}, \req{weakFormE}, and
\req{weakFormF}, the fully coupled variational multi-field problem to describing magnetostrictive induced fractures in ferromagnetic material is formulated in the compact form as:
\begin{equation}\label{compat_argmin}
	\begin{aligned}
		\mathcal{E}_{\BfrakU}(\BfrakU, \delta \BfrakU)
		=\mathcal{E}_\Bve(\BfrakU, \delta \Bu)&+\mathcal{E}_A(\BfrakU, \delta \BA)+\mathcal{E}_d(\BfrakU, \delta d)=0\\[2mm]
		&\forall\;(\delta \Bu,\delta \BA,\delta d)\in \Big( \calW_0^{\Bu}(\text{div};\calB), \calS_0^{\Bv}(\text{curl};\Omega),\calW_0^{d}(\text{div};\calB) \Big).
	\end{aligned}
\end{equation}
In order to solve the magnetostrictive induced fracture system \eqref{compat_argmin}, we first solve the first two equations monolithically (simultaneously obtain $(\Bu, \BA)$). Then, a staggered approach is employed to obtain the phase-field fracture $d$. To that end, we fix alternately  $(\bm u,
\BA)$ and estimate $d$ and vice versa. The procedure is continued until
convergence (using a given $\texttt{TOL}_\mathrm{Stag}$). The alternate minimization scheme applied to the \req{compat_argmin} is summarized in Algorithm 1.

\begin{algorithm}
    \small
	{\bf Input:} loading data $\bar{\bm t}_{n}$ on $\partial^{\Bu}_N\calB$ and $\BJ_s(t)$ on wires \\[2mm]
	\hspace{1.4cm}solution $(\bm u^{n-1},\BA^{n-1},d^{n-1})$ from step $n-1$.
	\\[2mm]
	\quad\quad Initialization, $k=1$:\\
	
	\quad\quad \textbullet\; set $(\bm u^0,\BA^0,d^0):=(\bm u^{n-1},\BA^{n-1},d^{n-1})$.\\
	
	Staggered iteration between $(\bm u,\BA)$ and $d$:\\
	
	\quad\quad \textbullet~ solve following system of equations in a monolithic manner by given $d^{k-1}$,
	\begin{align*}
		\begin{cases}
			\Delta t\int_{\Omega} \frac{1}{\mu_0(d)} \nabla \BA : \nabla (\delta \BA)
			\,\text{d}\Bx 
			+ \int_{\Omega} \sigma_0(d)\left[(\BA-\BA^n)+\Delta t\nabla \Phi \right] \cdot \delta \BA
			\,\text{d}\Bx   - \Delta t\int_{\Omega} \BJ_s(t) \cdot \delta \BA \, \text{d}\Bx
			= 0 ,\\\\
			\displaystyle\int_\calB \Big[
			\Btau(\Bve,\BB,d):\delta \Bve - \bar{\Bb} \cdot \delta \Bu \Big] \text{d}\Bx  - \int_{\partial^{\Bu}_N\calB} \bar{\Bt} \cdot \delta \Bu  \; \text{d}s
			= 0\quad
			with \; \;\BB=\nabla \times \BA,
		\end{cases}
	\end{align*}\\
	
	\quad\quad \quad for $(\bm u,\BA)$, set $(\bm u,\BA)=:(\bm u^k,\BA^k)$,\\
	
	\quad\quad \textbullet\; given $(\bm u^k, \BB^k)$, solve
	\begin{align*}
		(1-\kappa)\Delta t\displaystyle\int_{\calB} (d-1) \calH\cdot \delta d \;\text{d}\Bx +
		\displaystyle\int_{\calB}\Big[   \Delta t d\cdot\delta d+\eta_d(d-d^n)\cdot\delta d\Big]\;\text{d}\Bx
		+\displaystyle\int_{\calB}\Big[  l_d^2 \Delta t\nabla d\cdot\nabla(\delta d)\Big]\,
		\text{d}\Bx  = 0\\
	\end{align*}
	
	\quad\quad \quad for $d$, set $d=:d^k$,\\
	
	\quad\quad \textbullet\; for the obtained pair $(\bm u^k,\BA^k,d^k)$, check {staggered residual} by
	\begin{align*}
		\quad \quad \;\; \mathrm{Res}_\mathrm{Stag}^k:|\mathcal{E}_{\Bve}(\BfrakU^k;\delta \Bu)|
		+|\mathcal{E}_{\BA}(\BfrakU^k;\delta \BA)|
		\leq\texttt{TOL}_\mathrm{Stag}
	\end{align*}
	\quad\quad \textbullet\; if fulfilled, set $\BfrakU^k=(\bm u^k,\BA^k,d^k)=:(\bm u^n,\BA^n,d^n)=\BfrakU^n$ then stop; \\
	
	\quad\quad\; {\color{white}\textbullet} else $k+1\rightarrow k$. \\[2mm]
	{\bf Output:} solution $(\bm u^n,\BA^n,d^n)$ at $n^{\text{th}}$ time-step. \\[2mm]
	\caption{\em The staggered iterative solution process for \req{compat_argmin} at a fixed time-step $n$.}
	\label{alg_upd}
\end{algorithm}
\begin{Remark}
	\label{BMD}	
	Since the permeability of the ferromagnetic material is
	computed as a function of the crack phase-field with a
	degrading function, the magnetic vector potential $\BA$, the
	magnetic field $\BB$ and the magnetic flux $\BH$ are
	implicitly determined as a function of the crack
	phase-field. Therefore, when a crack starts to propagate in the material, the degradation function actively represents its degrading effect on the magnetic response of the material. Then, the magnetic vector potential $\BA$, the magnetic field $\BB$ and the magnetic flux $\BH$ decrease along the crack paths.
\end{Remark}

\newpage
\sectpb[Section29]{Space finite element discretization}

In this section, we provide spatial discretizations of the variational
forms given in Algorithm \ref{alg_upd}. This results in a discretized
multi-field problem to be solved for three-field unknowns represented
by $(\Bu, \BA, d)$. In the case of the magnetostatic problem, one can
use either N\'ed\'elec elements or alternatively $H^1(\Omega)$
conforming elements with continuous piecewise polynomials,
see \cite{li2019computational}. We use a Galerkin finite element
method to discretize the equations, employing second order
isoparametric curl-conforming triangular element $P_2$ for the ansatz
and test spaces of all primary fields see \cite{garcia02, jamelot05,
ciarlet+jamelot07, asadzadeh+beilina23,
meunier2010finite,cardoso2016electromagnetics}. The continuous solid
domains $\calB$, and the vacuum domain $\Omega$ are approximated by
$\calB_h$, and $\Omega_h$ such that $\calB\approx\calB_h$, and
$\Omega\approx\Omega_h$. The approximated domains $\calB_h$, and
$\Omega_h$ are decomposed with non-overlapping linear triangular finite element $\calB_e\subset\calB_h$ such that 
\begin{equation}
	\calB\approx\calB_h = \bigcup_{e}^{n_{e}}\calB_e \AND
	\Omega\approx\Omega_h = \bigcup_{e}^{n_{e}}\Omega_e.
\end{equation}
The unknown fields $(\Bu, \BA, d)$ are approximated by
\begin{equation}
	\Bu_h=\BN_\Bu\hat{\bm u}, \quad
	d_h=\bm{N}_d\hat{\bm d}, \quad
	\BA_h=\bm N_\BA\hat{\BA}.
	\label{Discr1}
\end{equation}
with a set of nodal solution values $(\hat{\bm u},\hat{\bm
d},\hat{\BA})$. Their derivatives are given by
\begin{equation}
	\bm\varepsilon(\bm u_h)=\bm B_u\hat{\bm u}, \quad
	\nabla d_h=\bm B_d\hat{\bm d}\;,
	\quad
	\nabla \BA_h=\bm B_\BA\hat{\bm A}\;,
	\label{Discr2}
\end{equation}
which represent constitutive state variables $(\Bve_{\Bu_h}, \nabla
d_h,\nabla \BA_h)$. Here, $\BB_{\Bu}$, $\BB_{d}$, $\BB_{\BA}$ are the
matrix representation of the shape function's derivatives,
corresponding to the global deformation, crack phase-field, and the
magnetic vector potential,  respectively. Additionally, the magnetic
flux densities (which is enter to deformation equation) are
approximated by the curl operator as:
\begin{equation}
	\BB(\BA_h)=\nabla \times \BA_h=\bm{\Pi}\hat{\BA}.
	\label{Discr3}
\end{equation}
Here, $\bm{\Pi}$ is the matrix representation of the shape function's
derivative, see  \cite{fonteyn10} (Section 3.2),
and \cite{bastos2003electromagnetic} for a detailed discussion. For
each primary field, the set of discretized weak forms which is based on the \textit{residual} force vector denoted by $\BR_{\bullet}$ has
to be determined. First, the mechanical weak form leads to
\begin{equation}
	{{\bf{R}}_\Bu} (\Bu,d,\BA)=
	\sum_{e} \int_{\calB_e}\big[ \BB^T_\Bu \Btau(\Bve_h,\BB_h,d_h)-  \overline{\Bb}\cdot\BN_\Bu\big]\,\text{d}\Bx - \int_{\partial\calB_N}  \overline{\Bt}\cdot\BN_\Bu\,ds ={\bm 0}, 
	\label{phf_fem_u}
\end{equation}
next, for the magnetic vector potential  yields,
\begin{equation}
		\begin{aligned}
	{{\bf{R}}_\BA}(\Bu,d,\BA)=
	\sum_{e}\int_{\Omega_e} \bigg[ \frac{\Delta t}{\mu_0(d)} \BB^T_\BA\nabla \BA 
	&+  \sigma_0(d)\big[(\BA-\BA^n)+\Delta t\nabla \Phi \big] \cdot \BN_{\BA}\\ 
	& - \Delta t \BJ_s(t) \cdot \BN_{\BA} \; \bigg]\,\text{d}\Bx ={\bm 0}\;,
	\label{phf_fem_phi}
		\end{aligned}
\end{equation}
and lastly the crack phase-field follows as:
\begin{equation}
	\begin{aligned}
	{{\bf{R}}_d}(\Bu,d,\BA)=
	\sum_{e}\int_{\calB_e} \bigg[ \bigg((1-\kappa)(d_h-1)\calH+  {d} &+ \frac{\eta_d}{\Delta t}(d_h-d_h^{n})\bigg)\cdot\BN_d\\
	&+   l_d^2\BB_d^T\cdot\nabla d_h \bigg]\,\text{d}\Bx ={\bm 0}. 
	\label{phf_fem_d}
	\end{aligned}
\end{equation}
%
\sectpa[Section5]{Numerical examples}
\label{Section5}
In this section, the capability of the proposed framework in predicting the magneto\-restrictively-induced cracking in the ferromagnetic materials is approved by demonstrating the numerical examples. 
At first, the evolution of the magnetic vector potential and the
magnetic field in a ferromagnetic material is investigated. Then,
cracking in the ferromagnetic materials is investigated by adopting
the magneto-mechanical model coupled with the phase-field
approach. The boundary-value problem of the first example is related
to a domain that includes copper wires and an iron beam that are
surrounded by vacuum. In the remaining examples, the domain of the
boundary-value problems contains predefined notches and copper wires
where a constant electric current source is enforced. These
representative numerical examples are solved by developing a finite
element software based on the FEniCS software library, see \cite{alnaes+etal15, logg+etal12}. The material parameters used in the numerical examples are listed in Table~\ref{material-parameters}. For the staggered approach, see Algorithm 1, we set $\texttt{TOL}_\mathrm{Stag}=10^{-4}$. We also note that the second term of anisotropic material constants in  \req{aniso_permability} is assumed for quasi-impermeable crack faces, hence we set $(\mu^f_0,\epsilon^f_0, \sigma^f_0,\rho^f_0)$ identical to vacuum counterparts. 

\begin{table}[!ht]
	\caption{Material parameters used in the numerical examples based on \cite{kovetz00}}	\vspace{2mm}
	\centering
	\begin{tabular}{ccllll}
		No.  &Parameter & Name                   & Example 1-2  & Example 3  & Unit            \\\hline 
		1.   &$E$        & Young's modulus       & $160$    & $16$ & $\mathrm{GPa}$ \\
		2.   &$\nu$      & Poisson's ratio       & $0.33$  & $0.2$  & --                 \\
		3.   &$\mu_{0,\Omega}$        &   Permeability of vacuum    & $4\pi\times10^{-7}$  & $--$  & $\mathrm{H/mm}$ \\
		4.   &$\mu_{0,\calW}$   &  Permeability of wires        & $1.26\times10^2$  & $1.26\times10^3$   & $\mathrm{H/mm}$ \\
		5.   &$\mu_{0,\calB}$        &   Permeability in solid  & $1000$  & $10^7$ & $\mathrm{H/mm}$ \\
		6.   &$\sigma_0$        &   Electromagnetic conductivity   & $1$  & $1$ & $\mathrm{S/{mm}}$ \\
		7.   &$B_{\text{ref}}$&  The magnetic flux density  & $2 \times 10^{6} $ & $2 \times 10^{6} $ & $\mathrm{T}$ \\
		8.   &$(\alpha_0, ... ,\alpha_3)$&  Magnetisation constants  & $2 \times 10^{-6}$  & $2 \times 10^{-6}$ & -- \\
		9.   &$(\alpha_4,\alpha_5)$&  Magnetostrictive constants  & $(2.5,7.75)$  & $(2.5,7.75)$ & -- \\
		10.   &$E_M$& Magnetostrictive viscosity  & $10^{-9} $ & $10^{-9}$ & -- \\
		11.   &$\zeta$& Permeability transition exponent  & $50 $ & $50 $ & -- \\
		12.   &$\eta_d$  & damage viscosity   & $1 \times 10^{-6}$ & $1 \times 10^{-6}$ & $\mathrm{N/(mm^2)s}$ \\
		13.   &$l_d=2h_e$  & fracture regularization parameter  & $2.196$  & $2.196$& -- \\
		14.   &$\kappa$  & Stabilization parameter   & $10^{-8}$ & $10^{-8}$ & -- \\
		15.   &$G_c$     & Griffith's energy release rate   & $0.0027$ & $0.0027$ &\\
		16.   &$(m,c_1,c_2)$  & degradation parameters   & $(2,0.1,0.9)$ & $(2,0.1,0.9)$ & -- \\  \hline
		\label{material-parameters}
	\end{tabular}
\end{table}

\sectpb[Section51]{Example 1: Magnetostatic problem for transient magnetic vector potential}
The first insight into the performance of the variational form of
Maxwell's equations in a finite element context are gained by a two-dimensional transient magnetostatic problem in a domain composed of copper
wires and an iron cylinder that are
surrounded by a vacuum. The geometrical configuration of the 2D test
problem and its finite element model are depicted in Figure \ref{example1_BVP}, in
which the subdomains for the iron cylinder, the copper wires, and the
surrounding vacuum with the radius of $r_v=5$ mm are
clearly visible.
The inner radius of the cylinder is $r_1=1$ mm
and its thickness is $n=0.2$ mm. The domain contains ten copper
wires with a radius of $r_w=0.1$ mm. The discretization of the domain is performed with
$64652$ triangular elements with the minimum element size of
$h_{\text{min}} \approx 0.04$ mm.

The wires are assumed to carry the electric current, which is varied linearly
with time. The electromagnetic waves are neglected, since the problem is supposed to be magnetostatics. To examine this problem, the Poisson-type equation \req{eq5} is
solved. The specific current $\BJ_s(t)$ along the $z$-direction is
given by
\begin{equation*}
	\begin{aligned}
	 \texttt{exterior wires:}  \;\; & \BJ_s(t,z)=-t \;\; \text{A/mm$^2$},\quad\texttt{interior wires:}  \;\; & \BJ_s(t,z)= t\;\; \text{A/mm$^2$},
	\end{aligned}		  
\end{equation*}
for the interior set of the winding copper wires.
The results are obtained in terms of the magnetic vector
potential in $z$-direction, $A_z(\Bx,t)$, for the different time
steps, and the magnetic flux
density $\BB(\Bx,t)$ for the last time step. These results are shown in Figure \ref{example1_A} and Figure \ref{example1_B}, respectively. 

\begin{figure}[!ht]
	\centering
	{\includegraphics[clip,trim=0cm 11.6cm 0cm 13cm, width=16cm]{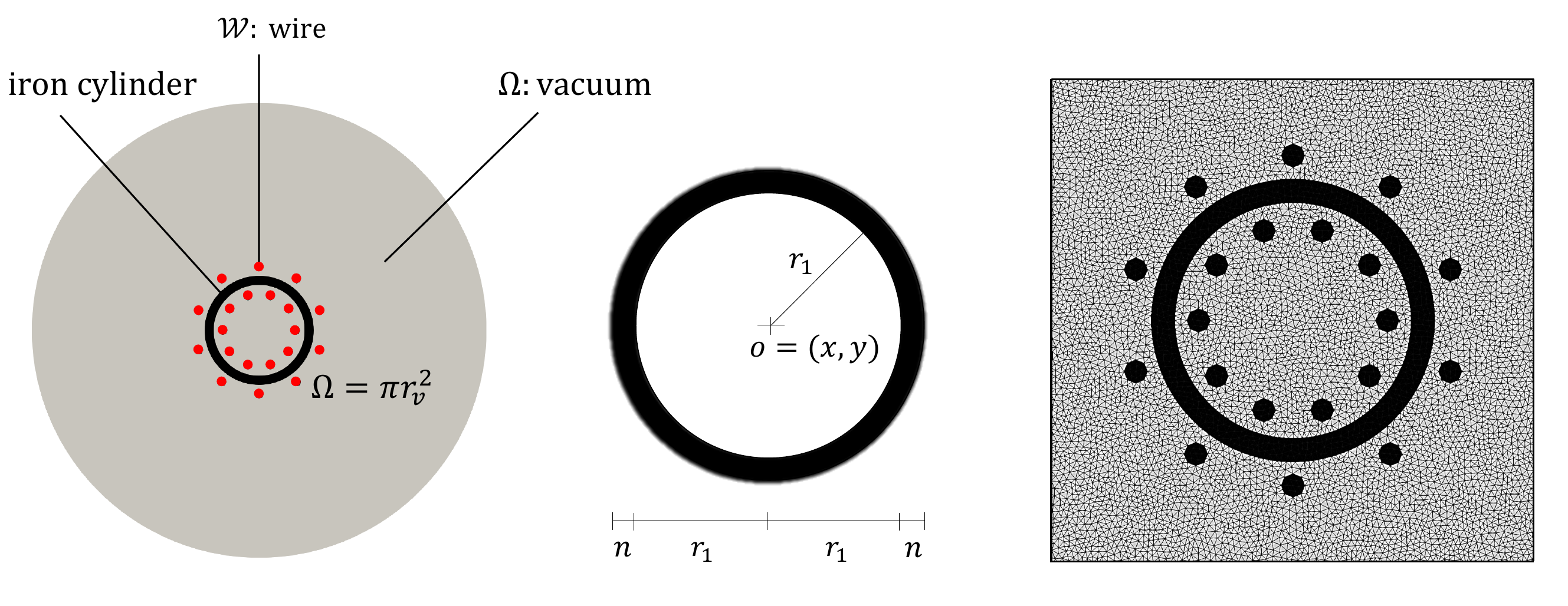}}  
	\vspace{-0.8cm}
	\caption*{\hspace*{1cm}(a)\hspace*{4.6cm}{(b)}\hspace*{4.9cm}{(c)}\hspace*{0.9cm}}
	\caption{Example 1. The representation of (a) whole domain of
    the magnetostatic problem including copper wires winding around
    the iron cylinder and vacuum, (b) geometry and dimensions of iron cylinder, and (c) finite element mesh generated for the test problem. The
    subdomains for the iron cylinder and copper wires are clearly
    demonstrated.  
	}
	\label{example1_BVP}
\end{figure}

\begin{figure}[!ht]
	\centering
	{\includegraphics[clip,trim=1cm 31.5cm 0cm 0cm, width=16.5cm]{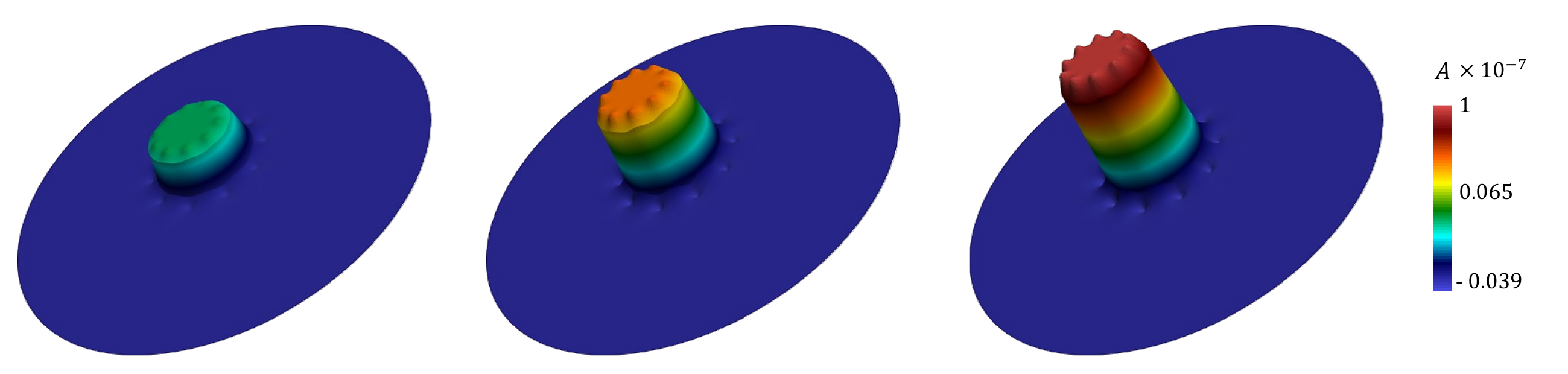}}  
	\vspace{-0.7cm}
	\caption*{\hspace*{1cm}(a)\hspace*{4.6cm}{(b)}\hspace*{4.9cm}{(c)}\hspace*{0.9cm}}
	\caption{Example 1. The representation of the $z$-component of
    the magnetic vector potential $A_z$ at
    (a) $t=40$ s
    (b) $t=80$ s
    (c) $t=120$ s.
	}
	\label{example1_A}
\end{figure}

\begin{figure}[!ht]
	\centering
	{\includegraphics[clip,trim=2cm 24.4cm 0cm 0cm, width=17cm]{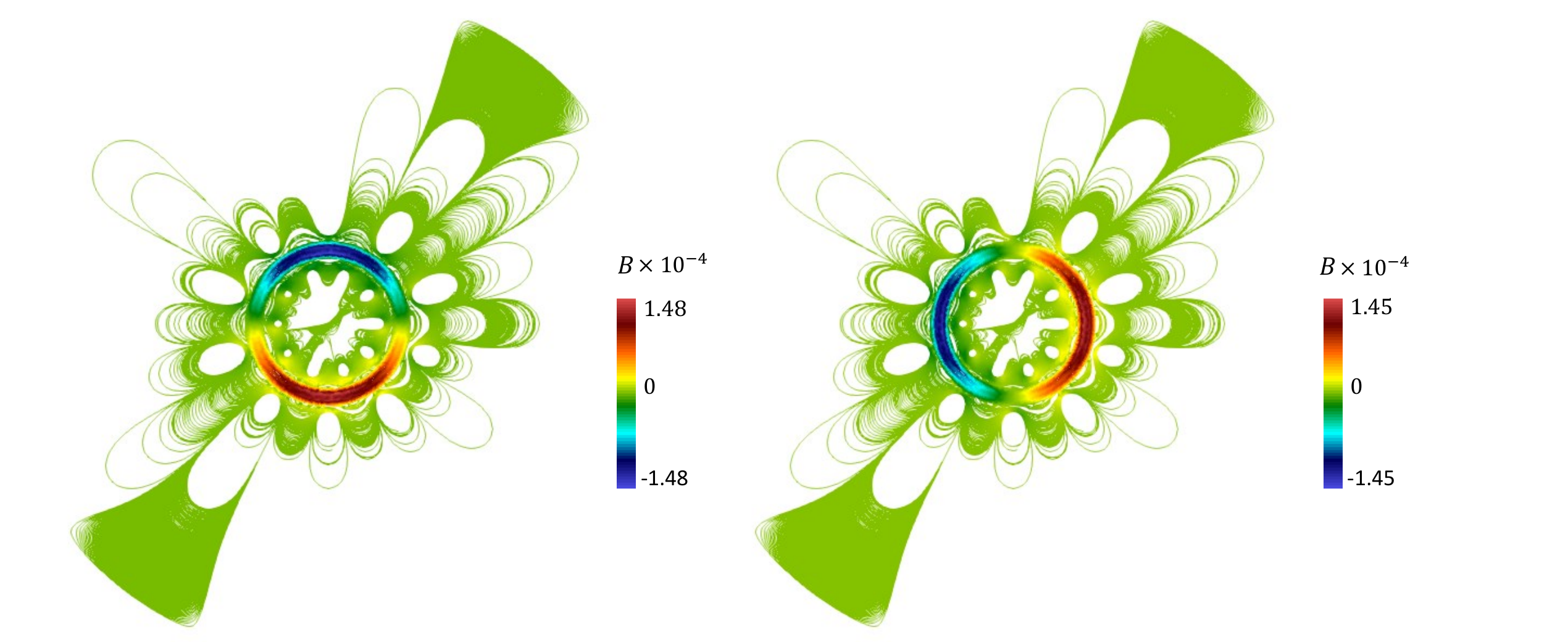}}  
				\vspace{-0.6cm}
	\caption*{\hspace*{0.2cm}(a)\hspace*{7.6cm}(b)\hspace*{1.1cm}}
	\caption{Example 1. The representation of (a) the $x-$ and
    (b) $y-$ components of the magnetic field
    field $\BB$ at $t=120$ s
	}
	\label{example1_B}
\end{figure}

\sectpb[Section53]{Example 2: Magneto-mechanically induced cracking in an
iron beam surrounded by vacuum}
The second example investigates the magnetostrictive-induced cracking in a
ferromagnetic notched beam. This problem is solved by
coupling the Poisson-type transient magnetic equation with the phase-field
model. The domain of the boundary value problem consists of vacuum, two
copper wires, and the single notched simply-supported iron beam. The electric current flows through the top and the bottom wires.
The geometrical configuration and finite element model of the
boundary value problem are displayed in Figure \ref{example3_BVP},
respectively. The whole domain is composed of vacuum with the
radius of $r_v=8$ mm, two copper wires with the radius of
$r_w=0.5$ mm, and the notched iron beam with $n=0.8$ mm and $m=3.6$ mm located at the middle of
the vacuum with the dimension of $A=4$ mm and $B=2$ mm hence $\calB=(0,4)\times(0,2)$ mm$^2$.
The whole domain
is discretized with $28325$ triangular elements. As the domain is
divided into three subdomains, and the crack phase-field equation is
only solved in the iron beam subdomain, the discretization is refined in that
region. In this case, the minimum element size is
$h_{\text{min}}\approx 0.005$ mm and the maximum element size is
$h_{\text{max}}\approx 1.098$ mm, see Figure \ref{example3_BVP}(c).

\begin{figure}[!ht]
	\centering
	{\includegraphics[clip,trim=0cm 11.5cm 0cm 13cm, width=16cm]{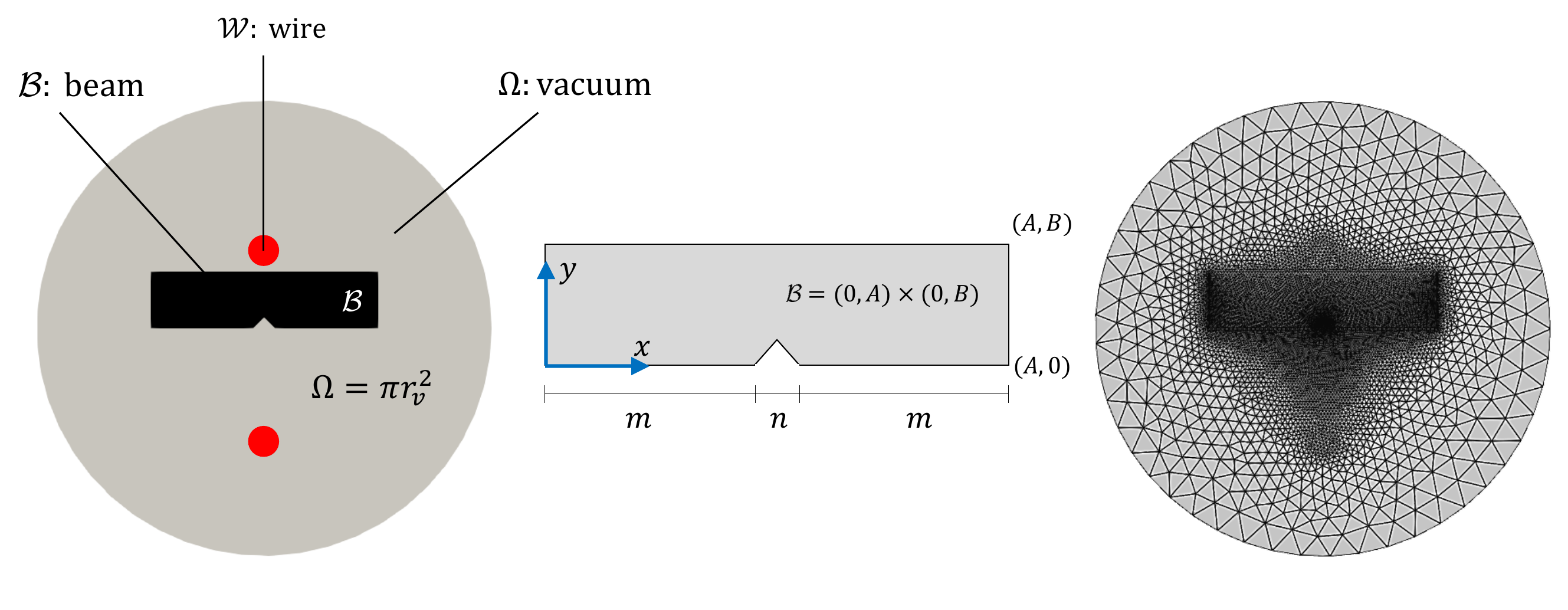}}  
		\vspace{-0.7cm}
	\caption*{\hspace*{1cm}(a)\hspace*{4.6cm}{(b)}\hspace*{4.9cm}{(c)}\hspace*{0.9cm}}
	\caption{Example 2. The representation of geometries and dimensions of (a) whole domain, (b) iron beam, and (c) finite element model of BVP. The subdomains are clearly invisible.
	}
	\label{example3_BVP}
\end{figure}

The wires located
on the top and bottom of the iron beam contain the electric
current. The specific current $\BJ_s(t)$ along the $z$-direction is
given in the form of a linear equation which varies in time as
wire. 
\begin{equation*}
\begin{aligned}
\texttt{top wire:}  \;\; & \BJ_s(t,z)=2.5\,t \;\; \text{A/mm$^2$}, \quad \texttt{bottom wire:}  \;\; & \BJ_s(t,z)= -1.5\,t\;\; \text{A/mm$^2$}.
\end{aligned}		  
\end{equation*}

At first, the
evolution of magnetic vector potential $A_{z}(\Bx,t)$ and magnetic
field $\BB(\Bx,t)$ are calculated through Maxwell's equations in the whole domain. 
The stress response regarding the magnetostrictive effects is determined in the iron beam. 
The electric current flow in the wires generates magnetization in the
notched iron beam that forces it to deform. The deformation causes
stress development around the notch. Therefore, cracking starts to
initiate when the maximum principal stress exceeds the critical value.
The crack initiation and propagation in the beam at different time
steps are demonstrated in Figure \ref{example3_res_d}. The magnetic vector potential $A_z(\Bx,t)$ and the magnetic field $\BB(\Bx,t)$ at the last time step $t=0.68$ s are exhibited in Figure \ref{example3_res_AB}. 

\begin{figure}[!tbh]
	\centering
	{\includegraphics[clip,trim=0.5cm 8cm 1cm 6.5cm, width=17cm]{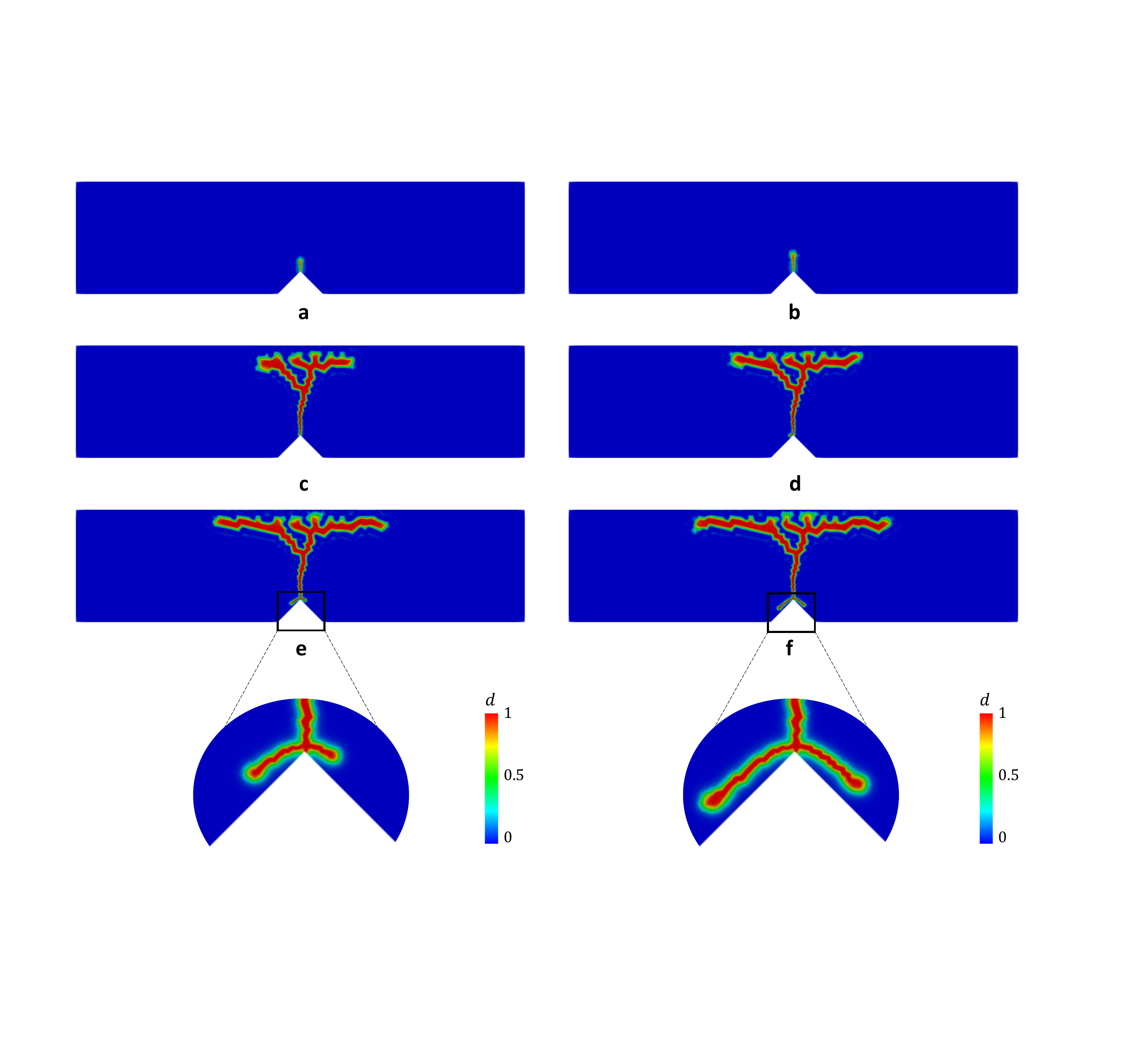}}   
	\caption{Example 2. The representation of the crack initiation and
		propagation in the iron beam under coupling
		electromagneto-mechanical effects at
		(a) $t=0.43$ s,
		(b) $t=0.45$ s,
		(c) $t=0.56$ s,
		(d) $t=0.59$ s,
		(e) $t=0.64$ s and
		(f) $t=0.68$ s.
	}
	\label{example3_res_d}
\end{figure}

\begin{figure}[!tbh]
	\centering
	{\includegraphics[clip,trim=1.7cm 26.3cm 5cm 1.2cm, width=15cm]{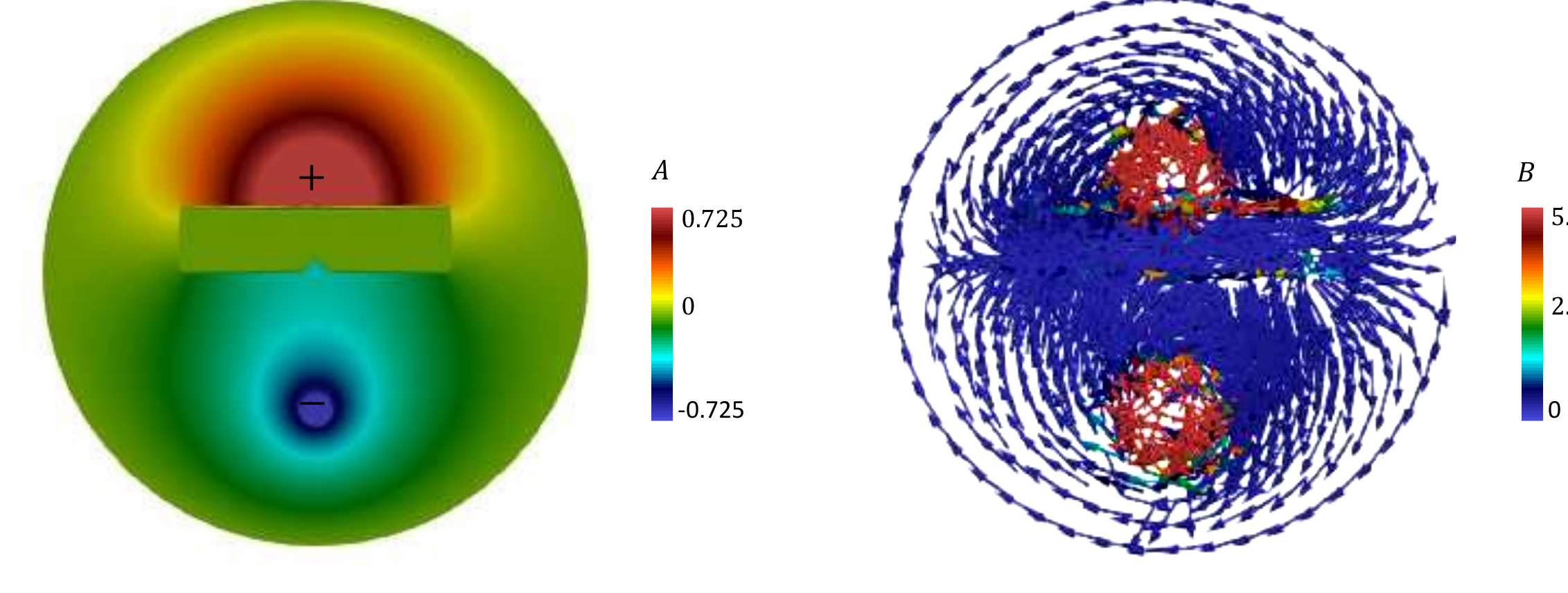}}   
			\vspace{-0.3cm}
	\caption*{\hspace*{0.2cm}(a)\hspace*{7.6cm}(b)\hspace*{1.1cm}}
	\caption{Example 2. The representation of (a) the magnetic vector potential $A_z$, and (b) the magnitude of the magnetic field $\BB$ at $t=0.68$ s (final failure state).
	}
	\label{example3_res_AB}
\end{figure}

\sectpb[Section54]{Example 3: Magneto-mechanically induced cracking in a
  ferromagnetic material containing predefined notches and wires}
This example comprises four sub-examples in which the
magneto-mechanically driven fracture in a ferromagnetic material is examined. In the first two examples, the electric current is enforced through
the notches, and in the last two examples, the electric current flows
through the copper wires. The
cracking response of the ferromagnetic material is examined by solving the coupled magneto-mechanically driven crack phase-field problem. The
following sub-examples investigate a
square plate of ferromagnetic material. The geometrical
configuration of the first two examples with the predefined notches is
shown in Figure \ref{example4_BVP1}. We set $A=40$ mm thus $\calB=(0,80)^2$ mm$^2$.
The discretization of the domain is performed with $51200$ triangular 
elements with an element size of $h \approx 0.707$ mm. In the last
two sub-examples, the square plates contain copper wires with the
radius of $r_w=0.3$ mm, see Figure \ref{example4_BVP2}. The third
sub-example with four copper wires is discretized with $31132$
triangular elements. In the domain, the uniform element size is $h \approx 1.02$ mm. The fourth sub-example with five copper wires is discretized with $30984$ triangular elements by the uniform element size as $h \approx 1.02$ mm. In these BVPs, all displacements
are fixed in both directions, and the electric vector potential is set
to zero at the boundaries
$\partial_{D}\cal{B}$. The notches and copper wires carry the constant electric
source current of $\bar{J}=0.03$ A/mm$^2$.
%
\sectpc{Sub-example 1: Transversely wired plate with three winding wires}
In this example, the domain
contains three predefined notches $\calC_1$, $\calC_2$ and
$\calC_3$ of length $n=8$.
The notch $\calC_1$  is located horizontally
in the center of the domain. The other vertical notches ($\calC_2$ and
$\calC_3$) are located with distances of $m=20$ mm from the right and
left boundaries, see Figure \ref{example4_BVP1}(a). The results in
terms of the crack phase-field, the magnetic vector potential $A_z$,
and von Mises stress are exhibited in
Figure \ref{example4_case1}. Moreover, the magnetic field $\BB$
developed in the plate is shown in Figure \ref{example4_res_M1}(a). The
variation of the magnetic vector potential average value with respect
to time it includes the effect of damage
is represented in Figure \ref{example4_res_LD1}.a. 

\sectpc{Sub-example 2: Transversely wired plate with nine winding wires}
In this example, the square plate considered in the BVP
contains nine predefined notches $\calC_i$ with $=1,2,3,...,9$. The
dimensions of the notches are the same as the previous sub-example, see
Figure \ref{example4_BVP1}(b). The results in terms of the crack
phase-field, magnetic vector potential, and von Mises stress are
presented in Figure \ref{example4_case2} at different time steps. The
magnetic field arising in the plate under the constant electric
current source is shown at the last time step in
Figure \ref{example4_res_M1}.b. The average value of the magnetic
vector potential over time in the square plate regarding the cracking
effect is depicted in Figure \ref{example4_res_LD1}(b).  
  
\begin{figure}[!t]
	\centering
	{\includegraphics[clip,trim=1cm 12.2cm 5cm 8cm, width=15cm]{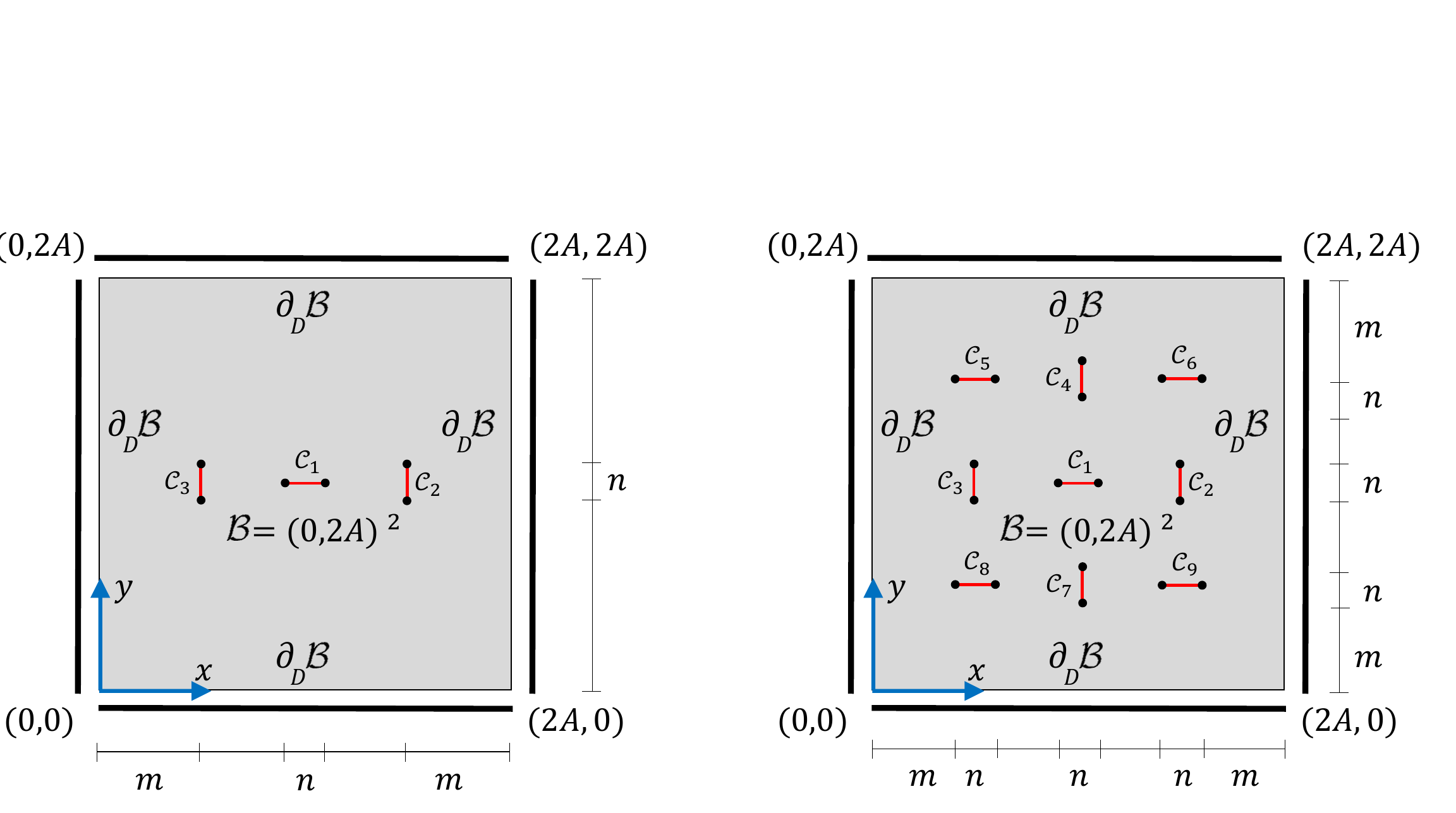}}  
			\vspace{-0.3cm}
	\caption*{\hspace*{0.1cm}(a)\hspace*{7.7cm}(b)\hspace*{0.3cm}}
	\caption{Transversely wired plates. The representation of geometry, the corresponding dimensions and the boundary conditions of the square plate contains (a) three predefined notches (Example 3.1), and (b) nine predefined notches (Example 3.2)  
	}
	\label{example4_BVP1}
\end{figure}
In both sub-examples, it is observed that cracks start to evolve in
the square plate when a constant electric
current is enforced in the predefined notches. As the crack propagates
over time, the electric vector potential increases along the
crack. However, its value decreases when it does not reach the
cracking state. The evolution of the average value of the
magnetic vector potential over time confirms that when the plate is
completely damaged, the value decreases. The comparison of the results
indicates that the average value of the magnetic vector potential in
the plate with nine notches reaches its maximum value earlier than in the plate with three notches.
\sectpc{Sub-example 3: Longitudinally wired plate with four winding wires}
In the third case, the square plate contains four copper
wires, $\calC_i$ with $=1,2,3,4$. The radius of the wires is $r_{\text{wire}}=3$ mm.
The dimensions of the plate are the same as the previous
sub-examples, see Figure \ref{example4_BVP1}(a). The wires are placed
of a distance of $m=20$ mm from the edge of the plate. The crack initiation and propagation, the evolution of the
magnetic vector potential, and the von Mises stress for the different time
steps are presented in
Figure \ref{example4_case2}. The magnetic field under the constant
electric current is presented
in Figure \ref{example4_res_M1}(a) for the last time step. The average
value of the magnetic vector potential over time in the square plate
regarding the cracking effect
is provided in Figure \ref{example4_res_LD2}.a.  

\begin{figure}[!tbh]
	\centering
	{\includegraphics[clip,trim=0.8cm 12.2cm 5cm 14cm, width=15cm]{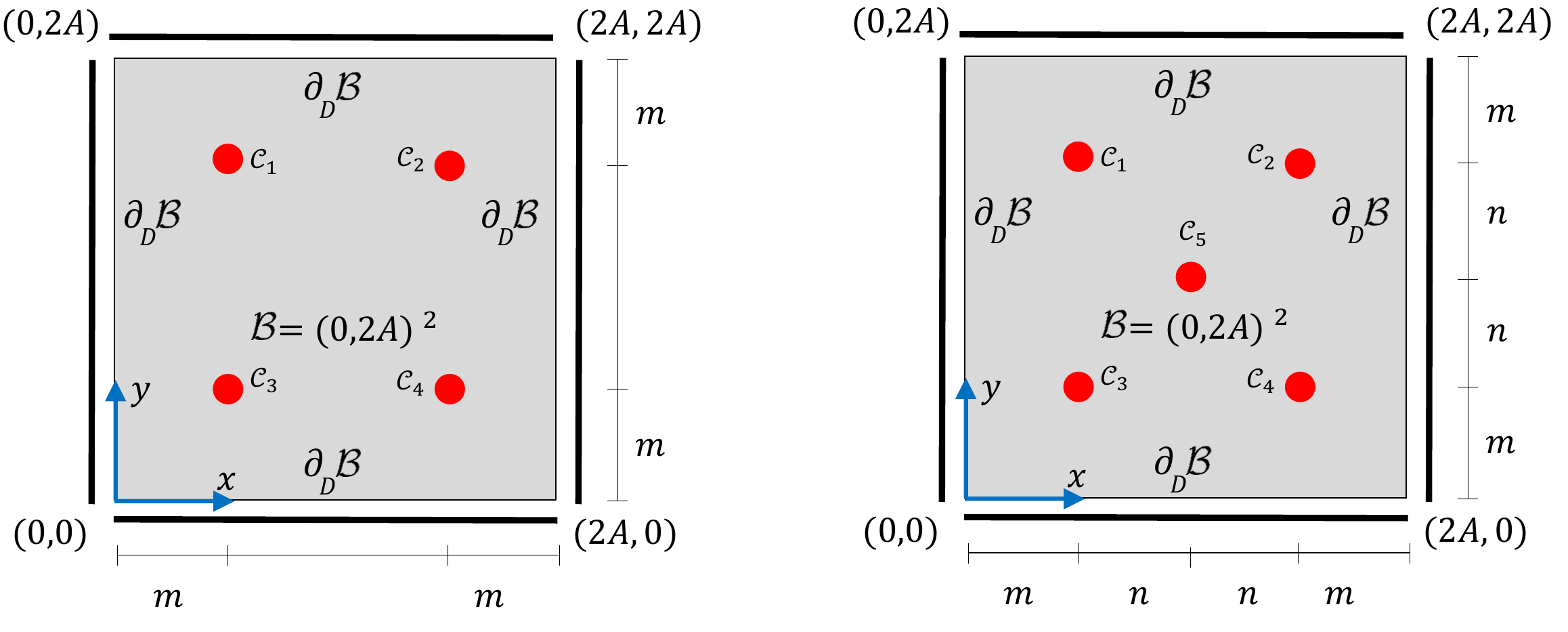}}  
		\vspace{-0.3cm}
	\caption*{\hspace*{0.1cm}(a)\hspace*{7.7cm}(b)\hspace*{0.3cm}}
	\caption{Longitudinally wired plates. The representation of geometry and the corresponding dimensions of the square plate contains (a) four copper wires (Example 3.3), and (b) five copper wires (Example 3.4).
	}
	\label{example4_BVP2}
\end{figure}

\begin{figure}[!ht]
	\centering
	{\includegraphics[clip,trim=0.5cm 6cm 3.2cm 3cm, width=16cm]{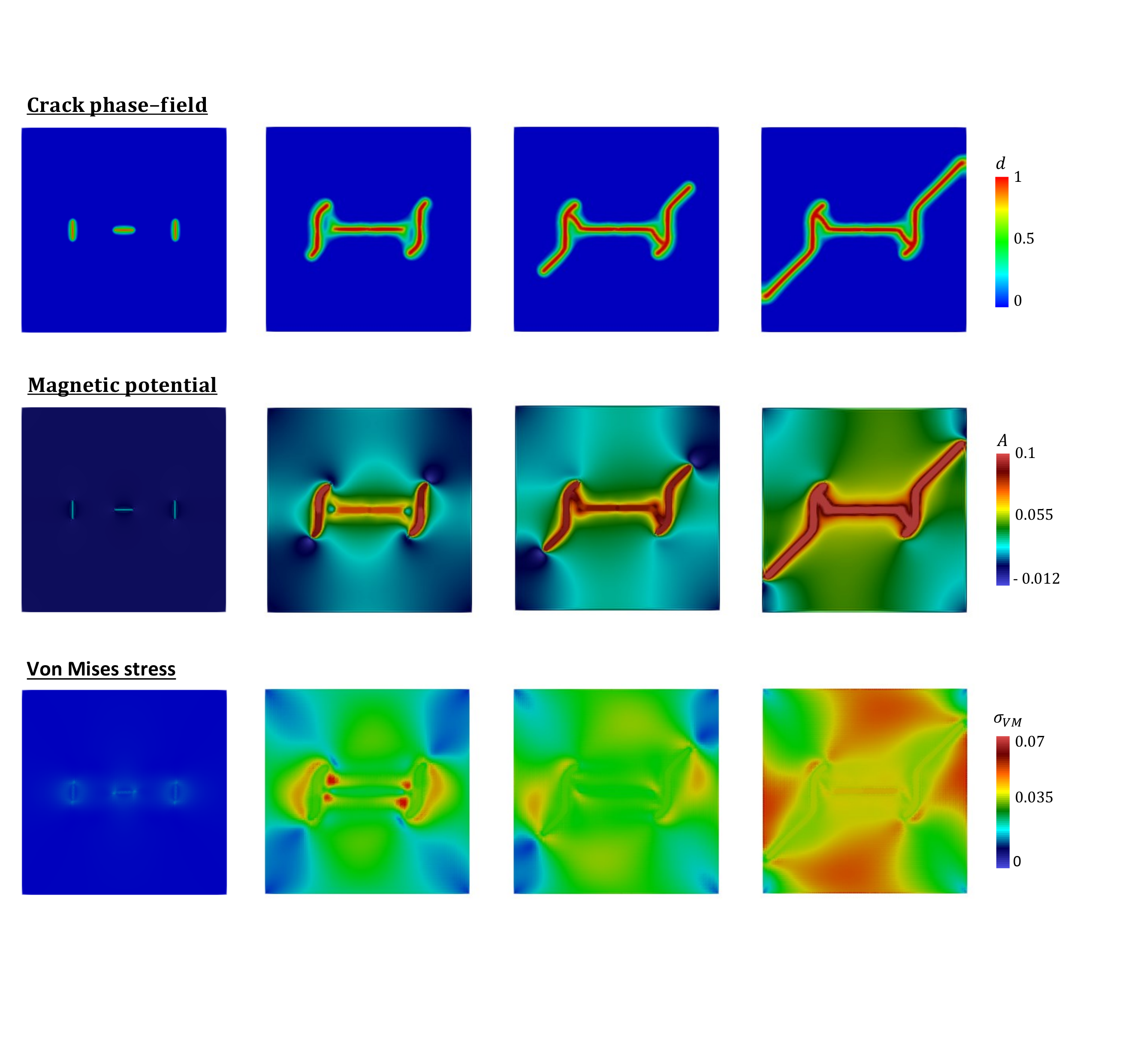}}  
	\caption{Example 3.1. Transversely wired plate with three winding wires. Reference results of the magnetostrictively induced crack driven by constant electric current source induction along the predefined notches. Evolution of the crack phase-field, magnetic vector potential and von Mises stress at different time steps. 
	}
	\label{example4_case1}
\end{figure}

\begin{figure}[!ht]
	\centering
	{\includegraphics[clip,trim=0.5cm 6cm 3.2cm 3cm, width=16cm]{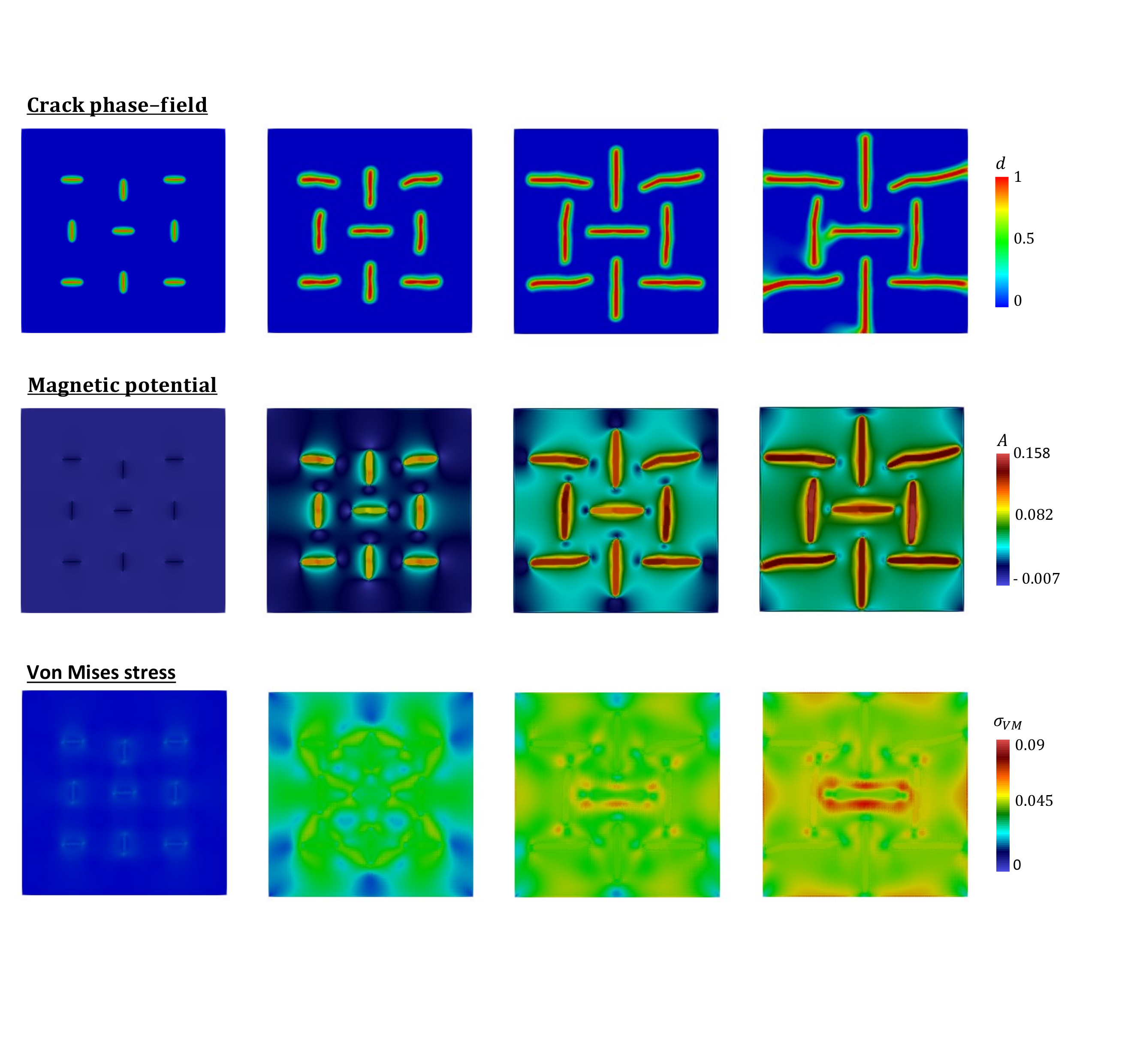}}  
	\caption{Example 3.2. Transversely wired plate with nine winding wires. Reference results of the magnetostrictively induced crack driven by constant electric current source induction along the predefined notches. Evolution of the crack phase-field, magnetic vector potential and von Mises stress at different time steps. 
	}
	\label{example4_case2}
\end{figure}

\begin{figure}[!ht]
	\centering
	{\includegraphics[clip,trim=0.5cm 6cm 3.2cm 3cm, width=16cm]{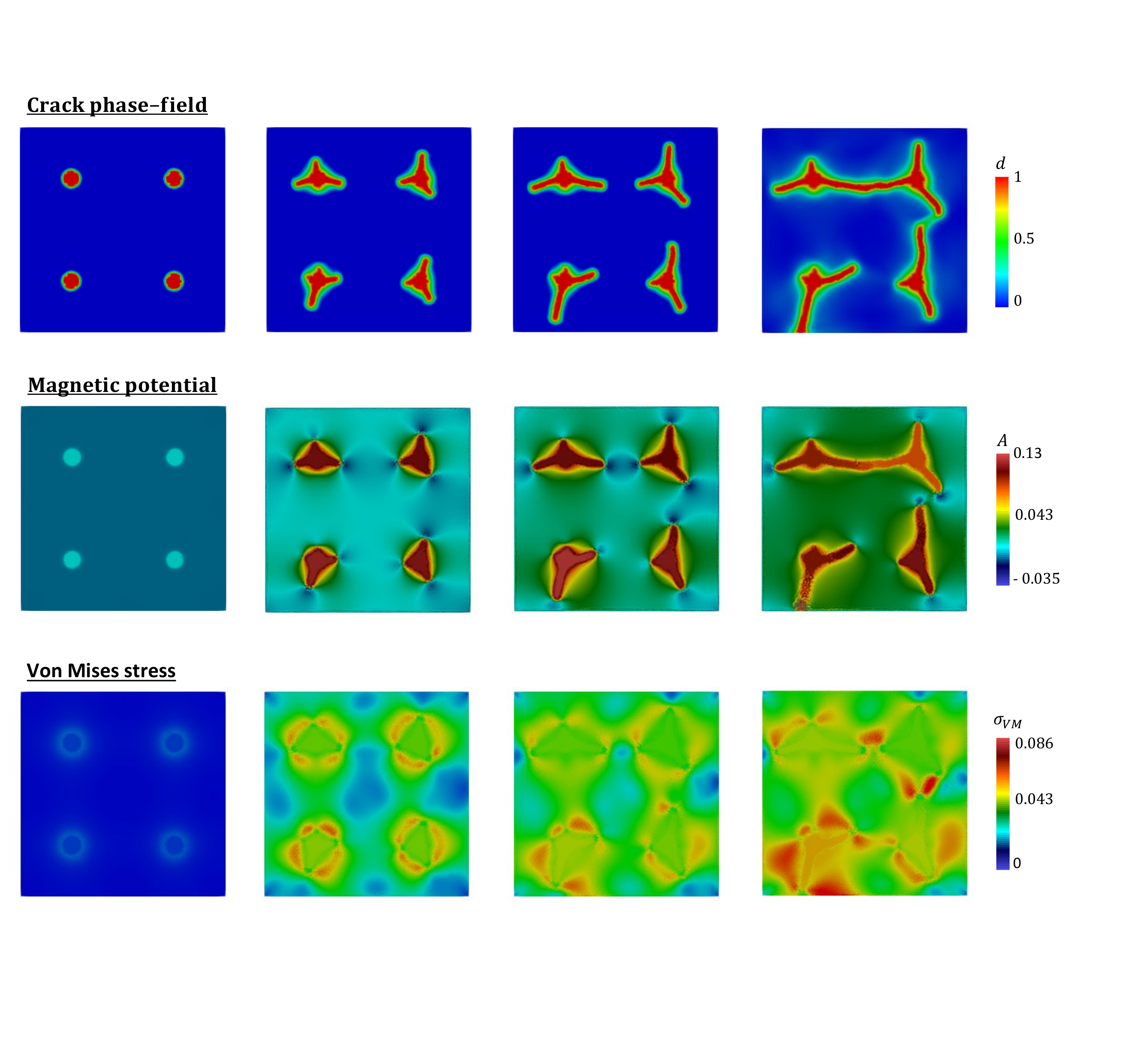}}  
	\caption{Example 3.3. Longitudinally wired plate with four winding wires. Reference results of the magnetostrictively induced crack driven by constant electric current source induction on the copper wires. Evolution of the crack phase-field, magnetic vector potential and von Mises stress at different time steps.	
	}
	\label{example4_case3}
\end{figure}

\sectpc{Sub-example 4: Longitudinally wired plate with five winding wires}
In the last example, the square plate contains five copper wires,
$\calC_i$ with $=1,2,3,4,5$. The geometries of the plate and wires and
their dimensions are as in the previous example. In addition to
Sub-example 3, this one has an additional copper wire located at the center of the plate with $n=20$ mm, see Figure \ref{example4_BVP2}(b). The results in terms of the crack phase-field, magnetic vector potential, and 
von Mises stress are depicted in Figure \ref{example4_case4}. The additional findings, including the magnetic field at the last time step and the average value of the magnetic vector potential, are shown in Figure \ref{example4_res_M1}(b) and Figure \ref{example4_res_LD2}(b), respectively.

\begin{figure}[!tbh]
	\centering
	{\includegraphics[clip,trim=0.5cm 6cm 3.2cm 3cm, width=16cm]{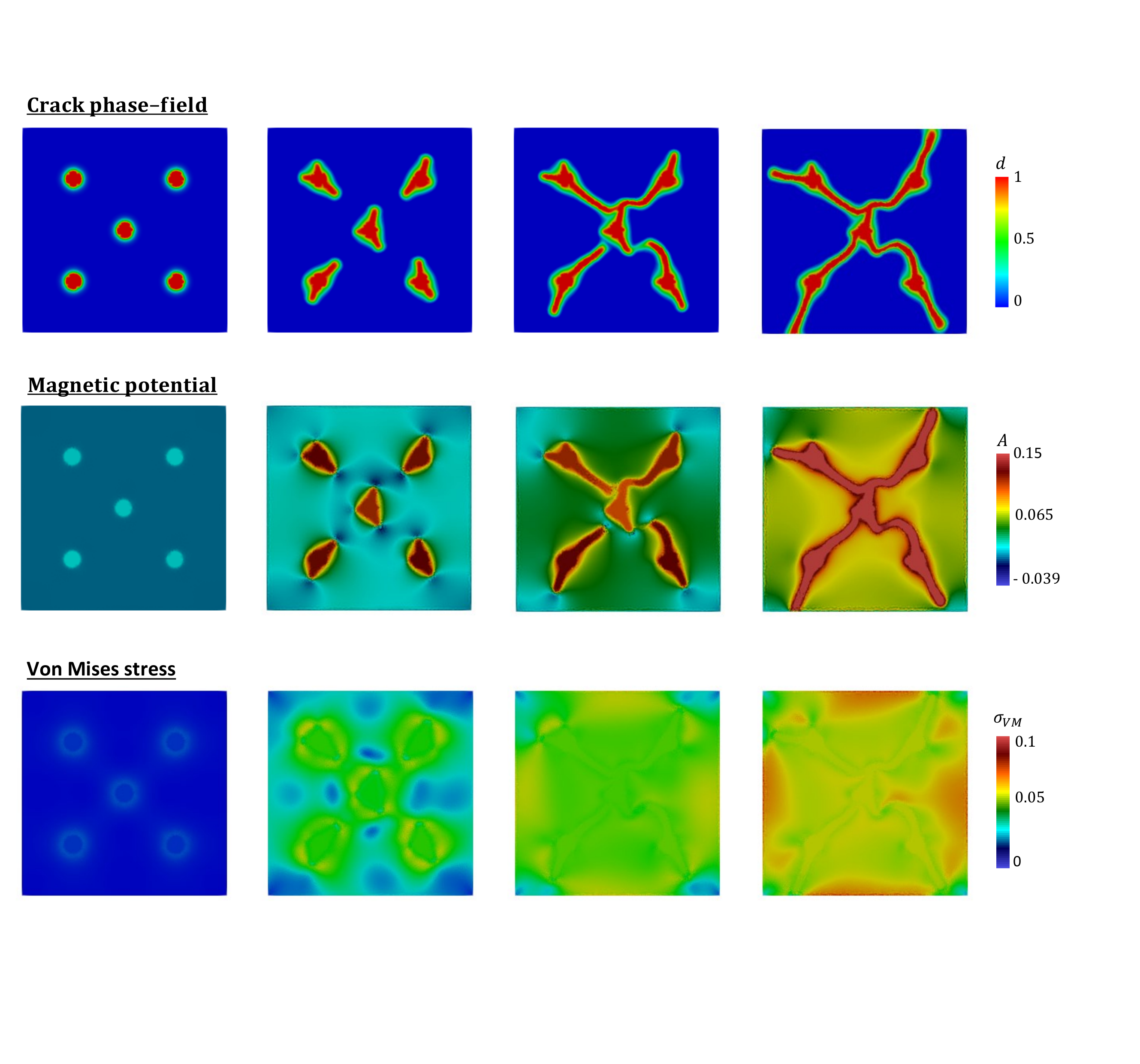}}  
	\caption{Example 3.4. Longitudinally wired plate with five winding wires. Reference results of the magnetostrictively induced crack driven by constant electric current source induction on the copper wires. Evolution of the crack phase-field, magnetic vector potential and von Mises stress at different time steps.
	}
	\label{example4_case4}
\end{figure}

\begin{figure}[!ht]
	\centering
	{\includegraphics[clip,trim=3cm 14.8cm 4.2cm 11cm, width=16cm]{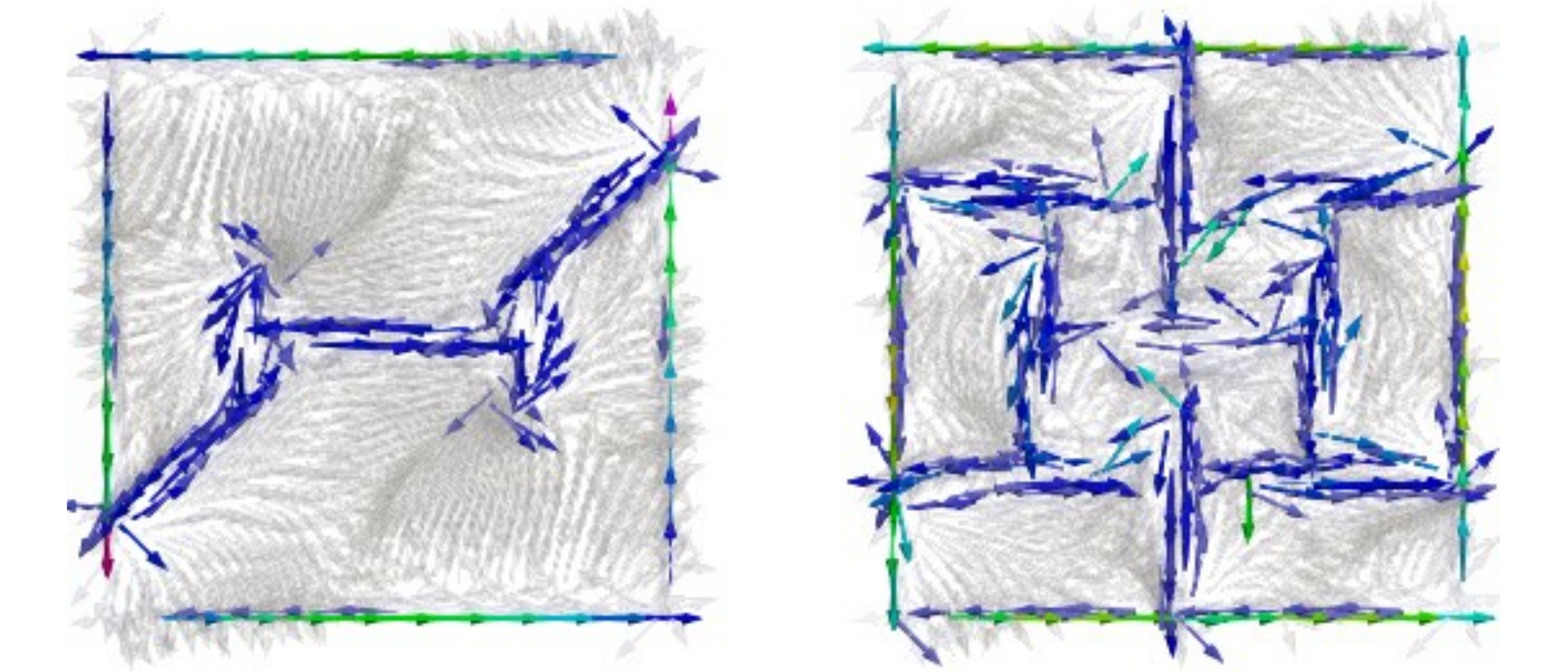}}  
	\vspace{-0.8cm}
	\caption*{\hspace*{1.3cm}(a)\hspace*{7.7cm}(b)\hspace*{0.3cm}}
	\caption{Transversely wired plates. The representation of magnetic field induced by the constant electric current source along the predefined notches on (a) square plate with three notches (Example 3.1), and (b) square plates with nine notches (Example 3.2).
	}
	\label{example4_res_M1}
\end{figure}

\begin{figure}[!ht]
	\centering
	{\includegraphics[clip,trim=3cm 13.8cm 4.2cm 11.5cm, width=15.8cm]{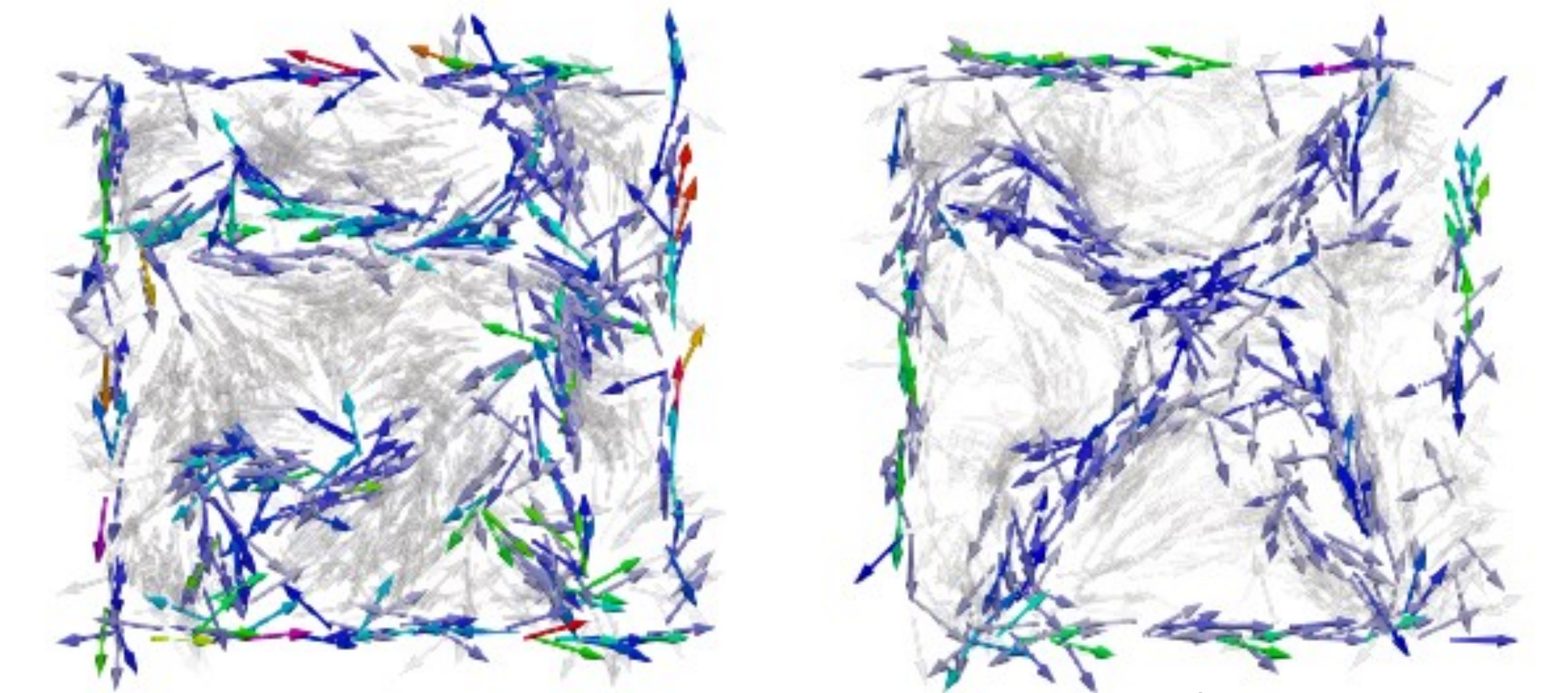}}  
	\vspace{-0.2cm}
	\caption*{\hspace*{1.3cm}(a)\hspace*{7.7cm}(b)\hspace*{0.3cm}}
	\caption{Longitudinally wired plates. The representation of magnetic field induced by the constant electric current source along the copper wires on square plate with (a) four copper wires (Example 3.3), and (b) five copper wires (Example 3.4).
	}
	\label{example4_res_M2}
\end{figure}

\begin{figure}[!t]
	\centering
	{\includegraphics[clip,trim=2cm 21.5cm 2cm 8cm, width=16cm]{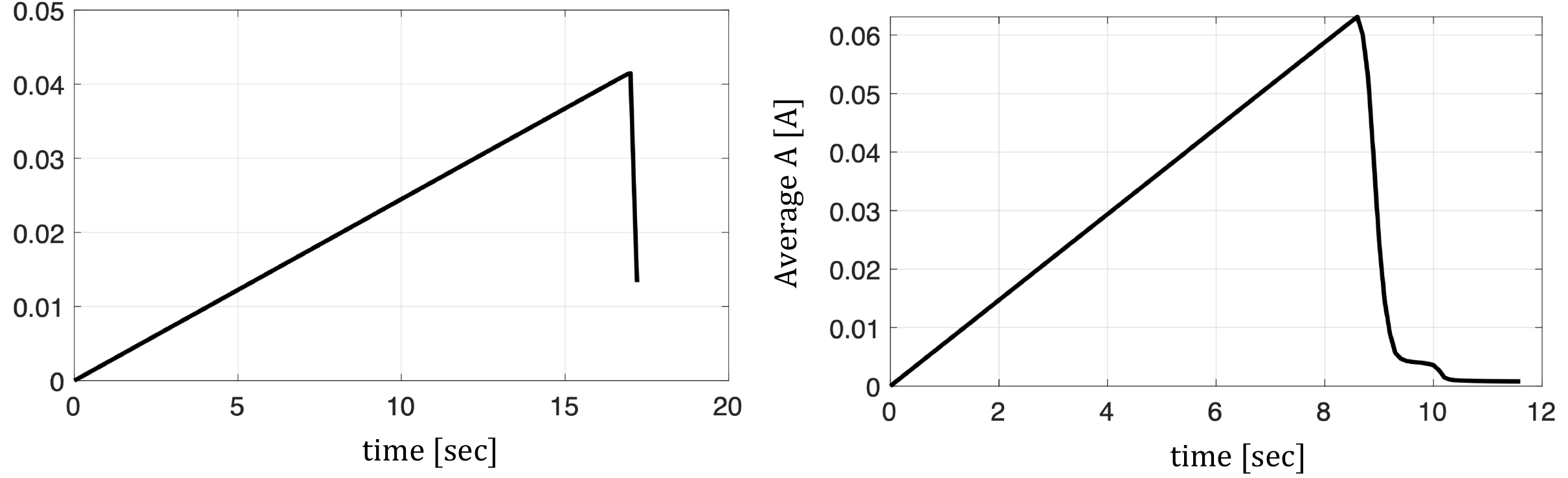}}  
				\vspace{-0.8cm}
	\caption*{\hspace*{1.3cm}(a)\hspace*{7.7cm}(b)\hspace*{0.3cm}}
	\caption{Longitudinally wired plates. The average value of magnetic vector potential $A_{\text{ave}}$ over time for the square plate contains 
	(a) three notches (Example 3.1), and (b) nine notches (Example 3.2).
}
	\label{example4_res_LD1}
\end{figure}

\begin{figure}[!t]
	\centering
	{\includegraphics[clip,trim=2.5cm 21.5cm 2cm 8cm, width=16cm]{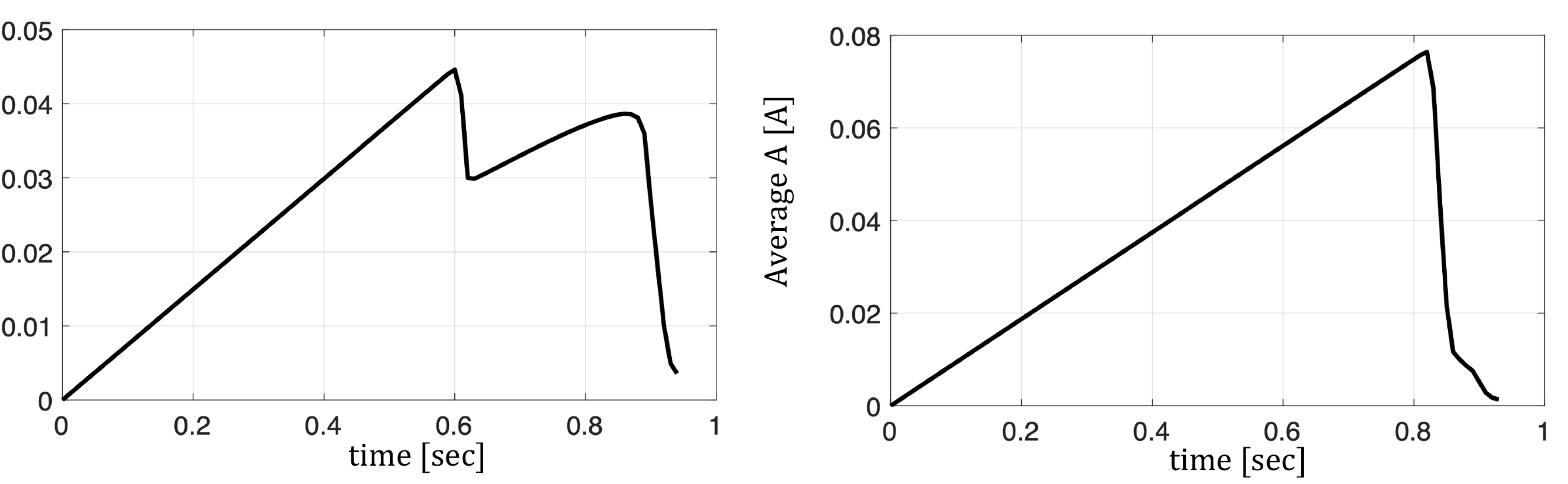}}  
			\vspace{-0.8cm}
	\caption*{\hspace*{1.3cm}(a)\hspace*{7.7cm}(b)\hspace*{0.3cm}}
	\caption{Transversely wired plates. The average value of magnetic vector potential $A_{\text{ave}}$ over time for the square plate contains 
	(a) four copper wires (Example 3.3), and (b) five copper wires (Example 3.4).}
	\label{example4_res_LD2}
\end{figure}

The results obtained for the last two examples show that induction of
a constant electric current in the copper wires causes mechanical
deformation and corresponding stresses in the plate. When the maximum principal stress produced under the magnetostrictive effects exceeds the critical stress value of a ferromagnetic material, the plate starts to degrade. The evolution of the magnetic vector potential confirms that when the material reaches the cracking state, its average value decreases. Furthermore, the average value of the magnetic vector potential in the plate with five wires reaches its maximum value earlier than in the plate with four wires. 

\noii{At this point it is necessary to remark, a further extension to three-dimensional setting is computationally time demanding. This is mainly due to the non-convexity of the energy functional to be minimized with respect to the displacement and the phase-field and non-linearity which calls for the resolution of small length scales. As such, an idea of using adaptive mesh refinement  \cite{yue2023efficient,noii2021quasi}, or alternatively using global/local approach \cite{noii2020adaptive} seems particularly appealing.}  
\cleardoublepage
\sectpa[Section7]{Conclusion}
A coupled electro-magneto-mechanical model, along with the phase-field
approach is developed for simulating crack growth in ferromagnetic material. The proposed coupled electro-magneto-mechanical model evaluates the evolution of the stress response and deformation of the material under the magnetostrictive effects. The magneto-mechanically-driven cracking is examined by applying the phase-field approach. To this end, we extended Maxwell's equation to a
variational-based electro-magneto-elastic model, then coupled it with the phase-field approach. In this model, the total stress is additively decoupled into three parts, containing purely mechanical,
purely-magnetic, and the coupled magneto-mechanical stresses. The magneto-driven
deformation in a solid body is computed as a function of the magnetic
field, which is implicitly defined through
the magnetic vector potential. 

The cracking response of the ferromagnetic material is computationally
formulated on the basis of the degradation of the effective stress tensor.
The degradation function is
formulated as a function of the damage variable. Its initiation
and propagation is developed through the thermodynamically consistent
variational form of the crack phase field. The total energy functional
includes the constitutive energy density functions that correspond to elastic, magnetostrictive, magnetization, and fracture contributions. In the phase-field method, the crack driving force is determined in terms of the mechanical free energy function. In the current work, the capability of the proposed model in predicting the
response of the ferromagnetic material is validated by a couple of representative numerical examples. At first, we used Maxwell's equations to investigate the variation of the electric
and magnetic fields in a solid body surrounded by a
vacuum. Then, we extended the coupled magneto-elastic model along with
the crack phase-field method to examine the cracking response under
the mechanical deformation of a ferromagnetic material induced by magnetization.  

Several topics for further research emerge from the present study. First, the assumption of monotonic increasing for specific current density $\BJ_s$ could be extended to time-dependent cyclic magnetization. So, the possibility of fatigue failure in the ferromagnetic materials with magnetostrictive effects could be
elaborated. For this purpose, the degradation of Griffith's energy release
rate can be considered. Moreover, thermal effects, which can
contribute as an additional source of failure in the ferromagnetic
material, can be taken into consideration. In that case, the
interaction between magnetic and thermal fields needs to be considered
(e.g., in an electromagnetic rail launcher). In addition, the
assumption of magnetostatics could be further relaxed (requiring the
use of N\'ed\'elec elements), such that the variation of the electric
displacement in time does not vanish anymore and hence causes non-zero permittivity in the electric equation.

\subsection*{Acknowledgment}
N. Noii  founded by the Priority Program \texttt{DFG-SPP 2020} within its second funding phase. P. Wriggers were funded by the Deutsche Forschungsgemeinschaft (DFG, German Research Foundation) under Germany's Excellence Strategy within the Cluster of Excellence PhoenixD, \texttt{EXC 2122} (project number: 390833453). 


\begin{Appendix}
	\setcounter{equation}{0}
	\renewcommand{\theequation}{A.\arabic{equation}}
\subsection*{Appendix A. Approximation of the five material constants in ${W}_{mag}$}

The objective is to derive the five material constants denoted as $g_i(\Bx)$ in ${W}_{mag}$ of \req{mag-part}. Recall the total form of the Helmholtz free energy function for the
description of the electromagneto-mechanical response of a
ferromagnetic material is introduced by
\begin{equation}
	W(\BfrakC):= 
		{W}_{elas}(\Bve ,d ) + {W}_{mag}(\Bve,\BB;d)+ {W}_{mos}(\Bve,\BB;d). 
\label{app_eqA10}
\end{equation}
and bases on five invariants 
%
\begin{equation}
	W(\BfrakC):= 
		{W}_{elas}(I_1,I_2,d ) + {W}_{mag}(I_1,I_4;d)+ {W}_{mos}(I_5,I_6;d). 
\label{app_eq1}
\end{equation}
%

This expression is based on invariants which express the
mechanics, the magnetization and the magnetostriction of the
material. The first two of
these invariants describe the isotropic characteristics of the solid
material as a function of the total strain tensor as follows
\eb
I_1(\Bve) = \tr(\Bve) \AND I_2(\Bve) = \frac{1}{2}(\Bve:\Bve).
\label{app_eq2}
\ee

The remaining invariants $I_4, I_5$ and $I_6$ describe
the single-valued magnetization
and the magnetostrictive curves as a function of the magnetic flux
density $\BB$. They are provided in the following expressions
\eb
I_4(\BB) = \BB \cdot \BB\,, \quad
I_5(\Bve, \BB) = \BB \cdot \Bve \cdot \BB \AND I_6(\Bve, \BB) = \BB\cdot \Bve^2 \cdot \BB.
\label{app_eq3}
\ee

The Cauchy stress tensor is derived by taking the partial derivative of the
total Helmholtz free energy function with respect to the total strain
tensor

\begin{equation}
   \Bsigma(\Bve,\BB,d)
	:= \frac{{\partial W(\Bve,\BB,d)}}{\partial \Bve}
	= \displaystyle\sum_{i} \frac{{\partial W(\Bve,\BB, d)}}{\partial I_i}\frac{{\partial I_i}}{\partial \Bve}.
	\label{app_eq4}
\end{equation}


Also, the magnitization, see (65), is determined by taking the partial derivative of the
total Helmholtz free energy function with respect to the magnetic flux density
\begin{equation}
	{\BM}(\Bve, \BB;d)
:= - \frac{{\partial W(\Bve,\BB,d)}}{\partial \BB}
=-\displaystyle\sum_{i} \frac{{\partial W(\Bve,\BB,d)}}{\partial I_i}\frac{{\partial I_i}}{\partial \BB}.
\label{app_eq5}
\end{equation}

In \req{app_eq4} and \req{app_eq5}, the partial derivatives of the
basic invariants with respect to the strain tensor and the magnetic
flux density vector are defined as

\eb
\begin{aligned}
	&\frac{\partial{I_1}}{\Bve} = \BI\, \quad \frac{\partial{I_2}}{\Bve}
	= \Bve\,,
	\quad \frac{\partial{I_5}}{\Bve} = \Bve \otimes \Bve,\\[1ex]
	& \frac{\partial{I_6}}{\Bve} = \BB \otimes \BB \cdot \Bve +
	\Bve \cdot \BB \otimes \BB, \\[1ex]
	&\frac{\partial{I_4}}{\BB} = 2\BB\,, \quad \frac{\partial{I_5}}{\BB} =
	\BB \cdot \Bve \AND  \frac{\partial{I_6}}{\BB} = \BB \cdot \Bve^2.
\end{aligned}
\label{app_eq6}
\ee
Therefore, the Cauchy stress tensor in follows with
\begin{equation}
\begin{aligned}
	{\Bsigma}(\Bve, \BB,d) &:=  \frac{{\partial W(\Bve, \BB,d)}}{\partial \Bve}\\[1ex]
	&=  g_e(d){\widetilde{\Bsigma}}^{+}({\bm{\varepsilon}})+{\widetilde{\Bsigma}}^{-}({\bm{\varepsilon}})
	+\frac{1}{2}\sum_{i=0}^{4}\frac{1}{i+1}\partial_{I_1}g_i
		\frac{I_4^{i+1}}{B^{2i}_{\text{ref}}}\textbf{I}\\[1ex]
	&+ \frac{g(d)}{2}\alpha_5
		\BB \otimes \BB + \frac{g(d)}{2}\alpha_6\left( \BB \otimes \BB \cdot \Bve
		+ \Bve \cdot \BB \otimes \BB \right),
\end{aligned}
\label{eq12A}
\end{equation}
see (72) and (73), as well. The constitutive equation for pure mechanics state is additively split to purely tensile contribution  ${\widetilde{\Bsigma}}^{+}({\bm{\varepsilon}})$ and purely compression contribution ${\widetilde{\Bsigma}}^{-}({\bm{\varepsilon}})$, reads
\begin{equation}
	{\bm \sigma}(\bm{\varepsilon},s ):=\frac{\partial {W}_{elas}(\Bve,d)}{\partial {\bm \varepsilon}} = g_e(d)\frac{\partial {\psi}^{+}}{\partial \bm \varepsilon}+\frac{\partial {\psi}^{-}}{\partial \bm \varepsilon}
	=g_e(d){\bm {\widetilde{\sigma}}^{+}}+{\bm {\widetilde{\sigma}}^{-}},
	\label{eq232_7}
\end{equation}
with
\begin{equation}\label{eq232_8}
	{\widetilde{\Bsigma}}^{+}({\bm{\varepsilon}})=K_nH{^+}(I_1)(\bm \varepsilon:\textbf{I})\textbf{I}+2\mu{\bm \varepsilon}^{dev},\quad \text{and} \quad {\widetilde{\Bsigma}^{-}}({\bm{\varepsilon}})=K_n\big(1-H{^+}(I_1)\big)(\bm \varepsilon:\textbf{I})\textbf{I}.
\end{equation}
Correspondingly magnetization vector denoted as $\BM$ yields:
\begin{equation}\label{eq16}
	\begin{aligned}
		{\BM}(\Bve, \BB;d)
		&:= - \frac{{\partial W(\Bve, \BB,d)}}{\partial \BB}\\
		&=-\sum_{i=0}^{4}\frac{g_i}{i+1}\partial_{I_4}\frac{I_4^{i+1}}{B^{2i}_{\text{ref}}}
		-g(d)\Big(\frac{1}{2}\alpha_5\BB
		\cdot \Bve -\frac{1}{2}\alpha_6\BB \cdot \Bve^{2}\Big),
	\end{aligned}
\end{equation}
see (65), as well.
%
Thus, the total stress tensor is derived as
%
\begin{equation}
	\begin{aligned}
		{\Btau}(\Bve, \BB,d) = &g_e(d){\widetilde{\Bsigma}}^{+}({\bm{\varepsilon}})+{\widetilde{\Bsigma}}^{-}({\bm{\varepsilon}})+
		\frac{1}{2}\sum_{i=0}^{4}\frac{1}{i+1}\partial_{I_1}g_i\frac{I_4^{i+1}}{B^{2i}_{\text{ref}}}\textbf{I} \\[1ex]
		& + \left( E_M\mu_{0}^{-1} + \sum_{i=0}^{4}\frac{g_i}{i+1}\partial_{I_4}\frac{I_4^{i+1}}{B^{2i}_{\text{ref}}} +
		\frac{1}{2}\alpha_5\right)\BB \otimes \BB \\[1ex]
		&- \left( \frac{E_M}{2}\mu_{0}^{-1} +
		\sum_{i=0}^{4}\frac{g_i}{i+1}\partial_{I_4}\frac{I_4^{i+1}}{B^{2i}_{\text{ref}}}
		\right)\BB \cdot
		\BB \textbf{I} \\[1ex]
		& - \Big(\frac{1}{2}\alpha_5\left(\BB \cdot \Bve \cdot \BB \right)\textbf{I} -
		\frac{1}{2}\alpha_6\left(\BB\cdot \Bve^{2} \cdot \BB \right)\textbf{I}\\
		&+ \frac{1}{2}\alpha_5\BB\otimes\BB\cdot\Bve + \frac{1}{2}\alpha_6\BB\otimes\BB\cdot\Bve^{2}\\[1ex]
		&+ \frac{1}{2}\alpha_6
		\left[\BB \otimes \BB \cdot \Bve  + \Bve  \cdot \BB
		\otimes \BB \right]\Big). \\[1ex]
	\end{aligned}
\label{app_eq8}
\end{equation}

To approximate the \textit{homogenous response} of $g_i(I_i)$ in \req{app_eq1}, we further assume (A.8) is in $(i)$ pure magnetic loading, $(ii)$ volume-preserving deformation and $(iii)$ intact region so $g(d)=1$. As a result, we  first compute the trace of the total stress tensor as follows:
\eb
\begin{aligned}
	\tr(\Btau) & = 3\lambda I_1 +
	\frac{3}{2}\sum_{i=0}^{4}\frac{1}{i+1}\frac{\partial{g_i}}{\partial
		I_1}I_4^{i+1} + 2 \mu I_1 \\ & -
	\left(\frac{1}{2\mu_0}+2\sum_{i=0}^{4}\frac{g_i}{i+1}\frac{\partial
		I_4^{i+1}}{\partial I_4} - \frac{1}{2}\alpha_5 \right)I_4 - \alpha_5I_5 + \alpha_6 I_6.
\end{aligned}
\label{app_eq9}
\ee

Considering the conditions of volume-preserving deformation in pure
magnetic loading, we get the following relation
\begin{equation}
		\begin{cases}
	\frac{3}{2}\displaystyle\sum_{i=0}^{4}\frac{1}{i+1}\frac{\partial{g_i}}{\partial
		I_1}I_4^{i+1} - \left(\frac{1}{2\mu_0}+2\sum_{i=0}^{4}\frac{g_i}{i+1}\frac{\partial
		I_4^{i+1}}{\partial I_4} - \frac{1}{2}\alpha_5 \right)I_4=0, \\\\
	- \alpha_5I_5 + \alpha_6 I_6 =0.
	\end{cases}
\end{equation}\label{app_eqA100}
%
The functions $g_i(I_i)$ is derived from the first expression of
\req{app_eqA100} for each set of $i=\{0,1,2,3,4\}$ as follows:

\begin{equation}
	\begin{cases}
	\frac{\partial g_0}{\partial
		I_1}  = \frac{1}{3}\left(\frac{1}{\mu_0}-\alpha_5 \right) +
	\frac{4}{3}g_0\\\\
	\frac{\partial g_i}{\partial
		I_1} = \frac{3}{4}\left(i+1\right)g_i, \quad \text{where} \quad i = 1,2,3,4
	\end{cases}
\;,
\end{equation}\label{app_eqA101}
which finally leads to
\begin{equation}\label{app_eqApp102}
	\begin{cases}
	g_0  = \frac{3}{4}\alpha_0\exp\left(\frac{3}{4} I_1\right) -
	\frac{1}{3}\left(\frac{1}{\mu_0} - \alpha_5\right)\\\\
	g_i  = \frac{3(i+1)}{4}\alpha_i\exp\left(\frac{4(i+1)}{3}
	I_1\right),  & \text{where} \quad i=1,2,3,4.
	\end{cases}
\;,
\end{equation}
%
Here, in (A.16) $\alpha_i$ for $i=\{0,1,2,3,4,5\}$ refers to the material parameters depending on the ferromagnetic structure due to magnetisation effect, see \cite{fonteyn10}.
Alternatively, the five material constants $g_i$ in ${W}_{mag}$ could be seen as material constants which could be calibrated in experiments.

\end{Appendix}
 {\normalsize
 	\begin{spacing}{0.8}
 \bibliographystyle{ieeetr}
\bibliography{./ref}

\begin{thebibliography}{100}

\bibitem{john+etal07}
S.~John, J.~Sirohi, G.~Wang, and N.~M. Wereley, ``Comparison of piezoelectric,
  magnetostrictive, and electrostrictive hybrid hydraulic actuators,'' {\em
  Journal of intelligent material systems and structures}, vol.~18, no.~10,
  pp.~1035--1048, 2007.

\bibitem{bar-cohen+zhang08}
Y.~Bar-Cohen and Q.~Zhang, ``Electroactive polymer actuators and sensors,''
  {\em MRS bulletin}, vol.~33, no.~3, pp.~173--181, 2008.

\bibitem{kim+tadokoro07}
K.~J. Kim and S.~Tadokoro, ``Electroactive polymers for robotic applications,''
  {\em Artificial Muscles and Sensors}, vol.~23, p.~291, 2007.

\bibitem{berlincourt71}
D.~Berlincourt, {\em Ultrasonic transducer materials}.
\newblock Springer, 1971.

\bibitem{bhansali+vasudev12}
R.~Guldiken and O.~Onen, ``5 - mems ultrasonic transducers for biomedical
  applications,'' in {\em MEMS for Biomedical Applications} (S.~Bhansali and
  A.~Vasudev, eds.), Woodhead Publishing Series in Biomaterials, pp.~120--149,
  Woodhead Publishing, 2012.

\bibitem{damjanovic+newham92}
D.~Damjanovic and R.~Newnham, ``Electrostrictive and piezoelectric materials
  for actuator applications,'' {\em Journal of intelligent material systems and
  structures}, vol.~3, no.~2, pp.~190--208, 1992.

\bibitem{sekhar+etal21}
B.~C. Sekhar, B.~Dhanalakshmi, B.~S. Rao, S.~Ramesh, K.~V. Prasad, P.~S. Rao,
  and B.~P. Rao, ``Piezoelectricity and its applications,'' {\em
  Multifunctional Ferroelectric Materials}, p.~71, 2021.

\bibitem{joule1842}
J.~P. Joule, ``On a new class of magnetic forces,'' {\em Ann. Electr. Magn.
  Chem}, vol.~8, no.~1842, pp.~219--224, 1842.

\bibitem{villari1864}
E.~Villari, ``Intorno alle modificazioni del momento magnetico di una verga di
  ferro e di acciaio, prodotte per la trazione della medesima e pel passaggio
  di una corrente attraverso la stessa,'' {\em Il Nuovo Cimento (1855-1868)},
  vol.~20, no.~1, pp.~317--362, 1864.

\bibitem{birk+reichel+schroeder22}
C.~Birk, M.~Reichel, and J.~Schr\"oder, ``Magnetostatic simulations with
  consideration of exterior domains using the scaled boundary finite element
  method,'' {\em Computer Methods in Applied Mechanics and Engineering},
  vol.~399, p.~115362, 2022.

\bibitem{brigadnov+dorfmann03}
I.~Brigadnov and A.~Dorfmann, ``Mathematical modeling of magneto-sensitive
  elastomers,'' {\em International Journal of Solids and Structures}, vol.~40,
  no.~18, pp.~4659--4674, 2003.

\bibitem{dorfmann+ogden03}
A.~Dorfmann and R.~Ogden, ``Magnetoelastic modelling of elastomers,'' {\em
  European Journal of Mechanics - A/Solids}, vol.~22, no.~4, pp.~497--507,
  2003.

\bibitem{dorfmann+ogden+saccomandi04}
A.~Dorfmann, R.~Ogden, and G.~Saccomandi, ``Universal relations for non-linear
  magnetoelastic solids,'' {\em International Journal of Non-Linear Mechanics},
  vol.~39, no.~10, pp.~1699--1708, 2004.

\bibitem{thomas+triantafyllidis09}
J.~D. Thomas and N.~Triantafyllidis, ``On electromagnetic forming processes in
  finitely strained solids: Theory and examples,'' {\em Journal of the
  Mechanics and Physics of Solids}, vol.~57, no.~8, pp.~1391--1416, 2009.

\bibitem{belahcen+etal06}
A.~Belahcen, K.~Fonteyn, S.~Fortino, and R.~Kouhia, ``A coupled magnetoelastic
  model for ferromagnetic materials,'' {\em Proc. of the IX Finnish Mechanics
  Days. von Hertzen R., Halme T.(eds.)}, pp.~673--682, 2006.

\bibitem{fonteyn10}
K.~A. Fonteyn, {\em Energy-based magneto-mechanical model for electrical steel
  sheets}.
\newblock PhD thesis, Aalto-yliopiston {T}eknillinen {K}orkeakoulu, 2010.

\bibitem{fonteyn2010fem}
K.~Fonteyn, A.~Belahcen, R.~Kouhia, P.~Rasilo, and A.~Arkkio, ``Fem for
  directly coupled magneto-mechanical phenomena in electrical machines,'' {\em
  IEEE Transactions on Magnetics}, vol.~46, no.~8, pp.~2923--2926, 2010.

\bibitem{rasilo+etal19}
P.~Rasilo, D.~Singh, J.~Jeronen, U.~Aydin, F.~Martin, A.~Belahcen, L.~Daniel,
  and R.~Kouhia, ``Flexible identification procedure for thermodynamic
  constitutive models for magnetostrictive materials,'' {\em Proceedings of the
  Royal Society A}, vol.~475, no.~2223, p.~20180280, 2019.

\bibitem{miehe+rosato+kiefer11}
C.~Miehe, D.~Rosato, and B.~Kiefer, ``Variational principles in dissipative
  electro-magneto-mechanics: a framework for the macro-modeling of functional
  materials,'' {\em International Journal for Numerical Methods in
  Engineering}, vol.~86, no.~10, pp.~1225--1276, 2011.

\bibitem{miehe2012geometrically}
C.~Miehe and G.~Ethiraj, ``A geometrically consistent incremental variational
  formulation for phase field models in micromagnetics,'' {\em Computer methods
  in applied mechanics and engineering}, vol.~245, pp.~331--347, 2012.

\bibitem{ethiraj14}
G.~Ethiraj, {\em Computational modeling of ferromagnetics and
  magnetorheological elastomers}.
\newblock 2014.

\bibitem{ethiraj+miehe16}
G.~Ethiraj and C.~Miehe, ``Multiplicative magneto-elasticity of
  magnetosensitive polymers incorporating micromechanically-based network
  kernels,'' {\em International Journal of Engineering Science}, vol.~102,
  pp.~93--119, 2016.

\bibitem{hanappier+charkaluk+triantafyllidis21}
N.~Hanappier, E.~Charkaluk, and N.~Triantafyllidis, ``A coupled
  electromagnetic-thermomechanical approach for the modeling of electric
  motors,'' {\em Journal of the Mechanics and Physics of Solids}, vol.~149,
  p.~104315, 2021.

\bibitem{zhang+etal21}
B.~Zhang, Y.~Kou, K.~Jin, and X.~Zheng, ``A multi-field coupling model for the
  magnetic-thermal-structural analysis in the electromagnetic rail launch,''
  {\em Journal of Magnetism and Magnetic Materials}, vol.~519, p.~167495, 2021.

\bibitem{ma+etal16}
Z.~Ma, H.~Zhao, W.~Liu, and L.~Ren, ``Thermo-mechanical coupled in situ fatigue
  device driven by piezoelectric actuator,'' {\em Precision Engineering},
  vol.~46, pp.~349--359, 2016.

\bibitem{zhou+etal21}
L.~Zhou, J.~Tang, W.~Tian, B.~Xue, and X.~Li, ``A multi-physics coupling
  cell-based smoothed finite element micromechanical model for the transient
  response of magneto-electro-elastic structures with the asymptotic
  homogenization method,'' {\em Thin-Walled Structures}, vol.~165, p.~107991,
  2021.

\bibitem{francfort+marigo98}
G.~A. Francfort and J.-J. Marigo, ``Revisiting brittle fracture as an energy
  minimization problem,'' {\em Journal of the Mechanics and Physics of Solids},
  vol.~46, no.~8, pp.~1319--1342, 1998.

\bibitem{bourdin+francfort+marigo08}
B.~Bourdin, G.~A. Francfort, and J.-J. Marigo, ``The variational approach to
  fracture,'' {\em Journal of Elasticity}, vol.~91, no.~1, pp.~5--148, 2008.

\bibitem{dal+toador02}
G.~Dal~Maso and R.~Toader, ``A model for the quasi-static growth of brittle
  fractures: Existence and approximation results,'' {\em Archive for Rational
  Mechanics and Analysis}, vol.~162, no.~2, pp.~101--135, 2002.

\bibitem{mumford+shah89}
D.~B. Mumford and J.~Shah, ``Optimal approximations by piecewise smooth
  functions and associated variational problems,'' {\em Communications on Pure
  and Applied Mathematics}, 1989.

\bibitem{hakim+karma09}
V.~Hakim and A.~Karma, ``Laws of crack motion and phase-field models of
  fracture,'' {\em Journal of the Mechanics and Physics of Solids}, vol.~57,
  no.~2, pp.~342--368, 2009.

\bibitem{miehe+hofacker+welschinger10}
C.~Miehe, M.~Hofacker, and F.~Welschinger, ``A phase field model for
  rate-independent crack propagation: Robust algorithmic implementation based
  on operator splits,'' {\em Computer Methods in Applied Mechanics and
  Engineering}, vol.~199, no.~45, pp.~2765--2778, 2010.

\bibitem{miehe+welschinger+hofacker101}
C.~Miehe, F.~Welschinger, and M.~Hofacker, ``A phase field model of
  electromechanical fracture,'' {\em Journal of the Mechanics and Physics of
  Solids}, vol.~58, no.~10, pp.~1716--1740, 2010.

\bibitem{linder+miehe12}
C.~Linder and C.~Miehe, ``Effect of electric displacement saturation on the
  hysteretic behavior of ferroelectric ceramics and the initiation and
  propagation of cracks in piezoelectric ceramics,'' {\em Journal of the
  Mechanics and Physics of Solids}, vol.~60, no.~5, pp.~882--903, 2012.

\bibitem{monk2003finite}
P.~Monk {\em et~al.}, {\em Finite element methods for Maxwell's equations}.
\newblock Oxford University Press, 2003.

\bibitem{gross2004electromagnetic}
P.~W. Gross, P.~W. Gross, P.~R. Kotiuga, and R.~P. Kotiuga, {\em
  Electromagnetic theory and computation: a topological approach}, vol.~48.
\newblock Cambridge University Press, 2004.

\bibitem{noii2021quasi}
N.~Noii, M.~Fan, T.~Wick, and Y.~Jin, ``A quasi-monolithic phase-field
  description for orthotropic anisotropic fracture with adaptive mesh
  refinement and primal--dual active set method,'' {\em Engineering Fracture
  Mechanics}, vol.~258, p.~108060, 2021.

\bibitem{noii2020adaptive}
N.~Noii, F.~Aldakheel, T.~Wick, and P.~Wriggers, ``An adaptive global--local
  approach for phase-field modeling of anisotropic brittle fracture,'' {\em
  Computer Methods in Applied Mechanics and Engineering}, vol.~361, p.~112744,
  2020.

\bibitem{shipman12}
J.~Shipman, J.~D. Wilson, and C.~A. Higgins, {\em An introduction to physical
  science}.
\newblock Cengage Learning, 2012.

\bibitem{emery+camps17}
W.~Emery and A.~Camps, ``Chapter 2 - basic electromagnetic concepts and
  applications to optical sensors,'' in {\em Introduction to Satellite Remote
  Sensing} (W.~Emery and A.~Camps, eds.), pp.~43--83, Elsevier, 2017.

\bibitem{spaldin2010magnetic}
N.~A. Spaldin, {\em Magnetic materials: fundamentals and applications}.
\newblock Cambridge university press, 2010.

\bibitem{fliegans2020modeling}
J.~Fliegans, O.~Tosoni, N.~Dempsey, and G.~Delette, ``Modeling of
  demagnetization processes in permanent magnets measured in closed-circuit
  geometry,'' {\em Applied Physics Letters}, vol.~116, no.~6, p.~062405, 2020.

\bibitem{bruzzese2022theory}
C.~Bruzzese, {\em Theory of Electrical Machines}.
\newblock Societ{\`a} Editrice Esculapio, 2022.

\bibitem{jackson1999classical}
J.~D. Jackson, ``Classical electrodynamics,'' 1999.

\bibitem{melia2001electrodynamics}
F.~Melia, {\em Electrodynamics}.
\newblock Chicago Lectures in Physics, University of Chicago Press, 2001.

\bibitem{van2021comment}
H.~van Hees, ``Comment on 'defining the electromagnetic potentials','' {\em
  European Journal of Physics}, vol.~42, no.~2, p.~028003, 2021.

\bibitem{maudlin2018ontological}
T.~Maudlin, ``Ontological clarity via canonical presentation: Electromagnetism
  and the aharonov--bohm effect,'' {\em Entropy}, vol.~20, no.~6, p.~465, 2018.

\bibitem{powell2008two}
J.~D. Powell and A.~E. Zielinski, ``Two-dimensional current diffusion in the
  rails of a railgun,'' tech. rep., Army Research LAB ABERDEEN Proving Ground
  MD Weapons and Materials Research, 2008.

\bibitem{zhang2021multi}
B.~Zhang, Y.~Kou, K.~Jin, and X.~Zheng, ``A multi-field coupling model for the
  magnetic-thermal-structural analysis in the electromagnetic rail launch,''
  {\em Journal of Magnetism and Magnetic Materials}, vol.~519, p.~167495, 2021.

\bibitem{bourdin2000numerical}
B.~Bourdin, G.~Francfort, and J.-J. Marigo, ``Numerical experiments in
  revisited brittle fracture,'' {\em Journal of the Mechanics and Physics of
  Solids}, vol.~48, no.~4, pp.~797--826, 2000.

\bibitem{francfort1998revisiting}
G.~Francfort and J.-J. Marigo, ``Revisiting brittle fracture as an energy
  minimization problem,'' {\em Journal of the Mechanics and Physics of Solids},
  vol.~46, no.~8, pp.~1319--1342, 1998.

\bibitem{kienle2019}
D.~Kienle, F.~Aldakheel, and M.-A. Keip, ``A finite-strain phase-field approach
  to ductile failure of frictional materials,'' {\em International Journal of
  Solids and Structures}, vol.~172, pp.~147--162, 2019.

\bibitem{dittmann2020}
M.~Dittmann, F.~Aldakheel, J.~Schulte, F.~Schmidt, M.~Kr{\"u}ger, P.~Wriggers,
  and C.~Hesch, ``Phase-field modeling of porous-ductile fracture in non-linear
  thermo-elasto-plastic solids,'' {\em Computer Methods in Applied Mechanics
  and Engineering}, vol.~361, p.~112730, 2020.

\bibitem{ruan2022thermo}
H.~Ruan, S.~Rezaei, Y.~Yang, D.~Gross, and B.-X. Xu, ``A thermo-mechanical
  phase-field fracture model: Application to hot cracking simulations in
  additive manufacturing,'' {\em Journal of the Mechanics and Physics of
  Solids}, p.~105169, 2022.

\bibitem{peng2023meso}
Z.~Peng, Q.~Wang, W.~Zhou, X.~Chang, Q.~Yue, and C.~Huang, ``Meso-scale
  simulation of thermal fracture in concrete based on the coupled
  thermal--mechanical phase-field model,'' {\em Construction and Building
  Materials}, vol.~403, p.~133095, 2023.

\bibitem{ambati+kruse+lorenzis16}
M.~Ambati, R.~Kruse, and L.~De~Lorenzis, ``A phase-field model for ductile
  fracture at finite strains and its experimental verification,'' {\em
  Computational Mechanics}, vol.~57, pp.~149--167, 2016.

\bibitem{yin2020ductile}
B.~Yin and M.~Kaliske, ``A ductile phase-field model based on degrading the
  fracture toughness: Theory and implementation at small strain,'' {\em
  Computer Methods in Applied Mechanics and Engineering}, vol.~366, p.~113068,
  2020.

\bibitem{noii2021bayesian1}
N.~Noii, A.~Khodadadian, J.~Ulloa, F.~Aldakheel, T.~Wick, S.~Fran{\c{c}}ois,
  and P.~Wriggers, ``Bayesian inversion for unified ductile phase-field
  fracture,'' {\em Computational Mechanics}, vol.~68, no.~4, pp.~943--980,
  2021.

\bibitem{heider2020phase}
Y.~Heider and W.~Sun, ``A phase field framework for capillary-induced fracture
  in unsaturated porous media: Drying-induced vs. hydraulic cracking,'' {\em
  Computer Methods in Applied Mechanics and Engineering}, vol.~359, p.~112647,
  2020.

\bibitem{ulloa2022variational}
J.~Ulloa, N.~Noii, R.~Alessi, F.~Aldakheel, G.~Degrande, and S.~Fran{\c{c}}ois,
  ``Variational modeling of hydromechanical fracture in saturated porous media:
  A micromechanics-based phase-field approach,'' {\em Computer Methods in
  Applied Mechanics and Engineering}, vol.~396, p.~115084, 2022.

\bibitem{noii2019phase}
N.~Noii and T.~Wick, ``A phase-field description for pressurized and
  non-isothermal propagating fractures,'' {\em Computer Methods in Applied
  Mechanics and Engineering}, vol.~351, pp.~860--890, 2019.

\bibitem{noii2022bayesian1}
N.~Noii, A.~Khodadadian, J.~Ulloa, F.~Aldakheel, T.~Wick, S.~Fran{\c{c}}ois,
  and P.~Wriggers, ``Bayesian inversion with open-source codes for various
  one-dimensional model problems in computational mechanics,'' {\em Archives of
  Computational Methods in Engineering}, pp.~1--34, 2022.

\bibitem{noii2022probabilistic}
N.~Noii, A.~Khodadadian, and F.~Aldakheel, ``Probabilistic failure mechanisms
  via monte carlo simulations of complex microstructures,'' {\em Computer
  Methods in Applied Mechanics and Engineering}, vol.~399, p.~115358, 2022.

\bibitem{noii2021bayesian}
N.~Noii, A.~Khodadadian, and T.~Wick, ``Bayesian inversion for anisotropic
  hydraulic phase-field fracture,'' {\em Computer Methods in Applied Mechanics
  and Engineering}, vol.~386, p.~114118, 2021.

\bibitem{noii2022bayesian}
N.~Noii, A.~Khodadadian, and T.~Wick, ``Bayesian inversion using global-local
  forward models applied to fracture propagation in porous media,'' {\em
  International Journal for Multiscale Computational Engineering}, vol.~20,
  no.~3, 2022.

\bibitem{noii2023level}
N.~Noii, H.~A. Jahangiry, and H.~Waisman, ``Level-set topology optimization for
  ductile and brittle fracture resistance using the phase-field method,'' {\em
  Computer Methods in Applied Mechanics and Engineering}, vol.~409, p.~115963,
  2023.

\bibitem{liu2022phase}
Z.~Liu, J.~Reinoso, and M.~Paggi, ``Phase field modeling of brittle fracture in
  large-deformation solid shells with the efficient quasi-newton solution and
  global--local approach,'' {\em Computer Methods in Applied Mechanics and
  Engineering}, vol.~399, p.~115410, 2022.

\bibitem{aldakheel2020global}
F.~Aldakheel, N.~Noii, T.~Wick, and P.~Wriggers, ``A global--local approach for
  hydraulic phase-field fracture in poroelastic media,'' {\em Computers and
  Mathematics with Applications}, 2020.

\bibitem{aldakheel2021multilevel}
F.~Aldakheel, N.~Noii, T.~Wick, O.~Allix, and P.~Wriggers, ``Multilevel
  global--local techniques for adaptive ductile phase-field fracture,'' {\em
  Computer Methods in Applied Mechanics and Engineering}, vol.~387, p.~114175,
  2021.

\bibitem{hageman2022electro}
T.~Hageman and E.~Mart{\'\i}nez-Pa{\~n}eda, ``An electro-chemo-mechanical
  framework for predicting hydrogen uptake in metals due to aqueous
  electrolytes,'' {\em Corrosion Science}, vol.~208, p.~110681, 2022.

\bibitem{wu2022crack}
J.-Y. Wu and Y.-F. Hong, ``Crack nucleation and propagation of
  electromagneto-thermo-mechanical fracture in bulk superconductors during
  magnetization,'' {\em Journal of the Mechanics and Physics of Solids},
  p.~105168, 2022.

\bibitem{zhao2022phase}
Y.~Zhao, R.~Wang, and E.~Mart{\'\i}nez-Pa{\~n}eda, ``A phase field
  electro-chemo-mechanical formulation for predicting void evolution at the
  li--electrolyte interface in all-solid-state batteries,'' {\em Journal of the
  Mechanics and Physics of Solids}, vol.~167, p.~104999, 2022.

\bibitem{kuhn2015}
C.~Kuhn, A.~Schl{\"u}ter, and R.~M{\"u}ller, ``On degradation functions in
  phase field fracture models,'' {\em Computational Materials Science},
  vol.~108, pp.~374--384, 2015.

\bibitem{wu2017}
J.-Y. Wu, ``A unified phase-field theory for the mechanics of damage and
  quasi-brittle failure,'' {\em Journal of the Mechanics and Physics of
  Solids}, vol.~103, pp.~72--99, 2017.

\bibitem{wu2018}
J.-Y. Wu, V.~Nguyen, C.~Nguyen, D.~Sutula, S.~Bordas, and S.~Sinaie, ``Phase
  field modeling of fracture,'' {\em Advances in Applied Mechancis: Multi-Scale
  Theory and Computation}, vol.~52, 2018.

\bibitem{miehe2010phase}
C.~Miehe, M.~Hofacker, and F.~Welschinger, ``A phase field model for
  rate-independent crack propagation: Robust algorithmic implementation based
  on operator splits,'' {\em Computer Methods in Applied Mechanics and
  Engineering}, vol.~199, no.~45-48, pp.~2765--2778, 2010.

\bibitem{miehe2015phase}
C.~Miehe, L.-M. Schaenzel, and H.~Ulmer, ``Phase field modeling of fracture in
  multi-physics problems. part i. balance of crack surface and failure criteria
  for brittle crack propagation in thermo-elastic solids,'' {\em Computer
  Methods in Applied Mechanics and Engineering}, vol.~294, pp.~449--485, 2015.

\bibitem{lee2016pressure}
S.~Lee, M.~F. Wheeler, and T.~Wick, ``Pressure and fluid-driven fracture
  propagation in porous media using an adaptive finite element phase field
  model,'' {\em Computer Methods in Applied Mechanics and Engineering},
  vol.~305, pp.~111--132, 2016.

\bibitem{miehe2016phase}
C.~Miehe and S.~Mauthe, ``Phase field modeling of fracture in multi-physics
  problems. part iii. crack driving forces in hydro-poro-elasticity and
  hydraulic fracturing of fluid-saturated porous media,'' {\em Computer Methods
  in Applied Mechanics and Engineering}, vol.~304, pp.~619--655, 2016.

\bibitem{pao+hutter75}
Y.-H. Pao and K.~Hutter, ``Electrodynamics for moving elastic solids and
  viscous fluids,'' {\em Proceedings of the IEEE}, vol.~63, no.~7,
  pp.~1011--1021, 1975.

\bibitem{daniel+bernard+hubert20}
L.~Daniel, L.~Bernard, and O.~Hubert, ``Multiscale modeling of magnetic
  materials,'' 2020.

\bibitem{flatau+dapino+calkins00}
A.~Flatau, M.~Dapino, and F.~Calkins, ``5.26 - magnetostrictive composites,''
  in {\em Comprehensive Composite Materials} (A.~Kelly and C.~Zweben, eds.),
  pp.~563--574, Oxford: Pergamon, 2000.

\bibitem{gao+etal22}
C.~Gao, Z.~Zeng, S.~Peng, and C.~Shuai, ``Magnetostrictive alloys: Promising
  materials for biomedical applications,'' {\em Bioactive Materials}, vol.~8,
  pp.~177--195, 2022.

\bibitem{dapino04}
M.~J. Dapino, ``On magnetostrictive materials and their use in adaptive
  structures,'' {\em Structural Engineering and Mechanics}, vol.~17, no.~3-4,
  pp.~303--330, 2004.

\bibitem{khodadadian2020bayesian}
A.~Khodadadian, N.~Noii, M.~Parvizi, M.~Abbaszadeh, T.~Wick, and C.~Heitzinger,
  ``A bayesian estimation method for variational phase-field fracture
  problems,'' {\em Computational Mechanics}, vol.~66, no.~4, pp.~827--849,
  2020.

\bibitem{wang2020phase}
Q.~Wang, Y.~Feng, W.~Zhou, Y.~Cheng, and G.~Ma, ``A phase-field model for
  mixed-mode fracture based on a unified tensile fracture criterion,'' {\em
  Computer Methods in Applied Mechanics and Engineering}, vol.~370, p.~113270,
  2020.

\bibitem{wang2023modeling}
Q.~Wang, Q.~Yue, W.~Zhou, Y.~Feng, and X.~Chang, ``Modeling of both
  tensional-shear and compressive-shear fractures by a unified phase-field
  model,'' {\em Applied Mathematical Modelling}, vol.~117, pp.~162--196, 2023.

\bibitem{kovetz00}
A.~Kovetz, {\em Electromagnetic theory}, vol.~975.
\newblock Oxford University Press Oxford, 2000.

\bibitem{wu2021wrinkling}
B.~Wu and M.~Destrade, ``Wrinkling of soft magneto-active plates,'' {\em
  International Journal of Solids and Structures}, vol.~208, pp.~13--30, 2021.

\bibitem{trimarco2002stresses}
C.~Trimarco, ``Stresses and momenta in electromagnetic materials,'' {\em
  Mechanics Research Communications}, vol.~29, no.~6, pp.~485--492, 2002.

\bibitem{zheng2011magnetic}
X.~Zheng and K.~Jin, ``Magnetic force models for magnetizable elastic bodies in
  the magnetic field,'' {\em From Waves in Complex Systems to Dynamics of
  Generalized Continua: Tributes to Professor Yih-Hsing Pao on His 80th
  Birthday}, pp.~353--383, 2011.

\bibitem{henjes1994traction}
K.~Henjes, ``The traction force in magnetic separators,'' {\em Measurement
  Science and Technology}, vol.~5, no.~9, p.~1105, 1994.

\bibitem{rinaldi2002body}
C.~Rinaldi and H.~Brenner, ``Body versus surface forces in continuum mechanics:
  Is the maxwell stress tensor a physically objective cauchy stress?,'' {\em
  Physical Review E}, vol.~65, no.~3, p.~036615, 2002.

\bibitem{mauthe2017variational}
S.~A. Mauthe, {\em Variational multiphysics modeling of diffusion in elastic
  solids and hydraulic fracturing in porous media}.
\newblock Stuttgart: Institut fur Mechanik (Bauwesen), 2017.

\bibitem{miehe2015minimization}
C.~Miehe, S.~Mauthe, and S.~Teichtmeister, ``Minimization principles for the
  coupled problem of darcy--biot-type fluid transport in porous media linked to
  phase field modeling of fracture,'' {\em Journal of the Mechanics and Physics
  of Solids}, vol.~82, pp.~186--217, 2015.

\bibitem{li2019computational}
J.~Li and Y.-T. Chen, {\em Computational partial differential equations using
  MATLAB{\textregistered}}.
\newblock Crc Press, 2019.

\bibitem{garcia02}
E.~Garcia, {\em Solution to the instationary Maxwell equations with charges in
  non-convex domains}.
\newblock PhD thesis, PhD thesis, Universit{\'e} Paris VI, France, 2002.

\bibitem{jamelot05}
E.~Jamelot, ``Solution to maxwell equations with continuous galerkin finite
  elements,'' {\em Ph. D. Thesis, Ecole Polytechnique, Palaiseau, France},
  2005.

\bibitem{ciarlet+jamelot07}
P.~Ciarlet and E.~Jamelot, ``Continuous galerkin methods for solving the
  time-dependent maxwell equations in 3d geometries,'' {\em Journal of
  Computational Physics}, vol.~226, no.~1, pp.~1122--1135, 2007.

\bibitem{asadzadeh+beilina23}
M.~Asadzadeh and L.~Beilina, ``A stabilized p1 domain decomposition finite
  element method for time harmonic maxwell's equations,'' {\em Mathematics and
  Computers in Simulation}, vol.~204, pp.~556--574, 2023.

\bibitem{meunier2010finite}
G.~Meunier, ``The finite element method for electromagnetic modeling,'' 2010.

\bibitem{cardoso2016electromagnetics}
J.~R. Cardoso, {\em Electromagnetics through the finite element method: A
  simplified approach using Maxwell's equations}.
\newblock Crc Press, 2016.

\bibitem{bastos2003electromagnetic}
J.~P.~A. Bastos and N.~Sadowski, {\em Electromagnetic modeling by finite
  element methods}.
\newblock CRC press, 2003.

\bibitem{alnaes+etal15}
M.~Aln{\ae}s, J.~Blechta, J.~Hake, A.~Johansson, B.~Kehlet, A.~Logg,
  C.~Richardson, J.~Ring, M.~E. Rognes, and G.~N. Wells, ``The fenics project
  version 1.5,'' {\em Archive of Numerical Software}, vol.~3, no.~100, 2015.

\bibitem{logg+etal12}
G.~N. W. e.~a. A.~Logg, {K.-A.}~Mardal, {\em Automated Solution of Differential
  Equations by the Finite Element Method}.
\newblock Springer, 2012.

\bibitem{yue2023efficient}
Q.~Yue, Q.~Wang, W.~Zhou, T.~Rabczuk, X.~Zhuang, B.~Liu, and X.~Chang, ``An
  efficient adaptive length scale insensitive phase-field model for
  three-dimensional fracture of solids using trilinear multi-node elements,''
  {\em International Journal of Mechanical Sciences}, vol.~253, p.~108351,
  2023.

\end{thebibliography}
\end{spacing}}
\end{document}